\documentclass[]{nmeauth}

\usepackage{moreverb}
\usepackage[colorlinks,bookmarksopen,bookmarksnumbered,citecolor=red,urlcolor=red]{hyperref}
\newcommand\BibTeX{{\rmfamily B\kern-.05em \textsc{i\kern-.025em b}\kern-.08em
T\kern-.1667em\lower.7ex\hbox{E}\kern-.125emX}}
\def\volumeyear{2016}

\usepackage{lineno,hyperref}
\modulolinenumbers[5]
\usepackage{url}
\usepackage{xstring}
\usepackage{amsmath}
\usepackage{mathtools}
\usepackage{amssymb}
\usepackage{color}
\usepackage{bigints}		
\usepackage{empheq}
\usepackage{cases}		
\usepackage{booktabs}
\usepackage{multirow}
\usepackage{footnote}
\usepackage{accents}
\usepackage{arcs}

\usepackage[boxed, ruled,vlined, lined, linesnumbered]{algorithm2e}
 
\usepackage{float}
\usepackage{relsize}
\usepackage{graphicx}
\usepackage{caption}
\usepackage{subcaption}
\usepackage{etoolbox}
\usepackage{xfrac}

\usepackage{xcolor}
\definecolor{bl}{rgb}{0.0,0.0,0.0} 
\definecolor{crimson}{rgb}{0.7,0.0,0.2}
\definecolor{darkGrey}{RGB}{98, 98, 98}
\definecolor{customBlue}{RGB}{0, 0, 204}

\newcommand{\boldZero}{\textbf {0}}

\newcommand{\RR}[1]{\mathbb{R}^{#1}}
\newcommand{\innerProduct}[2]{\langle #1, #2 \rangle}
\newcommand{\norm}[2]{{\lVert #1\rVert}_{#2}}
\newcommand{\defeq}{\vcentcolon=}
\newcommand{\orderOf}{\mathcal{O}}
\newcommand{\varWithSubSup}[3]{{#1}^{#2}_{#3} }

\newcommand{\numParam}{N_{\mu}}
\newcommand{\paramSpace}{\RR{\numParam}}

\newcommand{\sampleOutputState}{\mathbf{P}}
\newcommand{\sampleRomOutputState}{\widetilde{\sampleOutputState}}

\newcommand{\numHFSsymbol}{N_\text{train}}
\newcommand{\numHFS}[1]{\numHFSsymbol={#1}}

\newcommand{\numTestSamples}{N_\text{BHP}}

\newcommand{\x}{\textbf {x}}

\newcommand{\z}{\textbf {z}}

\newcommand{\nstateRom}{\ell}

\newcommand{{\perm}}{\textbf{k}}
\newcommand{\phaseVelocity}[1]{ \textbf{v}_{#1} }
\newcommand{\phaseDensity}[1]{ \rho_{#1} }
\newcommand{\phaseFlowRate}[1]{ q_{#1} }
\newcommand{\relperm}{k_{rj} }
\newcommand{\phaseFlowRateWellD}{(\phaseFlowRate{j})_{d} }

\newcommand{\nCells} {\ensuremath{N_{{c}}}}
\newcommand{\nSampleOutputState} {N_{P}}

\newcommand{\BHP} {\ensuremath{\textbf{u}}}

\newcommand{\BHPnoBf} {\ensuremath{u}}
\newcommand{\nWells} {\ensuremath{N_{u}}}

\newcommand{\wellName}[2] {\ensuremath {#1}_{#2}}

\newcommand{\surr} {\rm{surr}}
\newcommand{\corr} {\rm{corr}}

\newcommand{\dt} {\Delta t}

\newcommand{\timestep}{{n}}

\newcommand{\PODbasis}{\mathbf{\Phi}}

\newcommand{\LHSbasis}{\mathbf{\Psi}}
\newcommand{\PODMean}{\bar {\x}}

\newcommand{\LHSbasisIPlusOne} { \resTwo{ \LHSbasis} {i} }
\newcommand{\LHSbasisIPlusOneT} { \left( {\resTwo{ \LHSbasis} {i} } \right)^{T}}

\newcommand{\resTwo}[2] { {\bf {#1}}^{ #2} }

\newcommand{\JIPlusOne} { \resTwo{\acute J} {i} }
\newcommand{\BIPlusOne} { \resTwo{\acute B} {i} }
\newcommand{\CIPlusOne} { \resTwo{\acute C} {i} }
\newcommand{\g} {{\mathbf g}}
\newcommand{\gArgOne} {{\bar \x^1}}
\newcommand{\gArgTwo} {{\bar \x^2}}
\newcommand{\gArgThree} {{\bar\uVar}}

\newcommand{\xTraining} { \acute{ \x} }

\newcommand{\zIPlusOne} { \resTwo{\acute z} {i} }
\newcommand{\uTraining} {  \acute{\mathbf u}  }
\newcommand{\uIPlusOne} { \resTwo{\acute u} {i} }
\newcommand{\JrIPlusOne} {\resTwo{\acute  J}{i}_{r}}
\newcommand{\BrIPlusOne} {\resTwo{\acute B}{i}_{r}}
\newcommand{\CrIPlusOne} {\resTwo{\acute C}{i}_{r}}
\newcommand{\JI} { \resTwo {\acute J} {i} }
\newcommand{\BI} { \resTwo {\acute B} {i} }
\newcommand{\CI} { \resTwo {\acute C} {i} }

\newcommand{\xI} { \resTwo {{\xTraining}} {i} }
\newcommand{\xIGen} {\xTraining}
\newcommand{\zI} { \resTwo {{\acute z}} {i} }
\newcommand{\zIGen} {{\bf{\acute z}}}

\newcommand{\xIMinusOne} { \resTwo {{\xTraining}} {i-\rm{1}} }
\newcommand{\zIMinusOne} { \resTwo {{\acute z}} {i-\rm{1}} }

\newcommand{\gNPlusOne} { \resTwo{g} {n} }
\newcommand{\gNPlusOneTilde} { \resTwo{\tilde{g}} {n} }

\newcommand{\uVar} { {\bf{u}}}
\newcommand{\uNPlusOne} { \resTwo{u} {n}  }

\newcommand{\xNMinusOne} { \resTwo{x} {n \rm{-1}} }
\newcommand{\zNMinusOne} { \resTwo{z} {n \rm{-1}} }

\newcommand{\xN} { \resTwo {x} {n} }
\newcommand{\zN} { \resTwo {z} {n} }

\newcommand{\uArg}[1] { \BHP^{#1} }
\newcommand{\uLinArg}[1] { \acute{\textbf{u}}^{#1} }

\newcommand{\RL} {\rm{RL}}
\newcommand{\SDrl} {(S_d)_{\RL}}
\newcommand{\gRL} { \resTwo{g}{n}_{\RL} }

\newcommand{\perturbation}[1]{
\sum_{d\in\wellBlockSetProducer}\frac{\sum_{k=1}^{\numTimesteps}
|u_d^k- \acute u_d^k|\dt^k
}{\sum_{k=1}^{\numTimesteps}|\acute
u_d^k|\dt^k}
}

\newcommand{\perturbationProducer}{\perturbation{n_{p_w}}}

\newcommand{\letterForMapping}{r}

\newcommand{\romOperatorGenerateFeatures}{\features}
\newcommand{\romOperatorGenerateFeaturesArg}[1]{\romOperatorGenerateFeatures^{#1}}

\newcommand{\relerrorModelQoI}{\rel{\letterForMapping}_{q}}
\newcommand{\errorModelQoI}{\letterForMapping_{q}}
\newcommand{\errorModelPrimaryVar}[1]{\letterForMapping_{ss,#1}}
\newcommand{\relerrorModelPrimaryVar}[1]{\rel{\letterForMapping}_{ss,#1}}

\newcommand{\errorModelQoIApprox}{\hat{\letterForMapping}_{q}}
\newcommand{\relerrorModelQoIApprox}{\hat{\rel{\letterForMapping}}_{q}}
\newcommand{\errorModelPrimaryVarApprox}[1]{\hat{\letterForMapping}_{ss,#1}}

\newcommand{\errorGen}{ \delta}
\newcommand{\relerrorGen}{\rel\delta}

\newcommand{\errorOperatorQoI}{ \errorGen_{\fomOperatorTwo}}

\newcommand{\errorOperatorPrimaryVar}{ \boldsymbol\errorGen_{ss}}
\newcommand{\relerrorOperatorPrimaryVar}{ \rel{\boldsymbol\delta}_{ss}}

\newcommand{\errorOperatorQoIApprox}{ \hat{\delta}_{\fomOperatorTwo}}
\newcommand{\errorOperatorPrimaryVarApprox}{ \hat{\boldsymbol\delta}_{ss}}

\newcommand{\rel}[1]{\bar{#1}}

\newcommand{\errorOperatorQoIArg}[1]{\errorOperatorQoI^{#1}}
\newcommand{\errorOperatorQoITimeN}{\errorOperatorQoIArg{\timestep}}

\newcommand{\relerrorOperatorQoI}{ \relerrorGen_{\fomOperatorTwo}}
\newcommand{\relerrorOperatorQoIArg}[1]{ \relerrorOperatorQoI^{#1}}

\newcommand{\relerrorOperatorQoITimeN}{\relerrorOperatorQoIArg{\timestep}}
\newcommand{\relerrorOperatorQoIApproxTimeN}{\hat{\rel{\errorGen}}_{\fomOperatorTwo}^{\timestep}}

\newcommand{\relerrorOperatorPrimaryVarArg}[1]{\relerrorOperatorPrimaryVar^{#1}}
\newcommand{\relerrorOperatorPrimaryVarTimeN}{\relerrorOperatorPrimaryVarArg{\timestep}}
\newcommand{\relerrorOperatorPrimaryVarApproxTimeN}{\hat{\rel{\boldsymbol\errorGen}}_{\fomOperatorOne}^{\timestep}}

\newcommand{\errorOperatorPrimaryVarArg}[1]{\errorOperatorPrimaryVar^{#1}}
\newcommand{\errorOperatorPrimaryVarTimeN}{\errorOperatorPrimaryVarArg{\timestep}}

\newcommand{\errorOperatorQoIApproxTimeN}{\errorOperatorQoIApprox^{\timestep}}
\newcommand{\errorOperatorPrimaryVarApproxTimeN}{\errorOperatorPrimaryVarApprox^{\timestep}}

\newcommand{\featureSymbol}{{f}}
\newcommand{\errorFun}{{h}}

\newcommand{\features}{\boldsymbol{\featureSymbol}}
\newcommand{\featuresArg}[1]{\features^{#1}}
\newcommand{\featuresN}{{\featuresArg{\timestep}}}

\newcommand{\numFeatures}{{N_{f}}}
\newcommand{\featuresSpace}{\RR{1\times\numFeatures}}

\newcommand{\featuresClassify}{\features_{\rm{c}}}
\newcommand{\featuresClassifyArg}[1]{{\featuresClassify^{#1}}}
\newcommand{\featuresClassifyN}{{\featuresClassifyArg{\timestep}}}

\newcommand{\numFeaturesClassify}{N_{f_c}}

\newcommand{\errorMapNoise}{\varepsilon}

\newcommand{\howMuchMemory}[1]{\tau = {#1}}
\newcommand{\memory}{\tau}
\newcommand{\featuresMem}{\features_{\rm{mem}}}
\newcommand{\featuresMemArg}[1]{{\featuresMem^{#1}}}
\newcommand{\featuresMemN}{{\featuresMemArg{\timestep}}}
\newcommand{\numFeaturesMemory}{ \left( 1+\memory \right) \numFeatures }

\newcommand{\dSdTn}{{ (\unitVec{}_{2d}^T \PODbasis)(\zN-\zNMinusOne)}~/~{\dt^{\timestep}}}
\newcommand{\dSdTi}{{ (\unitVec{}_{2d}^T \PODbasis) (\zI-\zIMinusOne)}~/~{\dt^{i}}}

\newcommand{\paramsGen}{\boldsymbol\mu}
\newcommand{\paramsGenArg}[1]{\paramsGen^{#1}}
\newcommand{\paramsLinearize}{\acute\paramsGen}
\newcommand{\paramsTrain}{\boldsymbol\mu_\text{train}}
\newcommand{\paramsTrainArg}[1]{\paramsTrain^{#1}}
\newcommand{\paramsTest}{\boldsymbol\mu_\text{test}}
\newcommand{\paramsTestArg}[1]{\paramsTest^{#1}}
\newcommand{\nstateFom}{N_{x}}
\newcommand{\numTimesteps}{N_{t}}
\newcommand{\paramTrain}{\mathcal T_\text{EMML}}
\newcommand{\paramTest}{\mathcal S_\text{EMML}}
\newcommand{\paramTrainROMES}{\mathcal T_\text{ROMES}}
\newcommand{\paramTrainTPWL}{\mathcal T_\text{TPWL}}

\newcommand{\wellBlockSet}{\mathcal D}
\newcommand{\wellBlockSetDummy}{\overbar{\mathcal D}}
\newcommand{\wellBlockSetProducer}{\wellBlockSet_P}
\newcommand{\wellBlockSetInjector}{\wellBlockSet_I}
\newcommand{\wellBlockArg}[1]{{d_{#1}}}

\newcommand{\unitVec}[1]{{\textbf{e}_{#1}}}

\newcommand{\fomOperatorOne}{g}
\newcommand{\romOperatorOne}{g_{\surr}}
\newcommand{\fomOperatorOneArg}[1]{\fomOperatorOne^{#1}}
\newcommand{\romOperatorOneArg}[1]{\romOperatorOne^{#1}}

\newcommand{\fomOperatorOneTimeN}{\fomOperatorOneArg{\timestep}}
\newcommand{\romOperatorOneTimeN}{\romOperatorOneArg{\timestep}}

\newcommand{\operatorTwo}{ q}
\newcommand{\fomOperatorTwo}{ \operatorTwo}
\newcommand{\romOperatorTwo}{\operatorTwo_{\surr}}
\newcommand{\corrOperatorTwo}{\operatorTwo_{\corr}}

\newcommand{\stateFom}{{\x}}
\newcommand{\stateRom}{\z}

\newcommand{\stateFomTime}[1]{{\stateFom}^{#1}}
\newcommand{\stateRomTime}[1]{{\stateRom}^{#1}}

\newcommand{\stateFomSpace}{\RR{\nstateFom}}
\newcommand{\stateFomSampledSpace}{\RR{\nSampleOutputState}}
\newcommand{\stateRomSpace}{\RR{\nstateRom}}

\newcommand{\outputFom}{q}
\newcommand{\outputRom}{\outputFom_{\surr}}

\newcommand{\outputFomSpace}{\RR{}}
\newcommand{\outputRomSpace}{\RR{}}

\SetKwData{Left}{end}
\SetKwData{This}{this}
\SetKwData{Up}{up}
\SetKwFunction{Union}{Union}
\SetKwFunction{FindCompress}{FindCompress}
\SetKwInOut{InputSettings}{Input settings}
\SetKwInOut{Inputs}{Inputs}
\SetKwInOut{Output}{Output}

\newcommand{\figRefOne}[1] {Figure~\ref{#1}}

\newcommand{\figRefTwo}[2] {Figures~\ref{#1} and \ref{#2}}

\newcommand{\eqnRefOne}[1] {Equation~(\ref{#1})}

\newcommand{\nTrain}{N_\text{train}}

\newcommand{\colorText}[2]{{\leavevmode\color{#1} #2}}

\newcommand{\blackText}[1]{\colorText{black}{#1}}

\definecolor{myGreen}{RGB}{26 148 49}

\newcommand{\cmmnt}[1]{\ignorespaces}

\newcommand{\overbar}[1]{\mkern 1.5mu\overline{\mkern-1.5mu#1\mkern-1.5mu}\mkern 1.5mu}

\newcommand{\coarse}{{c}}

\newcommand{{\permC}}{\perm^{*}}

\newcommand{{\TC}}{T^{*}}
\newcommand{{\sC}}{S^{\coarse}}
\newcommand{{\fC}}{f^{\coarse}}
\newcommand{{\pC}}{p^{\coarse}}
\newcommand{{\vC}}{\phaseVelocity{}^{\coarse}}
\newcommand{{\phiC}}{\phi^{*}}

\newcommand{{\sBar}}{\overbar{S}}
\newcommand{{\pBar}}{\overbar{p}}
\newcommand{{\fBar}}{\overbar{f}}
\newcommand{{\vBar}}{\overbar{\phaseVelocity{t}}}

\newcommand{{\vPrime}}{\phaseVelocity{t}'}

\newcommand{{\sVar}}{\sigma_{S}}
\newcommand{{\pVar}}{\sigma_{p}}
\newcommand{{\fVar}}{\sigma_{f}}
\newcommand{{\vVar}}{\sigma_{\phaseVelocity{}}}
\newcommand{{\vsCovar}}{\sigma_{\phaseVelocity{}S}}

\newcommand{\inNaturalSequence}[1]{=1,\ldots,#1}

\newcommand{\reviewer}[1]{{\blackText{#1}}}

\begin{document}

\title{Error \reviewer{modeling for surrogates of dynamical systems using machine learning}}
\author{Sumeet Trehan\affil{1}\corrauth,
Kevin Carlberg\affil{2}\ and Louis J.\ Durlofsky\affil{1}}
\address{\affilnum{1} Department of Energy Resources Engineering, Stanford
University, 367 Panama Street, Stanford, CA 94035\break
\affilnum{2}Extreme-Scale Data Science and Analytics Department, Sandia
National Laboratories, 7011 East Ave, MS 9159, Livermore, CA 94550}
\corraddr{Sumeet Trehan, Department of Energy Resources Engineering, Stanford
University, 367 Panama Street, Stanford, CA 94035. Email: strehan@alumni.stanford.edu}

%\author[1]{Sumeet Trehan\corrauth}
%\author[2]{Kevin Carlberg}
%\author[1]{Louis J. Durlofsky}
%\affil[1]{Department of Energy Resources Engineering, Stanford University}
%\affil[2]{Sandia National Laboratories}
%\author{Sumeet Trehan\corrauth, Kevin T. Carlberg, Louis J. Durlofsky}
%\address{367 Panama Street, Stanford, CA 94035}
%\renewcommand\Authands{ and }

%\corraddr{strehan@stanford.edu}

\begin{abstract}
A machine-learning-based framework for \reviewer{modeling} the
error introduced by surrogate models of parameterized dynamical systems is proposed. The
framework entails the use of high-dimensional regression techniques (e.g., random
forests, LASSO) to map a large set of inexpensively computed `error indicators'
(i.e., features) produced by the surrogate model at a given time instance to a prediction of the
surrogate-model error in a quantity of interest (QoI). \reviewer{This
eliminates the need for the user to hand-select a small number of informative
features.} The
methodology requires a training set of parameter instances at which the
time-dependent surrogate-model error is computed by simulating both the
high-fidelity and surrogate models. Using these training
data, the method first \reviewer{determines regression-model locality} (via classification
or clustering), and subsequently constructs a `local' regression model to
predict the time-instantaneous error within each identified region of feature space.  We consider
two uses for the resulting error model: (1) as a correction to the surrogate-model
QoI prediction at each time instance, and (2) as a way to statistically model
arbitrary functions of the time-dependent surrogate-model error (e.g., time-integrated
errors). We apply the proposed framework to \reviewer{model} errors in reduced-order models of nonlinear oil--water subsurface flow simulations,
with time-varying well-control (bottom-hole pressure) parameters. The reduced-order models used in this work entail application of trajectory piecewise linearization in conjunction with proper orthogonal decomposition. When the first use of the method is considered, numerical experiments demonstrate consistent improvement in accuracy in the time-instantaneous QoI prediction relative to the original
surrogate model, across a large number of test cases. When the second use is
considered, results show that the
proposed method provides accurate statistical predictions of the 
time- and well-averaged errors.  

\end{abstract}

\maketitle

%----------------------------------------------
% Paper begins
%----------------------------------------------

%----------------------------------------------
\section{Introduction}
%----------------------------------------------

% General
Computational simulation is being increasingly employed for real-time and
many-query problems such as design, optimal control, uncertainty
quantification, and inverse modeling. However, these applications require the
rapid and repeated simulation of a (parameterized) computational model, which
often corresponds to a discretization of partial differential equations
(PDEs) that can be nonlinear, time dependent, and multiscale in nature. Accurate
predictive models can therefore incur substantial computational costs. While
recent advances in parallel computing have reduced simulation times for 
high-fidelity models, the rapid, repeated simulation of such models remains
a bottleneck in many applications. 

This computational challenge has motivated the development of a wide range of
surrogate-modeling methods. Surrogate models---which can be categorized as
data fits, lower-fidelity models, or reduced-order models (ROMs)---are
approximations of high-fidelity models (HFMs) that aim to provide large
computational savings while preserving accuracy.
Unfortunately, these models often introduce non-negligible errors due to the
assumptions and approximations employed in their construction, and these errors can have
deleterious effects on the resulting analysis. Thus, to 
more rigorously employ surrogate models, it is critical to quantify the
errors they
introduce. Without reliable error quantification\reviewer{---which can be accomplished
via statistical approaches, rigorous error bounds, error
indicators, or error models---}the accuracy of
surrogate-model predictions is unknown, and the trustworthiness of the
resulting analysis may be questionable.

For this reason, researchers have developed a variety of methods for
\reviewer{quantifying}
surrogate-model errors. Data-fit surrogate models construct a deterministic
function (e.g., polynomial fit \cite{knill1999response}, artificial neural
network \cite{watson1996ann}), or stochastic process (e.g., Gaussian
process/kriging \cite{kennedy2001bayesian,arendt2012quantification}) that
explicitly approximates the mapping from the model input parameters to the model
output quantity/quantities of interest (QoI). For this class of surrogate model,
statistical approaches such as the $R^2$ value~\cite{james2013introduction},
cross validation~\cite{james2013introduction, hastie2005elements}, confidence
intervals, and prediction intervals~\cite{james2013introduction,
hastie2005elements} can be applied to quantify surrogate error. When the data fit
corresponds to a stochastic-process model, the prediction variance can be
applied to quantify the uncertainty in the prediction directly
\cite{j2006design, forrester2007multi}. Although data-fit surrogates are nonintrusive to implement (they require only `black-box'
queries of the HFM), their predictions are not physics based, which can
lead to inaccurate predictions, especially for high-dimensional
input-parameter spaces. For this reason, many applications demand more
sophisticated physics-based
surrogates such as lower-fidelity or reduced-order models. 
%used (co-)kriging (also known as Gaussian-process
%regression) as a lower-fidelity model, and derived an error estimator for the
%model. 
%Similarly, to construct the error model, Watson and
%Gupta~\cite{watson1996ann} used artificial neural networks, Knill et
%al.~\cite{knill1999response} used response surfaces, Kennedy et
%al.~\cite{kennedy2001bayesian} and Arendt et
%al.~\cite{arendt2012quantification} used kriging, 

Lower-fidelity models are physics-based surrogates that apply simplifications
to the original HFM, such as coarser discretizations, lower-order
approximations or linearization, or they neglect some physics. In this case, one approach
for \reviewer{quantifying the error would be to explicitly model the error via a `data-fit'} mapping between
input parameters and lower-fidelity-model error. These
approaches---which have been pioneered in the field of multifidelity design
optimization---typically enforce `global' zeroth-order consistency between the
corrected low-fidelity-model QoI and the HFM QoI at training points
\cite{gano2005hybrid,huang2006sequential,march2012provably,rajnarayan2008multifidelity,NE12, lodoen2005assessment, omre2004improved}, 
or `local' first- or second-order consistency at trust-region centers
\cite{Alexandrov2001, eldred:soc}. Unfortunately, in the context of dynamical
systems (which we consider in this work), such corrections are only applicable
to scalar-valued QoI (e.g., time-averaged quantities). Quantifying the error in
time-dependent quantities has been largely
ignored, with the exception of recent work that interpolates time-dependent error models in the input-parameter space \cite{pagani2016reduced}. 
In addition, the error may exhibit a complex or oscillatory dependence on the input
parameters, which adds further challenges in high-dimensional input-parameter spaces and can cause the approach to fail~\cite{NE12,
drohmann2015romes}.
Alternatively, adjoint-based error estimation (i.e., dual-weighted-residual
error indicators) can also be applied to \reviewer{approximate} QoI errors in
the context of coarse
finite-element
\cite{babuvska1984post,becker1996weighted,rannacher1999dual,bangerth2003adaptive},
finite-volume \cite{venditti2000adjoint,venditti2002grid,park2004adjoint}, and
discontinuous Galerkin \cite{lu2005posteriori,fidkowski2007simplex}
discretizations. Unfortunately, for dynamical systems, dual-weighted residuals
require an additional (time-local or time-global) dual solve, which can incur
a non-negligible additional cost.

Projection-based ROMs apply projection to reduce the dimensionality
of the equations governing the HFM. Typically, the low-dimensional bases are
derived empirically by evaluating the HFM at training points or by performing
other analyses, e.g., solving Lyapunov equations or computing a Krylov subspace.
ROM error is typically \reviewer{estimated} by deriving rigorous \textit{a
priori} and \textit{a posteriori} error bounds for the state, QoI, or transfer function;
such bounds have been derived for the reduced-basis method \cite{BM12, GP05, RHP07},
system-theoretic approaches (e.g., balanced truncation, rational interpolation)
\cite{antoulas2005approximation}, and proper-orthogonal-decomposition (POD) Galerkin \cite{volkweinPODanalysis,rathinam:newlook}
and least-squares Petrov--Galerkin (LSPG) \cite{carlbergGalDiscOpt} methods; see, e.g., \cite{surveyWillcoxGugercin} for a review. However---for dynamical systems---such error bounds
typically grow exponentially in time, causing the bound to significantly
overpredict the error \cite{drohmann2012reduced}, which can limit the practical utility of
these bounds.
Note that the approaches developed for \reviewer{quantifying} the low-fidelity-model error could also
be adopted to quantify ROM errors \cite{NE12,carlberg2014adaptive}, as ROMs can be interpreted as lower-dimensional models with
empirically derived basis functions.

To address the above issues in the context of quantifying ROM errors, Drohmann
and Carlberg~\cite{drohmann2015romes} devised the reduced-order-model error
surrogates (ROMES) approach\reviewer{, which can be considered an error-modeling
approach, as it} constructs a Gaussian process that maps
\textit{error indicators} (e.g., error bounds, residual norms, dual-weighted
residuals) produced inexpensively by the ROM to a \textit{distribution over
the ROM QoI error.} This study demonstrated that---even in the presence of
high-dimensional input-parameter spaces---the ROM produces a small number of
inexpensively computable error indicators that can be employed to derive
accurate, low-variance predictions of the ROM error. Follow-on work also
investigated the use of statistical modeling and regression methods to
\reviewer{model} ROM errors from indicators \cite{manzoni2014accurate}. While
promising, the ROMES approach requires the user to hand-select the error
indicators that are most strongly correlated with the ROM error; i.e., feature
selection is left to the user.  \reviewer{This task can be challenging in
general applications, as the user may not have strong \textit{a priori}
knowledge of which (small number of) features inform the error.} Further, the
ROMES method was demonstrated only on a stationary (i.e., steady-state)
problem, and its extension to dynamical systems is nontrivial. Finally,
application of the ROMES method to other physics-based surrogates (e.g.,
coarsened or upscaled models) is not obvious, as different error indicators
will likely be required for such surrogates.

In this paper, we propose a machine-learning-based framework for
\reviewer{modeling} the
error introduced by physics-based surrogate models of dynamical systems. The
approach applies statistical techniques for high-dimensional regression (e.g.,
random forests) to map a large set of inexpensively computed quantities or
features generated by the surrogate model to a prediction of the
time-instantaneous surrogate-model QoI error.  This method is referred to as
error \reviewer{modeling} via machine learning (EMML). In contrast to the ROMES
method, the proposed EMML approach 
\reviewer{enables a large number of potential error indicators to be included
in the candidate feature set and thus}
does not require \reviewer{the user to
manually select a small number of features that inform the error, which can be
challenging as mentioned above}.
Thus, feature selection is included in the process of constructing the error
model, and extension to multiple types of physics-based surrogates is
straightforward---we assume only that the surrogate produces a large set of
features that can be mined for potential error indicators. 

While the proposed framework can be applied to any surrogate model in principle, in this
study we apply the method to \reviewer{model} errors introduced by a TPWL ROM
\cite{RewienskiTPWL,
dong2008general, CardosoTPWL, he2011enhanced} 
applied with a proper orthogonal decomposition (POD) basis and LSPG projection \cite{CarlbergPG,carlbergGalDiscOpt}, which we refer to hereon as
POD--TPWL. \reviewer{We employ this ROM rather than other approaches (e.g.,
(D)EIM with POD--Galerkin, GNAT) because TPWL is less intrusive: 
it simply requires extracting linear operators (i.e., Jacobians
of the residual with respect to the current state, previous state, and controls) from the HFM
simulation code during the offline training stage, and it is entirely
independent of 
the HFM simulation code during the online stage.} This particular ROM was also employed in previous subsurface-flow
studies involving oil--water and oil--gas compositional problems~\cite{he2014reduced,he2015constraint,trehanTPWQ}. At each time step during test
simulations, POD--TPWL performs
linearization around the nearest training solution; LSPG
projection is then applied---with a low-dimensional POD basis---to
reduce the dimensionality of the linearized system.  While the
approach has been shown to yield $\orderOf(10^{2})-\orderOf(10^{3})$ speedups~\cite{CardosoTPWL, he2011enhanced,he2014reduced}, POD--TPWL incurs non-negligible errors
due to the approximations it introduces, namely (1) linearization error, and
(2) `out-of-plane error' \cite{rathinam:newlook} arising from employing a
low-dimensional POD trial subspace, (3) `in-plane error' arising from
projection, and (4) error propagated from the previous
time step. See \cite{he2015constraint} for further discussion of POD--TPWL
error. We aim to apply the proposed EMML framework to \reviewer{model} the QoI error
resulting from these approximations. 

Finally, we note that while machine learning and its application across various
disciplines have been extensively studied, the use of machine learning within
the domain of physics-based modeling and simulation is relatively new, although
it is quite promising \cite{doeReport,doeReportML}. In the context of improving
surrogate-model predictions, researchers have recently investigated the use of
machine learning techniques to identify (via classification) the spatial
locations for high low-fidelity-model error~\cite{ling2015evaluation}. Machine
learning has also been used to quantify the inherent error with kriging
data-fit surrogates~\cite{tracey2013application} and to derive improved closure
models in the context of computational fluid
mechanics~\cite{weatheritt2015development, weatheritt2016novel,
tracey2015machine, duraisamy2015new, ling2016reynolds}. Machine learning was also applied to
derive the source term for the transport of an intermittency variable while
transitioning from laminar to turbulent flow~\cite{duraisamy2014transition}.
Similarly, regression techniques (e.g., LASSO~\cite{tibshirani1996regression})
have been used to calibrate and infer uncertainties in viscosity-model
coefficients for transonic flow applications~\cite{lefantzi2015eddy}. We note
that although machine learning has been used in a post-processing step to
identify the physical regions of high surrogate-model
error~\cite{ling2015evaluation}, to our knowledge, the direct
\reviewer{approximation} of
error through construction of a regression model has not been pursued. 

The paper proceeds as follows. In Section~\ref{section:problemDescription} we
describe the general EMML framework. Following the definition of the error in
the QoI, we introduce four methods to map this error to a set of inexpensively
computed features using high-dimensional regression techniques.
The particular HFM (subsurface flow) and surrogate model (POD--TPWL) considered as an application are discussed in Section~\ref{section:governingEquations}. In
Section~\ref{section:numericalResults}, we present numerical results for
\reviewer{modeling} the POD--TPWL error in flow quantities driven by time-varying control
variables over a large number of test cases. Different algorithmic treatments
are also considered. A summary and suggestions for future work are provided in
Section~\ref{section:conclusion}. \reviewer{Appendices~A and B present descriptions of random-forest and LASSO regression.}

%----------------------------------------------
\section{General Problem Statement}
\label{section:problemDescription}
%----------------------------------------------

In this section, we describe the overall EMML framework for error
\reviewer{modeling}. We
begin by introducing both the high-fidelity model (HFM) and the surrogate model in
a general setting. 
%
%
%------------------------------------
\subsection{Dynamical-system high-fidelity model}
%------------------------------------
Given input 
parameters $\paramsGen\in\paramSpace$, we assume that the HFM generates a time
history of states $\stateFomTime{n}\in\stateFomSpace$,
$n=1,\ldots,\numTimesteps$, and a scalar-valued output QoI $\outputFom$
that depends on the state, i.e.,
\begin{align}
 \label{eqn:FOMGeneralEquations}
	\begin{split}
		\fomOperatorOneTimeN &: \paramsGen\mapsto
		\xN,\quad n=1,\ldots,\numTimesteps,\\
		 &: \paramSpace\rightarrow\stateFomSpace,
	\end{split}\\
	\begin{split}
	\label{eq:fomOutput}
		\fomOperatorTwo &: \sampleOutputState\stateFom \mapsto
		\outputFom(\sampleOutputState\stateFom),\\
		 &: \stateFomSampledSpace\rightarrow\outputFomSpace.
	\end{split}
\end{align}
Here, $\sampleOutputState\in \{0,1\}^{\nSampleOutputState\times\nstateFom}$
with $\nSampleOutputState\leq\nstateFom$ is a sampling matrix comprising
selected rows of the identity matrix. This operator extracts the elements of
the state vector required to compute the QoI $\outputFom$. 
%
%
%------------------------------------
\subsection{Dynamical-system surrogate model}
%------------------------------------
We assume that the inexpensive surrogate model generates a time history of
surrogate-model states $\stateRomTime{n}\in\stateRomSpace$,
$n=1,\ldots,\numTimesteps$ (ideally with $\nstateRom\ll \nstateFom$), given the input parameters
$\paramsGen\in\paramSpace$. The QoI can be computed from the surrogate model
states $\stateRomTime{}$ using the function $\outputRom$. Our critical
assumption is that the surrogate model also produces \textit{auxillary data} in the form
of  a time history of `features' $\featuresN\in\featuresSpace$,
$n=1,\ldots,\numTimesteps$, given these parameters $\paramsGen$. The surrogate
model is described as follows:
\begin{align}
\label{eqn:ROMGeneralEquations}
	\begin{split}
		\romOperatorOneTimeN &: \paramsGen \mapsto\zN,\quad n=1,\ldots,\numTimesteps,\\
		 &: \paramSpace\rightarrow\stateRomSpace,
	\end{split}
	\\
	\begin{split}
	\label{eqn:romOuput}
	\romOperatorTwo &: \stateRom \mapsto \outputFom (\sampleRomOutputState\stateRom),\\
		 &: \stateRomSpace\rightarrow\outputRomSpace,
	\end{split}
	\\
	\begin{split}
	\label{eqn:ROMFeatureConstruction}
		\romOperatorGenerateFeaturesArg{\timestep} &: \paramsGen \mapsto
		\romOperatorGenerateFeaturesArg{\timestep}(\paramsGen),\quad n=1,\ldots,\numTimesteps, \\
		 &: \paramSpace\rightarrow\featuresSpace.
	\end{split}
\end{align}
Here, $\sampleRomOutputState\in\RR{\nSampleOutputState\times\nstateRom}$
denotes a prolongation operator that transforms the surrogate model state into
the elements of the high-fidelity state required to compute the QoI. 
%Note that
%the features may include information such as the surrogate model state and
%parameters. 
The decision of what to include in the set of features is motivated by the underlying form of
$\romOperatorOneTimeN$, as discussed in Section \ref{subsection:Features}. For
notational simplicity, from hereon we use
$\romOperatorGenerateFeaturesArg{\timestep}$ in place of
$\romOperatorGenerateFeaturesArg{\timestep}(\paramsGen)$. 
%
%------------------------------------
\subsection{Error modeling} 
\label{subsection:ErrorModeling}
%------------------------------------
Our objective is to predict the error in the QoI at each time step $n$, which
we define as
\begin{align}
\label{eq:qoiError}
		\errorOperatorQoITimeN(\paramsGen) \defeq \fomOperatorTwo^{\timestep}(\paramsGen)
		 - \romOperatorTwo^{\timestep}(\paramsGen),\quad n=1,\ldots,\numTimesteps,
\end{align}
where $\errorOperatorQoITimeN \in \RR{}$, $\fomOperatorTwo^{\timestep}
(\paramsGen)\defeq \fomOperatorTwo
\circ\sampleOutputState\fomOperatorOneTimeN(\paramsGen)$, denotes the QoI
computed using the HFM at time instance $n$, and $\romOperatorTwo^{\timestep}
(\paramsGen)\defeq \romOperatorTwo \circ\romOperatorOneTimeN(\paramsGen)$ 
denotes the QoI computed using the surrogate model at time instance $n$.  We
propose to \reviewer{approximate} this QoI error as
$\errorOperatorQoIApproxTimeN\approx\errorOperatorQoITimeN,\
n=1,\ldots,\numTimesteps$, by constructing error surrogates using
high-dimensional regression methods developed in the context of machine
learning.  

In particular, we propose four regression-based approaches that construct a mapping
from the surrogate-model features to a QoI-error prediction
$\errorOperatorQoIApproxTimeN$. While it is possible to construct an
error prediction
$\errorOperatorQoIApproxTimeN$
as a function of 
only the input parameters $\paramsGen$, this approach can fail due to the oscillatory behavior of certain surrogate-model errors in the input space \cite{NE12, drohmann2015romes}, as discussed in the Introduction.  For notational
simplicity, hereon $\errorOperatorQoITimeN =
\errorOperatorQoITimeN(\paramsGen)$. 
%--------------
\subsubsection{Method 1: QoI error.}
\label{subsubsection:Method1}
%--------------
The first method models the nondeterministic mapping 
$\romOperatorGenerateFeaturesArg{\timestep}\mapsto\errorOperatorQoITimeN$ 
as a sum of a deterministic function $\errorModelQoI:
\featuresSpace\rightarrow\RR{}$, and nondeterministic noise $\errorMapNoise$ as follows: 
\begin{equation}
\label{eq:errorOperatorQoITimeNNoise}
	\errorOperatorQoITimeN =
	\errorModelQoI(\romOperatorGenerateFeaturesArg{\timestep}) +
	\errorMapNoise,\quad n=1,\ldots,\numTimesteps.
\end{equation}
Here, $\errorMapNoise$ is a zero-mean random variable that
accounts for  potentially unknown features, deficiencies in the model form of
$\errorModelQoI$,
and the error introduced due to sampling variability. Thus,
$\errorMapNoise$ denotes irreducible error induced by the error
model \eqnRefOne{eq:errorOperatorQoITimeNNoise}, and it may in principle depend on the features (i.e., $\errorMapNoise =
\errorMapNoise(\romOperatorGenerateFeaturesArg{\timestep})$); this
heteroscedasticity can occur, for example, when the mean and variance of the
predicted error are larger for larger-magnitude features. For simplicity, we
neglect this dependence in the current study. 

We note that employing the error model in
\eqnRefOne{eq:errorOperatorQoITimeNNoise} allows us to (1)
\reviewer{approximate} the form of
$\errorModelQoI$ using data, which in turn enables us to express the error in
\reviewer{the} QoI as a function of the features only, and (2) account for the possibility that the same
feature vector may yield different values of the QoI error, i.e.,
$\errorOperatorQoIArg{n} \neq \errorOperatorQoIArg{m}$ with $n\neq m$ but
$\romOperatorGenerateFeaturesArg{n} = \romOperatorGenerateFeaturesArg{m}$.

Next, we construct a model
$\errorModelQoIApprox:\featuresSpace\rightarrow\RR{}$ of the function
$\errorModelQoI$, such that
$\errorModelQoIApprox(\romOperatorGenerateFeaturesArg{\timestep}(\paramsGen))\approx
\errorModelQoI(\romOperatorGenerateFeaturesArg{\timestep}(\paramsGen)), n=1,\ldots,\numTimesteps$, $\forall \paramsGen\in\paramSpace$. This model allows the error to be approximated as a function of the features as
\begin{gather}
	\errorOperatorQoIApproxTimeN =\errorModelQoIApprox(\romOperatorGenerateFeaturesArg{\timestep})
	,\quad n=1,\ldots,\numTimesteps.
\end{gather} 

Note that 
we consider \reviewer{modeling} $\errorOperatorQoITimeN$ as a prediction problem rather
than as a time-series-analysis problem. This is because, in the problem under
consideration, we perform numerical experiments for different input-parameter
instances over the same time interval.
Thus, we include time in
the feature set, as described in Section~\ref{subsection:Features}.
Implicitly, we assume that samples are independent and identically distributed
(i.i.d.), with each sample corresponding to the quantity at a given time step.

%--------------
\subsubsection{Method 2: relative QoI error.}
\label{subsubsection:Method2}
%--------------
In many cases, the QoI errors $\errorOperatorQoITimeN(\paramsGen)$ can
exhibit a wide range of observed values. This can make the machine-learning
task more challenging, as the associated regression model must be predictive
across this entire range of values. To address this, we can instead apply
regression to the \textit{relative} QoI errors---which typically exhibit a
narrower range of values---and subsequently approximate the QoI 
errors in a postprocessing step.

We define the relative QoI error at time step $n$ by
\begin{align}
\label{eq:relErrorQoI}
		\relerrorOperatorQoITimeN \defeq
		\dfrac{\errorOperatorQoITimeN}{\fomOperatorTwo(\sampleOutputState\xN)
		},\quad n=1,\ldots,\numTimesteps,
\end{align}
and---following Method~1---express the mapping
$\romOperatorGenerateFeaturesArg{\timestep}\mapsto\relerrorOperatorQoITimeN$ as
\begin{equation}
\label{eq:errorOperatorQoITimeNNoiseRel}
	\relerrorOperatorQoITimeN =
	\relerrorModelQoI(\romOperatorGenerateFeaturesArg{\timestep}) +
	\errorMapNoise,\quad n=1,\ldots,\numTimesteps,
\end{equation}
where $\relerrorModelQoI: \featuresSpace\rightarrow\RR{}$ is an unknown
deterministic function. We \reviewer{model} 
$\relerrorModelQoI$ by constructing an approximation
$\relerrorModelQoIApprox:\featuresSpace\rightarrow\RR{}$ such that
$\relerrorModelQoIApprox(\romOperatorGenerateFeaturesArg{\timestep}(\paramsGen))\approx
\errorModelQoI(\romOperatorGenerateFeaturesArg{\timestep}(\paramsGen))$,
$n=1,\ldots,\numTimesteps$, $\forall \paramsGen\in\paramSpace$. This 
allows the \reviewer{approximated} relative error to be expressed as
\begin{gather} 
	\relerrorOperatorQoIApproxTimeN =\relerrorModelQoIApprox(\romOperatorGenerateFeaturesArg{\timestep})
	,\quad n=1,\ldots,\numTimesteps.
\end{gather} 
From \eqnRefOne{eq:qoiError}, it follows that
\begin{align}
	\fomOperatorTwo (\sampleOutputState\xN) = 
	\dfrac{\fomOperatorTwo(\sampleRomOutputState\zN)}{1- \relerrorOperatorQoITimeN},\quad
	n=1,\ldots,\numTimesteps,
\end{align}
which allows the QoI error to be related to the relative QoI error by
\begin{equation} 
\label{eq:relToAbs}
	\errorOperatorQoITimeN = \fomOperatorTwo(\sampleRomOutputState\zN) 
	\left( \dfrac{\relerrorOperatorQoITimeN}
	{1-\relerrorOperatorQoITimeN} \right),\quad
	n=1,\ldots,\numTimesteps. 
\end{equation} 
Therefore, we can \reviewer{model} the QoI error as $\errorOperatorQoIApproxTimeN$ from the relative QoI regression
model $\relerrorOperatorQoIApproxTimeN$ in a postprocessing step as
\begin{equation} 
\label{eq:relToAbsModel}
	\errorOperatorQoIApproxTimeN = \fomOperatorTwo(\sampleRomOutputState\zN)
	\left( \dfrac{\relerrorOperatorQoIApproxTimeN}{1-\relerrorOperatorQoIApproxTimeN} \right),
	\quad n=1,\ldots,\numTimesteps.
\end{equation}

%--------------
\subsubsection{Method 3: state error.}\label{sec:method3}
%--------------
As an alternative to modeling the QoI error, we can model the error in the
relevant state(s), and then use this quantity to \reviewer{approximate} the QoI
error. This method is advantageous when the QoI error exhibits a more complex
behavior than the state error, as may be the case for highly nonlinear QoI. We
denote the error in the sampled state at time step $\timestep$ by 
\begin{align}
\label{eqn:MappingPrimaryVar}
		[\errorOperatorPrimaryVarTimeN]_{i} &\defeq [\sampleOutputState\xN]_i-
		[\sampleRomOutputState\zN]_i,\quad n=1,\ldots,\numTimesteps,~i
		\inNaturalSequence{\nSampleOutputState},
\end{align}
where $\errorOperatorPrimaryVarTimeN \in \RR{\nSampleOutputState}$ and 
$[\cdot]_i$ denotes the $i$th element of the vector-valued argument. 

The next steps follow closely the derivation of Method~1. As before, we model
the mappings $\romOperatorGenerateFeaturesArg{\timestep}\mapsto
[\errorOperatorPrimaryVarTimeN]_i$ as
\begin{equation}
\label{eq:errorOperatorPrimaryTimeNNoise}
	[\errorOperatorPrimaryVarTimeN]_i = \errorModelPrimaryVar{i}(\romOperatorGenerateFeaturesArg{\timestep}) + \errorMapNoise,
	\quad n=1,\ldots,\numTimesteps,~i
		\inNaturalSequence{\nSampleOutputState},
\end{equation}
where $\errorModelPrimaryVar{i}:\featuresSpace\rightarrow\RR{}$,
$i\inNaturalSequence{\nSampleOutputState}$, denote unknown
deterministic functions that allow the state-variable errors to be computed
as a function of the features. Analogous to
\eqnRefOne{eq:errorOperatorQoITimeNNoise}, we construct regression models
$\errorModelPrimaryVarApprox{i}:\featuresSpace\rightarrow\RR{}$,
such that
$\errorModelPrimaryVarApprox{i}(\romOperatorGenerateFeaturesArg{\timestep}(\paramsGen))\approx
\errorModelPrimaryVar{i}(\romOperatorGenerateFeaturesArg{\timestep}(\paramsGen))$,
$n=1,\ldots,\numTimesteps$, $i\inNaturalSequence{\nSampleOutputState}$,
$\forall\paramsGen\in\paramSpace$. This model allows the
\reviewer{approximated} state error to be expressed as 
\begin{gather} 
	[\errorOperatorPrimaryVarApproxTimeN]_i =\errorModelPrimaryVarApprox{i}(\romOperatorGenerateFeaturesArg{\timestep})
	,\quad n=1,\ldots,\numTimesteps,~i
		\inNaturalSequence{\nSampleOutputState}.
\end{gather} 

Note that 
instead of pursuing multi-response multivariate regression, we execute
$\nSampleOutputState$ independent multivariate regressions, i.e.,
we construct a unique and independent mapping
$\errorModelPrimaryVarApprox{i}$ for each of the
$\nSampleOutputState$ sampled states. Because the QoI error can be expressed as
\begin{align} 
\label{eq:methodThreeFinal}
	\errorOperatorQoITimeN&=
	\fomOperatorTwo \left(\sampleRomOutputState\zN + \errorOperatorPrimaryVarApproxTimeN
\right) - \fomOperatorTwo (\sampleRomOutputState\zN),\quad
	n=1,\ldots,\numTimesteps,
\end{align} 
we can \reviewer{model} the QoI error from the \reviewer{modeled} state error in a postprocessing step as 
\begin{align} 
\label{eq:methodThreeFinalEnd}
	\errorOperatorQoIApproxTimeN&=
	\fomOperatorTwo \left(\sampleRomOutputState\zN + \errorOperatorPrimaryVarApproxTimeN
\right) - \fomOperatorTwo (\sampleRomOutputState\zN),\quad
	n=1,\ldots,\numTimesteps.
\end{align}
%

%--------------
\subsubsection{Method 4: relative state error.}\label{sec:method4}
%--------------
Finally, if the errors in a particular sampled state exhibit a wide range of
observed values, we can construct a regression model for the relative state
error, which we define as
\begin{align}
	\begin{split}
		[\relerrorOperatorPrimaryVarTimeN]_i &\defeq
		\dfrac{[\errorOperatorPrimaryVarTimeN]_i}{[\sampleOutputState\xN
		]_i},\quad n=1,\ldots,\numTimesteps,~i
		\inNaturalSequence{\nSampleOutputState},
	\end{split}
\end{align}
with $\relerrorOperatorPrimaryVarTimeN\in \RR{\nSampleOutputState}$, and subsequently use this model to predict the QoI error.
Analogous to Methods 1--3, we construct a
regression model, which maps the features $\featuresN$ to the relative error
in the sampled state $[\relerrorOperatorPrimaryVarTimeN]_i$, thus enabling the
computation of the \reviewer{approximated} relative error
$[\relerrorOperatorPrimaryVarApproxTimeN]_i$. Analogous to Method~2, we relate
the relative error in the sampled state to the actual error. After using
Method~3 or Method~4 or some combination thereof to compute the value of
$[\errorOperatorPrimaryVarApproxTimeN]_i,~\forall i
\inNaturalSequence{\nSampleOutputState}$, we apply
\eqnRefOne{eq:methodThreeFinalEnd} to determine the QoI error \reviewer{model}
$\errorOperatorQoIApprox$.

%--------------
\subsection{Training data}
%--------------
Each of the four methods described above entails the use of a regression model to
predict the output (response)---which corresponds to the actual or relative
error in the QoI or in the sampled state---given the inputs (features). Constructing any such regression model relies on training data. We denote points in the input-parameter space used to collect these data as
\begin{equation}
	\nonumber\paramTrain\defeq\{\paramsTrainArg{1},\ldots,\paramsTrainArg{\numHFSsymbol}\}\subseteq\paramSpace,
\end{equation}
where $\paramsGenArg{i} \in \paramSpace$, $i=1,\ldots,\numHFSsymbol$ denotes
the $i$th EMML training instance of the input parameters, and $\numHFSsymbol$
denotes the number of training points. 

Next, we simulate both the HFM and the surrogate model for input-parameter
instances in the training set $\paramTrain$.
This produces the EMML training data, which comprise errors $\errorOperatorQoITimeN,
\errorOperatorPrimaryVarTimeN$, and features $\featuresN$ over all time steps
and training simulations, i.e., 
$$
	\nonumber \{  \left(\errorOperatorQoITimeN(\paramsGen),
	~\errorOperatorPrimaryVarTimeN(\paramsGen), ~
	\romOperatorGenerateFeaturesArg{\timestep}(\paramsGen)\right)
	\}_{\paramsGen\in\paramTrain,\ n=1,\ldots,\numTimesteps}.
$$
As mentioned previously, we assume that the associated samples are i.i.d. In the
following, we denote a general training error by $\errorGen^n(\paramsGen)$,
$\paramsGen\in\paramTrain,\ n=1,\ldots,\numTimesteps$.
% In the next sections, we first describe the feature design to estimate the POD-TPWL error. Next, we describe space partitioning, followed by the construction of the EMML training data (involving simulating HFMs and POD-TPWL models for some instances of $\paramsGen$). Finally, we describe machine learning techniques such as random forest and LASSO regression.

%-------------------------------------------
\subsection{\reviewer{Regression-model locality}}
\label{subsection:spacePartitioning}
%-------------------------------------------
As an alternative to constructing a single global regression function,
$\errorModelPrimaryVar{i}, \relerrorModelPrimaryVar{i}, \errorModelQoI$, or
$\relerrorModelQoI$, that is valid over the entire feature space, we can instead
construct multiple `local' \reviewer{regression} models that are tailored to
particular 
\reviewer{physical regimes or}
feature-space regions  with the intent of improving prediction
accuracy.  To realize this, we
partition the training data into subsets corresponding to different
feature-space regions and construct separate regression functions for each
subset. We consider two methods for \reviewer{determining regression-model locality}:
classification and clustering.

Classification is a supervised machine learning
technique~\cite{hastie2005elements} that predicts the label (or category) to
which an observation belongs. In this work it entails constructing a statistical model from
EMML training data containing samples whose category membership is known, along
with categorization criteria for those samples. In the current context, we
propose applying classification using `classification features'
$\featuresClassify\in \RR{1\times\numFeaturesClassify}$ that may in general be
different from the EMML features $\features$.  We employ these features to
identify the subsets of the EMML training data associated with different
physical regimes of the problem, for which different regression models are
appropriate. Then, given a new observation for which we require
error prediction, we first identify its category using the classification
model, and subsequently apply the associated local regression model for error
prediction.

Clustering is an unsupervised machine learning method~\cite{hastie2005elements}
that can be applied to partition the training data according to the (e.g.,
Euclidean) distance in feature space between the training samples. We propose
partitioning feature space (or a lower-dimensional space embedded in feature
space computed, e.g., via principal component analysis) according to the
Voronoi diagram produced by the cluster means. \reviewer{In this work, we
employ $k$-means clustering where the number of clusters $k$ is
determined by 
identifying the elbow of the 
curve reporting the 
sum of squared errors as a function of the number of clusters
\cite{hastie2005elements}}.  Given a new observation, we
first identify its cluster from the Voronoi diagram, and then apply the local
regression model.
%
%-------------------------------------------
\subsection{High-dimensional regression methods}\label{sec:regressionMethods}
While, in principle, standard regression methods (e.g., linear regression)
could be used to construct error surrogates $\errorModelPrimaryVar{i},
\relerrorModelPrimaryVar{i}, \errorModelQoI$, and $\relerrorModelQoI$, such
approaches may be ineffective when the number of candidate features
$\numFeatures$ produced by the surrogate model is large. This occurs, for
example, when a projection-based reduced-order model is employed as the
surrogate. This ineffectiveness arises due to (1) the lack of available
guidelines for feature-subset selection, (2) the time-consuming and
challenging nature of \textit{a priori} identification of the relevant subset
of features, and (3) the fact that the response--feature relationship may depend on
nonlinear interactions between a large number of features. We therefore
propose applying high-dimensional regression methods that incorporate
\reviewer{automatic feature selection.} 

A wide range of methods---such as tree-based
methods (gradient boosting, random forests), support vector machines,
$K$-nearest neighbors, \reviewer{elastic-net}, and artificial neural network---fits
within this category. While the
specific choice of regression technique depends on the problem at hand, we
pursue two specific methods in this work: random-forest regression (RF) and 
LASSO (least
absolute shrinkage and selection operator~\cite{tibshirani1996regression})
regression (LS). For completeness, Appendices \ref{subsection:randomForest} and \ref{subsection:lasso} provide brief
summaries of these two techniques.

\subsection{Application of error models} 
\label{subsection:application}
We propose two practical ways to use a QoI-error prediction
$\errorOperatorQoIApproxTimeN$: (1) as a \textit{correction} to the
surrogate-model QoI prediction at a given time instance, or (2) as an \textit{error indicator} to be used within the
Gaussian-process-based ROMES framework \cite{drohmann2015romes} for
statistically \reviewer{modeling} arbitrary functions of the time-dependent surrogate-model error. These approaches are now described in turn.
\subsubsection{QoI correction}
The most obvious way in which to employ the QoI-error prediction
$\errorOperatorQoIApproxTimeN$ is simply to apply it as a correction to the
time-instantaneous QoI computed using the surrogate model, i.e., employ
 \begin{equation} 
\corrOperatorTwo^{\timestep}(\paramsGen)\defeq
\romOperatorTwo^{\timestep}(\paramsGen) + \errorOperatorQoIApproxTimeN(\paramsGen),
  \end{equation} 
	as a corrected QoI at time instance $\timestep$. Of course, our expectation is that the corrected QoI error
	has smaller magnitude than the surrogate QoI error, i.e., 
	$
	|\fomOperatorTwo^{\timestep}(\paramsGen) - 
	\corrOperatorTwo^{\timestep}(\paramsGen) | = |
\errorOperatorQoITimeN(\paramsGen) - \errorOperatorQoIApproxTimeN(\paramsGen)
	|
 < 
	|\fomOperatorTwo^{\timestep}(\paramsGen) - 
	\romOperatorTwo^{\timestep}(\paramsGen)|  =
	|\errorOperatorQoITimeN(\paramsGen)|
	$.
\subsubsection{QoI error modeling}\label{sec:EMML_ROMES}
%\ST{Add a note here that in general, machine learning algorithms can be used
%to equip the error estimate predictions at each time step with confidence
%interval and prediction interval, which in turn can be interpreted as a source
%of epistemic uncertainty. For example, ordinary least squares method
%predictions are accompanied with confidence interval and prediction interval
%derived from normal theory assumptions. Similarly, machine learning algorithms
%like random forest can be used to predict confidence intervals using
%bootstrapping technique []. However, such an approach will be explored in
%future work. } \KTC{I'm not sure that's relevant in this section. I prefer to
%leave this out unless asked by a reviewer. The reason is that I've focused on
%functions of the error instead of the error itself, so it wouldn't be obvious
%how to do the above anyway.}\\

%
Alternatively, we can adopt the perspective of the ROMES method
\cite{drohmann2015romes}, which aims to
construct a \textit{statistical model} of the surrogate-model error via
Gaussian-process regression. The key insight of the method is that one
particular type of surrogate---reduced-order models---produce inexpensively
computable \textit{error indicators} such as error bounds, residual norms, and
dual-weighted residuals that correlate strongly with the surrogate-model
error. The method exploits such error indicators by constructing a Gaussian
process that maps the chosen error indicator to a (Gaussian) distribution over
the true surrogate-model error.

In the present context, given an arbitrary function of the time-dependent
surrogate-model error that we would like to predict,
$\errorFun(\errorOperatorQoIApprox^{1}(\paramsGen),\ldots,\errorOperatorQoIApprox^{\numTimesteps}(\paramsGen))$,
we propose employing the 
same function applied to the EMML-predicted surrogate QoI errors,
$\errorFun(\errorOperatorQoIArg{1}(\paramsGen),\ldots,\errorOperatorQoIArg{\numTimesteps}(\paramsGen))$,
as an error indicator in the ROMES framework. We expect this to perform well if
the 
EMML QoI-error predictions
$\errorOperatorQoIApproxTimeN$, $n=1,\ldots,\numTimesteps$, are accurate
representations of the true QoI errors
$\errorOperatorQoITimeN$, $n=1,\ldots,\numTimesteps$. 

Specifically, defining a probability space $(\Omega, \mathcal F, P)$, we
approximate the deterministic mapping 
$\paramsGen\mapsto\errorFun(\errorOperatorQoIArg{1}(\paramsGen),\ldots,\errorOperatorQoIArg{\numTimesteps}(\paramsGen))$
by a stochastic mapping 
$\errorFun(\errorOperatorQoIApprox^{1}(\paramsGen),\ldots,\errorOperatorQoIApprox^{\numTimesteps}(\paramsGen))\mapsto
\hat\errorFun$, with $\hat\errorFun:
\Omega\rightarrow \RR{}$ a scalar-valued Gaussian random variable that can be
considered a statistical model of the true error function
$\errorFun(\errorOperatorQoIArg{1}(\paramsGen),\ldots,\errorOperatorQoIArg{\numTimesteps}(\paramsGen))$. The
stochastic mapping is constructed via Gaussian-process regression using
ROMES training data
$$
	\{ \left( 
	\errorFun(\errorOperatorQoIArg{1}(\paramsGen),\ldots,\errorOperatorQoIArg{\numTimesteps}(\paramsGen)),
\errorFun(\errorOperatorQoIApprox^{1}(\paramsGen),\ldots,\errorOperatorQoIApprox^{\numTimesteps}(\paramsGen))\right)
	\}_{\paramsGen\in\paramTrainROMES},
$$
where $\paramTrainROMES\subset\paramSpace$ denotes the ROMES training set,
which should be distinct from the EMML training set $\paramTrain$. 

%----------------------------------------------
\section {Application to subsurface flow}
\label{section:governingEquations}
%----------------------------------------------

In this section, we first present the governing equations for a two-phase
oil--water subsurface flow system, followed by the POD--TPWL reduced-order (surrogate) model used in this work. We then describe the specialization of the EMML components (error-modeling approach,
feature design, training/test data, \reviewer{determining locality for the regression model}) employed for this application. Please refer to \cite{CardosoTPWL} for a description of the oil--water flow equations and associated finite-volume discretization, and for a detailed development of POD--TPWL for such systems. The use of LSPG projection  with POD--TPWL is described in \cite{he2014reduced,he2015constraint,trehanTPWQ}.

%We refer the reader to \cite{AzizAndSettari} for
%more details on subsurface flow modeling, and to \cite{CardosoTPWL} for a
%detailed development of POD-TPWL.

%--------------------------------
\subsection {High-fidelity model: two-phase oil--water system}\label{sec:HFMOilWater}
%--------------------------------
The HFM for the two-phase oil--water problem entails statements of conservation of mass for the oil and water components combined with Darcy's law for each phase. Assuming 
immiscibility (which means that components only exist in their corresponding phase), and neglecting capillary pressure and gravitational effects,
the equations for phase $j$ can be written as
\begin{subequations}
\label{eqn:governingEquationTogether}
\begin{align}
	\label{eqn:MassConservation}
	\frac{\partial} {\partial t} \left( {\phi\phaseDensity{j} S_{j}} \right) + \nabla \cdot \left(\phaseDensity{j} \phaseVelocity{j} \right) + \phaseDensity{j} {\tilde q}_{j} = 0,\\
	\label{eqn:Darcy}
	\phaseVelocity{j}=  -\lambda_{j} {\perm} \nabla p,
\end{align}
\end{subequations}
where $j$=$o$ designates the oil phase and $j$=$w$ the water phase, $t$ is time, $\phi$
is porosity (void fraction within the rock), $\phaseDensity{j}$
denotes phase density, $S_{j}$ is phase saturation (i.e., phase fraction), $\phaseVelocity{j}$
is the Darcy phase velocity,  ${\tilde q}_{j}$ denotes the well phase flow rate per unit volume
(${\tilde q}_{j}>0$ for production/sink wells and ${\tilde q}_{j}<0$ for
injection/source wells), ${\perm}$ denotes the permeability tensor (taken to be isotropic in the examples here), $\lambda_{j}= \relperm(S_j) / \mu_{j}$ designates the phase mobility, with $\relperm$ the relative
permeability to phase $j$ and $\mu_j$ the phase viscosity, and $p$ denotes pressure (note
that $p_{o} = p_{w} = p$ because capillary pressure is neglected). 
We additionally have the saturation constraint $S_{w}+S_{o} = 1$. 
For subsurface flow problems, we are often interested in predicting the phase flow rates for all production and injection wells. 

The two-phase system described by \eqnRefOne{eqn:governingEquationTogether} is
discretized using a finite-volume method with pressure $p$ and water saturation
$S\defeq S_w$ in each grid block as the primary unknowns. Thus, $\nSampleOutputState=2$. Then, the
(time-dependent) states can be represented as
 \begin{equation} 
\x = \left[p_1\ S_1\ \cdots\ p_{\nCells}\ S_{\nCells}\right]^T\in\RR{2\nCells},
 \end{equation} 
where %a subscript $d$ denotes the value of a variable at grid block $d$,  
$\nCells$ denotes the number of grid blocks; thus, $\nstateFom=2\nCells$ for this application.
At each time step, we consider the set of QoI to be the well phase flow rates
$q_{j}$, at a subset of $\nWells\ll\nCells$ grid blocks referred to as the well
blocks, which we represent by the indices
$\wellBlockSet\defeq\{\wellBlockArg{1},\ldots,\wellBlockArg{\nWells}\}\subset\{1,\ldots,\nCells\}$ with $\wellBlockSet = \wellBlockSetProducer\cup\wellBlockSetInjector$ and $\wellBlockSetProducer\cap\wellBlockSetInjector=\emptyset$, where 
$\wellBlockSetProducer\subseteq\wellBlockSet$ denotes the set of
producer wells and $\wellBlockSetInjector\subseteq\wellBlockSet$ denotes the set of injector wells.

We compute these flow-rate QoI using the standard Peaceman well model~\cite{peaceman1983}:
\begin{align}
\label{eqn:PeacemanModel}
\phaseFlowRateWellD=	 	 (T_{\rm{well}})_d \left( \lambda_{j} \right)_{d} \left( p_{d}- u_d \right),\quad j=o,w\ \text{for}\ d\in\wellBlockSetProducer,\quad
j=w\ \text{for}\ d\in\wellBlockSetInjector
 . %\quad d=1,\ldots,\nWells.
\end{align}
Here, subscript $d$ indicates the value of a variable at grid block $d$, $T_{\rm{well}}\in\RR{}_+$ denotes the well index, which depends on the wellbore radius and well-block permeability and geometry ($T_{\rm{well}}$ is essentially the transmissibility linking the well to the well block), and $u\in\RR{}_+$ denotes the
prescribed wellbore pressure, also referred to as the bottom-hole pressure
(BHP). \eqnRefOne{eqn:PeacemanModel} is written for a discrete finite-volume model. Thus, $\phaseFlowRateWellD={\tilde q}_j V_d$, where $V_d$ is the volume of grid block $d$ $\left( \phaseFlowRateWellD \text{ is of units volume/time} \right)$. 

In the systems considered here, only water is injected. The output-function $\fomOperatorTwo$ introduced in
\eqnRefOne{eq:fomOutput} is defined by
\eqnRefOne{eqn:PeacemanModel}, 
%\footnote{Note that $(\lambda_{j})_d$ depends on $S_d$ and $(\phaseDensity{j})_d$ depends on $p_d$.}, 
where the sampling matrix
$\sampleOutputState$ simply extracts the pressure and saturation from the
appropriate well block, i.e., the output corresponding to $(q_j)_d$, $j=o,w$,
employs a sampling matrix $\sampleOutputState_d = \left[\unitVec{2d-1}\
\unitVec{2d}\right]^T$, where $\unitVec{i}$ denotes the $i$th canonical unit
vector. Thus, \reviewer{the QoI depend on both pressure and saturation,
i.e.,} $\nSampleOutputState = 2$ for each QoI. Note that the treatment described here is for cases where a particular well penetrates only a single grid block. Our procedures can be generalized, however, for multiblock well penetrations and for cases where rates (rather than BHPs) are specified.

In this work, we employ the time-varying well BHPs as the control
variable. We denote the (time-dependent) control vector  as
 \begin{equation} 
\uVar = \left[u_{\wellBlockArg{1}}\ \cdots\
u_{\wellBlockArg{\nWells}}\right]^T\in\RR{\nWells}.
 \end{equation} 
The time-varying BHP profiles constitute the input
parameters, i.e.,
 \begin{equation}
 \label{eqn:BHPschedule} 
\paramsGen =
\left[
\begin{array}{c}
\uArg{1}\\
\vdots\\
\uArg{\numTimesteps}
\end{array}
\right]\in\RR{\numParam},
 \end{equation} 
thus, $\numParam = \nWells\numTimesteps$. Alternatively, well flow rates $(q_j)_d$, $d\in\wellBlockSet$, or a combination of rates and BHPs, could be prescribed as the control variables.

Following ~\cite{CardosoTPWL,he2011enhanced,he2014reduced}, the discretized
set of nonlinear algebraic equations (obtained using fully implicit
discretization\footnote{Note that this assumes that a linear multistep method with $k=1$ step is applied for time integration.}) describing the
HFM is represented as
%
% Expression for residual
\begin{align}
\label{eqn:discretizedPDE}
			\gNPlusOne \defeq \textbf{g} \left( \xN, \xNMinusOne, \uNPlusOne \right) = \boldZero,\quad n=1,\ldots,\numTimesteps,
\end{align} 
where $\textbf{g}:(\gArgOne,\gArgTwo,\gArgThree)\mapsto \textbf{g}(\gArgOne,\gArgTwo,\gArgThree)$ and $\textbf{g}: \RR{2\nCells}\times \RR{2\nCells}\times\RR{\nWells}\rightarrow \RR{2\nCells}$ 
designates the residual vector we seek to drive to zero, and superscript $n$
denotes the value of a variable at time step $n$. The state
operator $\fomOperatorOneTimeN$ defined in \eqnRefOne{eqn:FOMGeneralEquations}
is implicitly defined by the sequential solution to
\eqnRefOne{eqn:discretizedPDE}. 

%--------------------------------
\subsection {Surrogate model: POD--TPWL}\label{sec:surrPODTPWL}
%--------------------------------
We now briefly describe the POD--TPWL formulation for oil--water systems, which
will be our surrogate model for this application. For further details, the
reader is referred to \cite{CardosoTPWL,he2011enhanced,he2014reduced,he2015constraint,trehanTPWQ}. Given a set of `test'
controls $\uNPlusOne$, $n=1,\ldots,\numTimesteps$, the POD--TPWL model 
linearizes the residual around a previously saved `training' simulation
solution. Then, at time step $n$ in the test simulation, rather than solving
the system of nonlinear algebraic equations \eqref{eqn:discretizedPDE} using, e.g., Newton's method, we
instead solve the system of linear algebraic equations 
%
%\begin{align}
	%%\label{eqn:DiscretizedLinear}
	%\nonumber \gNPlusOne_{\text{L}} &\defeq \gIPlusOne + \frac{\partial \gIPlusOne} {\partial \xIPlusOne}(\xNPlusOne - \xIPlusOne) + \frac{\partial \gIPlusOne} {\partial \xI}(\xN - \xI) + \frac{\partial \gIPlusOne} {\partial \uIPlusOne} (\uNPlusOne - \uIPlusOne) = \boldZero,
%\end{align}
\begin{align}
	\label{eqn:DiscretizedLinear}
	\gNPlusOne_{\text{L}} &\defeq \JIPlusOne(\xN - \xI) +
	\BIPlusOne(\xNMinusOne - \xIMinusOne) + \CIPlusOne (\uNPlusOne - \uIPlusOne) =
	\boldZero,\quad n=1,\ldots,\numTimesteps,
\end{align}
where we have used the fact that $\g(\xI,\xIMinusOne,\uIPlusOne)=\boldZero$ and defined
\begin{align}
	% TPWL
	\nonumber &\JIPlusOne  \defeq  \left. \frac{\partial \g} {\partial\gArgOne}\right\rvert_{(\xI,\uIPlusOne)}\in \RR{2\nCells\times2\nCells}, \quad 
	\nonumber \BIPlusOne  \defeq \left.  \frac{\partial \g} {\partial\gArgTwo}\right\rvert_{\xIMinusOne}\in \RR{2\nCells\times2\nCells},\\
	\nonumber &\CIPlusOne \defeq  \left. \frac{\partial \g} {\partial\gArgThree}\right\rvert_{(\xI,\uIPlusOne)}\in \RR{2\nCells\times\nWells}.
\end{align}
%\KTC{partial subscripts should be with respect to a variable, not an instance of the variable.}
%
Here, $\xTraining$ denotes the (saved) training simulation state, $\uTraining$
denotes the training controls, and $i$ denotes the time step associated with the
training state about which we linearize. Note that an accent $\acute\cdot$ indicates that the quantity
has been saved during training simulations. While the criterion for
determining the `closest' training configuration $(\xI,\xIMinusOne,\uIPlusOne)$ about which to linearize is application
dependent, here we use pore volume injected (PVI) to determine the appropriate
training solution. PVI quantifies the fraction of the system pore space that
has been filled by injected fluid (water in our case), and as such corresponds
to a dimensionless time. Thus, we seek to linearize around a solution that has
progressed to the same PVI as the current test solution. See
\cite{he2014reduced} for further discussion and details on the computation of
PVI.

To reduce the computational cost associated with solving
\eqnRefOne{eqn:DiscretizedLinear}, we approximate the state $\x$ in a
low-dimensional affine subspace, using POD, as
 ${\x}\approx{\PODbasis} {\z} + \PODMean$, where $\PODbasis \in
\RR{2\nCells \times \nstateRom}$ denotes a POD basis, $\z\in \RR{\nstateRom}$ designates the reduced
state, and
$\PODMean\in\RR{2\nCells}$ indicates a reference state, which is often taken to
be the mean of the snapshots. 
Replacing $\x$ with ${\PODbasis} {\z} + \PODMean$ in 
\eqnRefOne{eqn:DiscretizedLinear} yields
\begin{align}
\label{eqn:DiscretizedLinearReducedIntermediate}
		 \gNPlusOneTilde_{\text{L}} &= \JIPlusOne {\PODbasis}(\zN
		 -\zI) + \BIPlusOne {\PODbasis}(\zNMinusOne -\zIMinusOne)+ \CIPlusOne(\uNPlusOne - \uIPlusOne) = \boldZero,
\end{align}
where $\zIGen \defeq \PODbasis^T(\xIGen - \PODMean)$.

Because \eqnRefOne{eqn:DiscretizedLinearReducedIntermediate} is overdetermined ($2\nCells$
equations and $\nstateRom<2\nCells$ unknowns), it may not have a solution. Thus, we reduce the number of equations to $\nstateRom$ by forcing the residual in
\eqnRefOne{eqn:DiscretizedLinearReducedIntermediate} to be orthogonal to the
range of a test basis ${\LHSbasis} \in \RR{2\nCells \times \nstateRom}$. In
line with previous studies on the application of POD--TPWL for subsurface flow
models~\cite{he2014reduced, he2015constraint, trehanTPWQ}, we employ the least-squares
Petrov--Galerkin (LSPG) test
basis~\cite{CarlbergPG,carlberg2013gnat,carlbergGalDiscOpt}, i.e.,
$\LHSbasisIPlusOne = \JIPlusOne{\PODbasis}$. Premultiplying
\eqnRefOne{eqn:DiscretizedLinearReducedIntermediate} by
$(\LHSbasisIPlusOne)^T$, the linear system of equations in the low-dimensional
space is now expressed as
\begin{align}
	\label{eqn:DiscretizedLinearReduced}
	  \gRL = 	\JrIPlusOne (\zN -\zI) + \BrIPlusOne(\zNMinusOne -\zIMinusOne)+\CrIPlusOne(\uNPlusOne - \uIPlusOne) = \boldZero, 
\end{align}
where
\begin{subequations}
\label{eqn:PODTPWLcoeff}
\begin{align}
 	%   
	%\label{eqn:JrIPlusOne}
	&\JrIPlusOne \defeq \LHSbasisIPlusOneT \JIPlusOne {\PODbasis}\in \RR{\nstateRom\times \nstateRom }, \quad
    %
	%\label{eqn:BrIPlusOne}
	\BrIPlusOne  \defeq   \LHSbasisIPlusOneT \BIPlusOne {\PODbasis}\in \RR{\nstateRom \times \nstateRom },\\
    %
	%\label{eqn:CrIPlusOne}
	&\CrIPlusOne \defeq  \LHSbasisIPlusOneT \CIPlusOne\in \RR{\nstateRom \times
	\nWells},
	%\\ 
	%&\JrIPlusOne \in \RR{\nstateRom\times \nstateRom }, \quad 
	%\BrIPlusOne \in \RR{\nstateRom \times \nstateRom }, \quad 
	%o\CrIPlusOne \in \RR{\nstateRom \times \numParam}. 
\end{align}
\end{subequations}
and
the subscript $\RL$ indicates that this is the POD--TPWL representation.
Thus, the surrogate-state
operator $\romOperatorOneTimeN$ defined in \eqnRefOne{eqn:ROMGeneralEquations}
is implicitly defined by the sequential solution to
\eqnRefOne{eqn:DiscretizedLinearReduced} for this application. Further, this implies
 $\sampleRomOutputState = \sampleOutputState\PODbasis$. We also define the
 prolongation operator associated with the output in well block $d$ as
 $\sampleRomOutputState_d \defeq \sampleOutputState_d\PODbasis$
 for this system. During online POD--TPWL computations, well-block saturations computed as $\sampleRomOutputState_d\zN$ can fall outside of the physical range. In this case, any $S<0$ is mapped to $S=0$, and $S>1$ is mapped to $S=1$. 

POD--TPWL requires some number of training runs to be performed during an offline
stage. These training simulations involve solving
\eqnRefOne{eqn:discretizedPDE} for prescribed time-varying BHPs
$\paramsGen\in\paramTrainTPWL\subseteq\paramSpace$, where $\paramTrainTPWL$ denotes the
TPWL training points. The state snapshots generated from the
training simulations are saved and used to construct the POD basis $\PODbasis$
by performing SVD on the (centered) snapshots. The constructed POD basis is
then used to perform offline processing, which involves computing and saving
quantities such as $\JrIPlusOne, \BrIPlusOne$ and $\CrIPlusOne$. Although several training runs are used to construct the POD basis, only one of these, referred to as the primary training run, is used for linearization (i.e., $\JIPlusOne$, $\BIPlusOne$ and $\CIPlusOne$ all derive from the primary training run). 
%
%
%------------------------------------
\subsection{Error modeling}
\label{subsection:Features}
%------------------------------------
To \reviewer{model} QoI errors we adopt two approaches. 
Approach 1 is a hybrid treatment wherein we apply 
Method 4 (Section \ref{sec:method4}) to \reviewer{model} \textit{relative errors} in the well-block
pressure and
Method 3
(Section \ref{sec:method3}) to \reviewer{model} \textit{errors} in the well-block saturation.
 We pursue this approach because the error in the well-block pressure can exhibit
a wide range of values, which makes modeling the relative error an easier
task. On
the other hand, the error in the well-block saturation spans a narrow
range of absolute values because $0 \le S\le 1$; thus, directly modeling the state
error is appropriate.
Approach 2 simply applies 
Method~2 (Section \ref{subsubsection:Method2}) to model the relative QoI error directly.

%In the first approach, we apply two of the methods
%discussed earlier. We estimate the error in sampled states
%$\errorOperatorPrimaryVarApproxTimeN, n=1,\ldots,\numTimesteps,\
% g = {p_d,S_d},\ d\in\wellBlockSet
% $, and then use
%\eqnRefOne{eq:methodThreeFinalEnd} to compute $\errorOperatorQoIApprox$. The
%well-block pressure and saturation provided by POD-TPWL are first
%reconstructed from the reduced state using $[\pNWellBlockTPWL,
%\sNWellBlockTPWL] = \sampleRomOutputState \zN$. Due to scaling issues, the
%error in the well-block pressure is modeled as a relative error, i.e., 
%%
%\begin{align}
%	\errorPN \defeq \dfrac{\pNWellBlockHF - \pNWellBlockTPWL}{\pNWellBlockHF},
%\end{align}
%%
%while error in the well-block saturation is treated in terms of the actual error, i.e., 
%%
%\begin{align}
%	\errorSN \defeq \sNWellBlockHF - \sNWellBlockTPWL, 
%\end{align}
%%
%since $0 \le S_d \le 1$. Thus, this is a hybrid approach involving Methods 3
%and 4. This information is then used to compute
%$\errorOperatorQoIApproxTimeN$. We will also consider the direct application
%of Method~2, where we calculate the $\errorOperatorQoIApproxTimeN$ from
%\eqnRefOne{eq:relToAbsModel}.
%

\subsection{Feature design}\label{sec:feature_design}
The definition of the QoI in well block $d$ \eqref{eqn:PeacemanModel}, 
the HFM governing equations \eqref{eqn:discretizedPDE}, and the
surrogate-model governing equations \eqref{eqn:DiscretizedLinearReduced}
suggest that the corresponding sampled-state error $\errorOperatorPrimaryVarTimeN$ and 
QoI error $\errorOperatorQoITimeN$ will
likely depend on data such as the states/controls about which POD--TPWL is
linearized $(\sampleRomOutputState_d\zIPlusOne,\sampleRomOutputState_d\zIMinusOne,\uIPlusOne)$, the current well-block state
$\sampleRomOutputState_d\zN$, the previous well-block state
$\sampleRomOutputState_d\zNMinusOne$, and operators associated with the linearized
system such as 
$\sampleOutputState_d\JI \sampleOutputState_d^T\in \RR{2\times2}$,
$\sampleOutputState_d\BI\sampleOutputState_d^T\in \RR{2\times2}$, and
$\sampleOutputState_d\CI\sampleOutputState_d^T\in \RR{2\times2}$, where the $2\times2$ matrix is converted to a $1\times4$ row vector. Similarly,
these errors may also depend on
other quantities such as the PVI in the primary training (${\rm{PVI}}^{i}$) and the test case
(${\rm{PVI}}^{\timestep}$), as well as data such as time-step sizes $
\dt^{\timestep}$ and $\Delta \acute{t}^{i}$, and time instances
$t^{\timestep}$ and $\acute t^{i}$. Some of the features used in this work are shown in
Table~\ref{table:features}. \reviewer{We note that most of these features are
already computed during the course of the ROM simulation (e.g., the control input
$\BHPnoBf^{\timestep}_{d}$), while others (e.g., $ \dfrac{\innerProduct{\zN}{\zI }}
			{\norm{\zN}{2}~\norm{\zI}{2}}$) can be computed via inexpensive
			computations, with cost scaling with 
			the ROM dimension $\nstateRom$. Further, while Hessian information is
			useful for informing the linearization error of TPWL, we do not include
			it in the feature set due to the high 
			cost of its computation \cite{trehanTPWQ}.}

\begin{table} []
\centering
\caption{A subset of the features $\featuresN$ used in EMML} % title of Table
\label{table:features}
\begin{tabular}{l l l l l l}
\toprule \hline 
	     {\bf No. }  &
	     {\bf Feature} & 
	     {\bf No.}  & 
	     {\bf Feature}  \\[0.1cm]  
	    \hline
	    1. & $\sampleRomOutputState_d\zN$ &
	    2.  & $\sampleRomOutputState_d\zI$\\
	    3. & $\sampleRomOutputState_d\zNMinusOne$ &
	    4. & $\sampleRomOutputState_d\zIMinusOne$ \\
	    5.  & $\dSdTn$  &
	    6.  & $\dSdTi$   \\
	    7.  & $\BHPnoBf^{\timestep}_{d}$ &
	    8. & $\acute{\BHPnoBf}^{i}_{d}$ \\
	    9. & $\sampleOutputState_d\JI \sampleOutputState_d^T $ & 
	    10. & $\sampleOutputState_d\BI\sampleOutputState_d^T $ \\ [1mm]
	    11. & $\sampleOutputState_d\CI\sampleOutputState_d^T $ &
	    12. & $ \dfrac{\innerProduct{\zN}{\zI }}
			{\norm{\zN}{2}~\norm{\zI}{2}}$ \\
     	13. & $ \dt^{\timestep}$ & 
     	14. & $\Delta \acute{t}^{i}$  \\
     	15. & ${\rm{PVI}}^{\timestep}$&
	    16. & ${\rm{PVI}}^{i} $  \\
  		17. & $t^{\timestep}$ &
  		18. & $\acute t^{i}$ \\ 
  		19. & $\sampleRomOutputState_{2d}\zN - \sampleRomOutputState_{2d}\zI$ &
  		20. & $\sampleRomOutputState_{2d}\zNMinusOne - \sampleRomOutputState_{2d}\zIMinusOne$ \\
		\bottomrule
 \end{tabular} 
\end{table}

Note that solving \eqnRefOne{eqn:DiscretizedLinearReduced} for time step
$\timestep$ requires information from time step $\timestep-1$; however, we could also
include data from multiple previous time steps in the set of features. To include
features generated over a `memory' of $\memory$ previous time steps, we define
a new feature vector $\featuresMemN(\memory)\defeq \left [\featuresN\
\features^{\timestep-1}\ \cdots\ \features^{\timestep-\memory} \right]\in\RR{1
\times \numFeaturesMemory}$. Note that some features in
$\featuresMemN(\memory)$ will be strongly correlated (some will be identical);
we remove such highly correlated features in a preprocessing step by
computing the feature--feature Pearson correlation coefficients for all pairs.
Further details on correlation-criteria-based feature selection can be found
in~\cite{guyon2003introduction}. 
%
%
%------------------------------------
\subsection{EMML training and test data}
\label{subsection:TrainingTestDataset:EMML}
%------------------------------------
Following~\cite{trehanTPWQ}, we generate a set of $\numTestSamples$ BHP controls
by adding unique random perturbations to the `primary' training BHP control
$\paramsLinearize\equiv\left[(\uLinArg{1})^T\ \cdots\
(\uLinArg{\numTimesteps})^T\right]^T\in\paramTrainTPWL$
used in the primary training run.  For a given set of BHP controls, we define the perturbation in producer BHPs as 
\begin{equation}
	\Delta \BHPnoBf^{P}(\paramsGen) = \perturbationProducer.
\end{equation}
The perturbation in injection BHPs $\Delta \BHPnoBf^{I}(\paramsGen)$ is defined by
instead performing summation over the set of injector wells
$\wellBlockSetInjector$.

We partition the $\numTestSamples$ BHP schedules into $\numHFSsymbol<\numTestSamples$ clusters according to
their representation in the (two-dimensional) space defined by 
$\Delta \BHPnoBf^{P}$ and $\Delta \BHPnoBf^{I}$. We select the BHP schedules closest to the
cluster centers as `representative' schedules for which we  simulate both the
POD--TPWL and the high-fidelity models.
This subset of controls
$\paramTrain\defeq\{\paramsTrainArg{1},\ldots,\paramsTrainArg{\numHFSsymbol}\}
\subseteq\paramSpace$ constitutes the EMML training data. The
remaining set of $\numTestSamples-\numHFSsymbol$ test schedules comprises the EMML test data
$\paramTest\defeq\{\paramsTestArg{1},\ldots,\paramsTestArg{\numTestSamples-\numHFSsymbol}\}
\subseteq\paramSpace$. We note that the test set is not used in the
construction of the EMML model; it is simply used to assess EMML
performance after the model has been constructed.
%-------------------------------------------
\subsection{\reviewer{Regression-model locality} for POD--TPWL}
\label{subsection:spacePartitioning:TPWL}
%-------------------------------------------
As described in Section~\ref{subsection:spacePartitioning}, we
\reviewer{determine regression-model locality} using classification and clustering to construct
tailored local regression models for error prediction.  Recall that QoI,
described in~\eqnRefOne{eqn:PeacemanModel}, correspond to the oil and water
flow rates at the production wells $\wellBlockSetProducer$ and water flow rates
at the injection wells $\wellBlockSetInjector$. 
We \reviewer{construct a local regression model} only for production-well QoI
$(q_{j})_d$, $j=o,w$, $d\in\wellBlockSetProducer$; in fact, we
\reviewer{determine regression-model locality} 
for each production well \reviewer{independently}, which is valid for both the oil- and water-flow rate
QoI at that well.  We do not \reviewer{construct local regression models} for injection-well
QoI $(q_{w})_d$, $d\in\wellBlockSetInjector$, as a global regression model
performs well in this case due to the relatively simple behavior of the
associated QoI errors. 
%
%-----------
\subsubsection{Classification}
\begin{figure}[H] \centering
		\includegraphics[page= 3, width=0.55 \textwidth]	{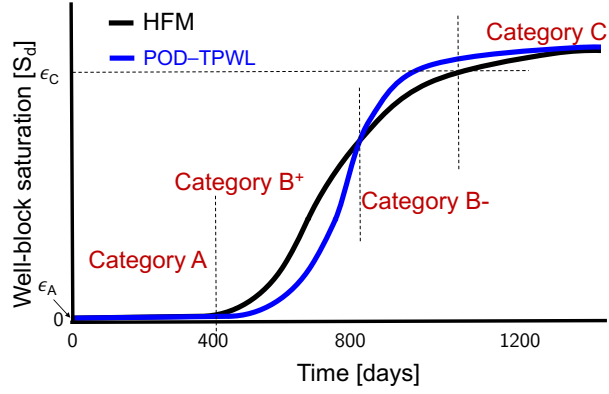}	       
\caption{Category assignment for EMML training data for a producer well.}
\label{fig:labels}
\end{figure}

%-----------
 For production-well QoI, we partition the EMML
	training data into four categories. This partitioning is
	based on well-block saturation $S_d$, $d\in\wellBlockSetProducer$, as shown in \figRefOne{fig:labels}, where the
	blue curve represents the POD--TPWL prediction and the black curve 
	the corresponding HFM prediction. 
%
% Labels
The four categories correspond to different stages of the system and are referred to as $A, B^{+}, B^{-},$ and
$C$. All samples with
$\SDrl\le\epsilon_{A}$ and $S_{d}\le\epsilon_{A}$ are assigned to category $A$, where  $\SDrl$ denotes the well-block saturation predicted by POD--TPWL and, as before, $S_d$ denotes well-block saturation from the high-fidelity simulation. Thus, all the samples in category $A$ have close
agreement between the POD--TPWL and HFM solutions; this category corresponds to the state before water `breakthrough' occurs at a particular production well. 
Samples with $\epsilon_{A} <
\SDrl\leq S_{d} \le \epsilon_{C}$ are assigned 
to category $B^{+}$, while samples with
$\epsilon_{A}
\leq S_{d} < \SDrl\le \epsilon_{C}$
are assigned 
to category $B^{-}$.
Finally, samples with $S_{d}>\epsilon_{C}$ are assigned to 
category $C$. This category corresponds to significant water
production. The actual values used in this work for $\epsilon_{A}$ and $\epsilon_{C}$ are given in Section~\ref{section:numericalResults}.

As mentioned in Section \ref{subsection:spacePartitioning}, we perform
classification using classification features $\featuresClassify\in
\RR{1\times\numFeaturesClassify}$. 
These features quantify (1) the 
perturbation in the prescribed control variables $u^k_d,\ k=1,\ldots,n,~d \in
\wellBlockSetProducer$ relative to the primary training BHP controls
$\acute{u}^k_d,\ k=1,\ldots,n,~d \in \wellBlockSetProducer$, and (2)
the differences in well-block pressure for producer--injector pairs,
which in turn may impact the velocity field as indicated
by~\eqnRefOne{eqn:Darcy}.
For a production well located in grid block
$d \in \wellBlockSetProducer$, classification features include 
 quantities such
as the difference between the test BHP schedule and the primary training run
BHP controls $\paramsLinearize$, i.e.,  
$\Big ( \sum\limits_{k=1}^n \left(\varWithSubSup{\BHPnoBf}{k}{d} - \varWithSubSup{\acute{\BHPnoBf}}{k}{d} \right)^2 \Big)^{\sfrac{1}{2}}$, the average well-block pressure difference between all producer-injector pairs,
$\frac{1}{n} \sum\limits_{k=1}^n (\unitVec{}_{2d-1}^{T} - \unitVec{}_{2d'-1}^{T}) \PODbasis \z^k
$, $d' \in \wellBlockSetInjector$, and the average well-block pressure difference between the test case and the primary training simulation, represented by 
$\frac{1}{n}\sum\limits_{k=1}^n ( \unitVec{}_{2d-1}^T \PODbasis ) ( \z^k -
\acute{\z}^k)$. Table \ref{table:featuresClassify} reports some of the classification features employed in the current application.
\begin{table}[H]
\centering
\caption{Classification features $\featuresClassifyN$ corresponding to a
production well in grid block $d \in \wellBlockSetProducer$} 
\label{table:featuresClassify}
\begin{tabular}{l l l l l l}
\toprule \hline 
	     {\bf No. }  &
	     {\bf Feature} &
	     {\bf No.}  & 
	     {\bf Feature}   \\[0.1cm]  
	    \hline 
	    1. & $\Big ( \sum\limits_{k=1}^n \left(\varWithSubSup{\BHPnoBf}{k}{d} - \varWithSubSup{\acute{\BHPnoBf}}{k}{d} \right)^2 \Big)^{\sfrac{1}{2}}$   &	
		2. & $\frac{1}{n}\sum\limits_{k=1}^n \left( (\varWithSubSup{u}{k}{d})_{\RL} - \varWithSubSup{\acute u}{k}{d} \right)$ \\ [1mm]
		3. & $\Big ( \sum\limits_{k=1}^n \left( (\unitVec{}_{2d-1}^{T} - \unitVec{}_{2d'-1}^{T}) \PODbasis \z^k \right)^2 \Big)^{\sfrac{1}{2}}$, $d' \in \wellBlockSetInjector$ &
		4. & $\Big ( \sum\limits_{k=1}^n \left( ( \unitVec{}_{2d-1}^T \PODbasis ) ( \z^k - \acute{\z}^k)   \right)^2 \Big)^{\sfrac{1}{2}}$ \\
		5. & $\frac{1}{n} \sum\limits_{k=1}^n (\unitVec{}_{2d-1}^{T} - \unitVec{}_{2d'-1}^{T}) \PODbasis \z^k$, $d' \in \wellBlockSetInjector$ &
		6. & $\frac{1}{n}\sum\limits_{k=1}^n ( \unitVec{}_{2d-1}^T \PODbasis ) ( \z^k - \acute{\z}^k)$  \\
		\bottomrule
 \end{tabular}
\end{table}
%

%----------------------------------------------
\section{Numerical results}
\label{section:numericalResults}
%----------------------------------------------

\begin{figure}[H]
 \centering
	\includegraphics[width= 0.6\textwidth]{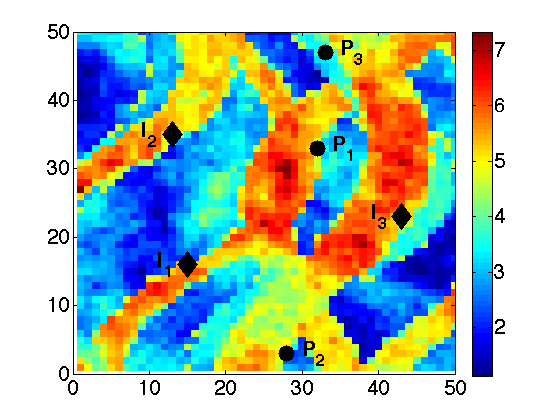}
	\caption {Permeability field ($\log k$, with $k$ in md) and well locations (Model~1,
	from~\cite{ObiThesis}).}
	 \label{fig:permFieldAgain}
\end{figure}

In this section, we present numerical results for the application
described in Section~\ref{section:governingEquations}. The specific
problem involves flow simulation in a synthetic two-dimensional horizontal reservoir.
The reservoir model contains $50\times50$ grid blocks such that $\nCells = 2500$ and $\nstateFom = 5000$. It contains three production wells
$|\wellBlockSetProducer|=3$, which we label as
$P_{1}$, $P_{2}$, and $P_{3}$, and three injection wells
$|\wellBlockSetInjector| = 3$, which we label as $I_{1}$, $I_{2}$, and $I_{3}$. The six wells ($\nWells=6$) are shown in \figRefOne{fig:permFieldAgain}. The permeability field is isotropic, i.e., $\perm=\mathrm{diag}(k)$, and the porosity is set to $\phi=0.2$. The relative permeability
functions are prescribed to be $k_{rw}(S) = S^{2}$ and $k_{ro}(S) = (1 -
S)^{2}$. We apply a backward Euler time integrator with adaptive time-step
selection. 
%
%--- Permeability figue ---

\begin{figure}  [H]
	\begin{subfigure}  {0.5\linewidth} \centering
		\includegraphics[width= 0.75\textwidth]	{{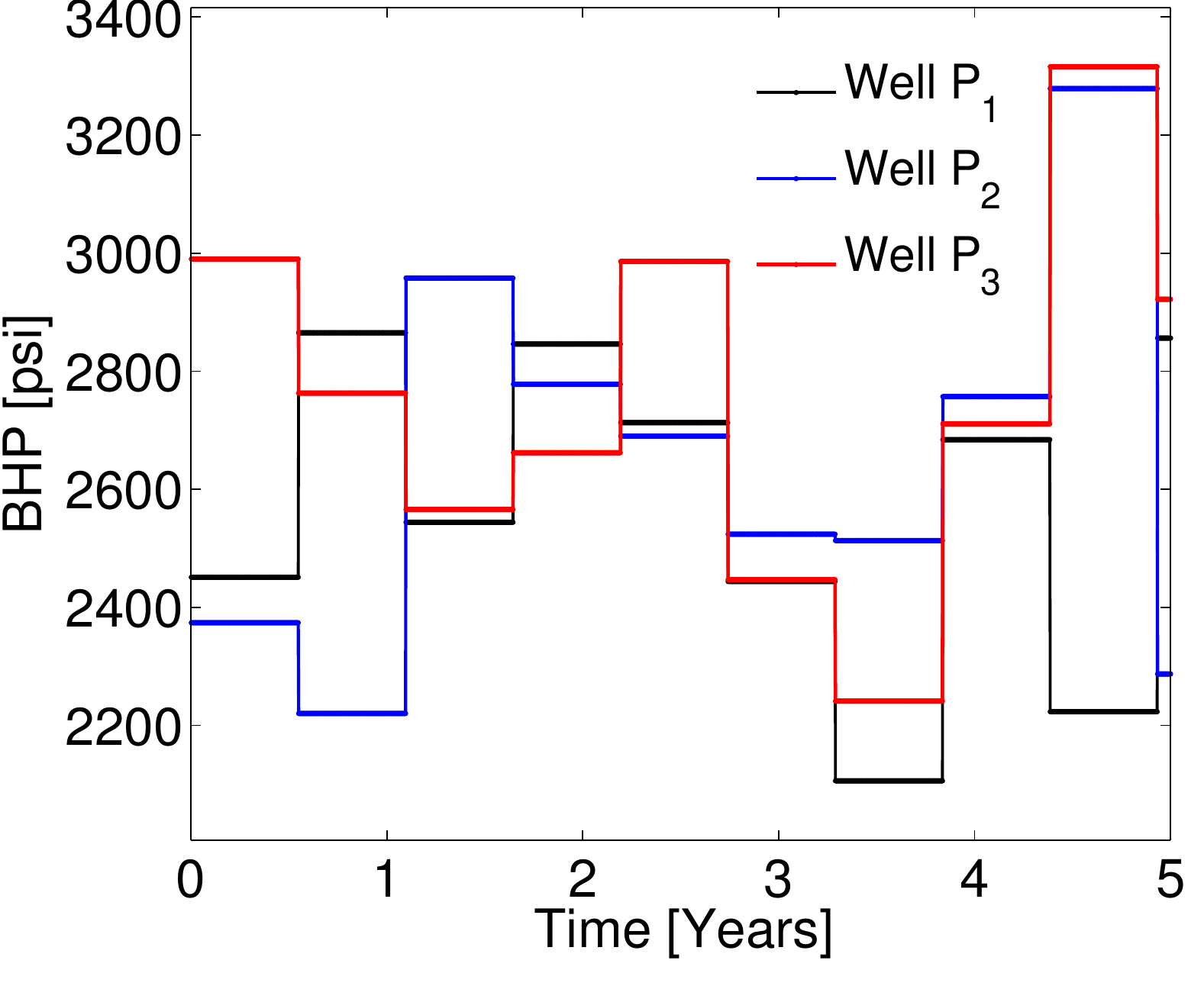}}	       
		\caption{Producer BHP profiles}
	\end{subfigure}
\quad % Injector: Trng 1
	\begin{subfigure}  {0.5\linewidth} \centering
		\includegraphics[width= 0.75\textwidth]	{{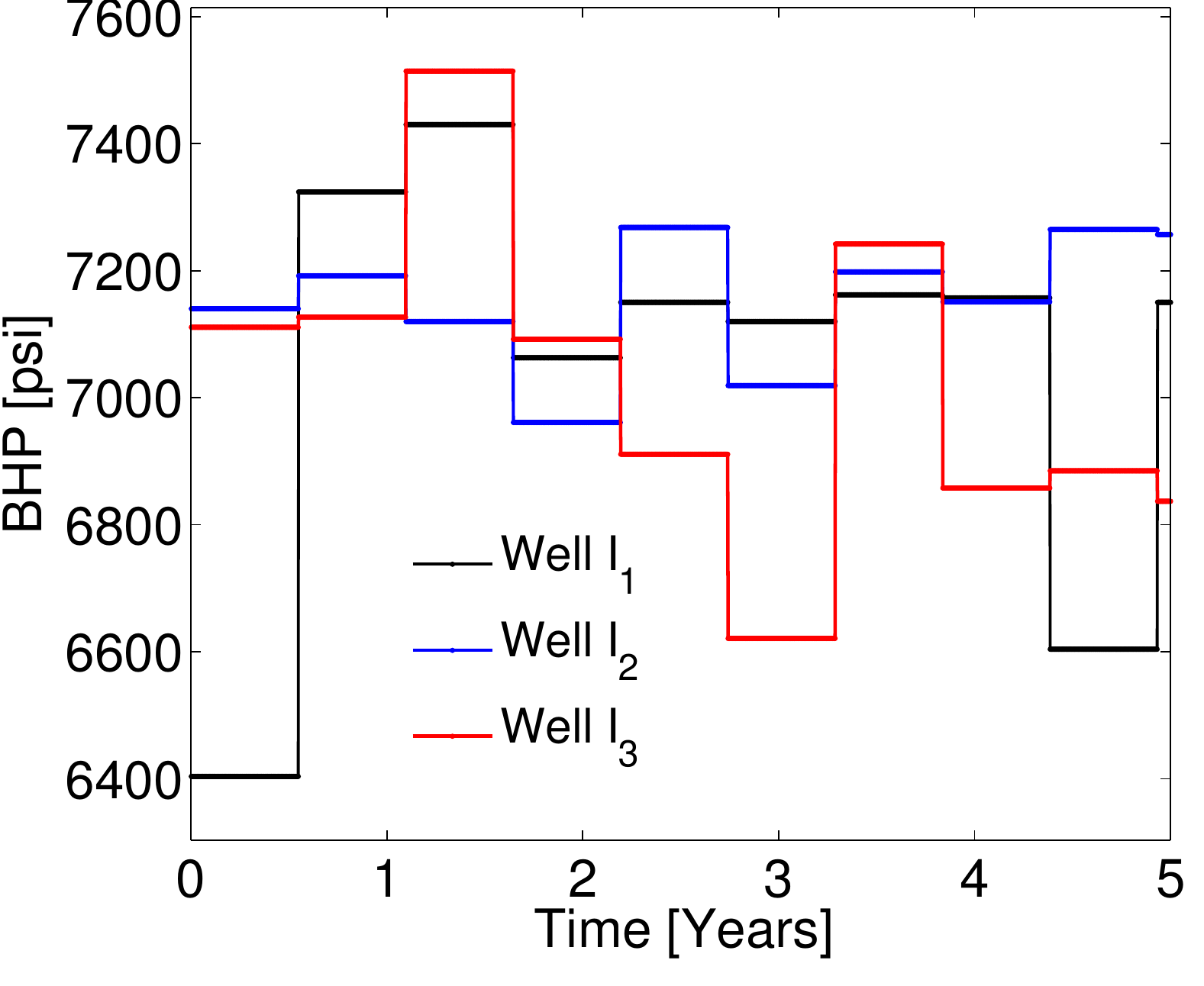}}	       
		\caption{Injector BHP profiles}
	\end{subfigure}
\caption {BHP profiles for the primary training run.}
 \label{fig:case1BHPTraining}
\end{figure}

Three training simulations, $\lvert \paramTrainTPWL \rvert = 3$, are performed
to construct the POD--TPWL model (the three runs provide a sufficient number of snapshots for the POD
basis), from which $\nstateRom=150$ POD basis vectors are extracted. Of these, 90
correspond to the saturation state variables and 60 to the
pressure state variables. \figRefOne{fig:case1BHPTraining} depicts the BHP controls $\acute{\paramsGen}\in\paramTrainTPWL$ applied in the primary training simulation (recall that this is the run used for linearization). These time-varying BHPs, as well as those considered in the test runs, are meant to be representative of the BHP schedules that can arise during oil production optimization computations. In such optimizations, the goal is to determine the time-varying BHPs that maximize an economic metric, or the cumulative oil recovered from the reservoir.

%
% Training BHP profile pictures
%--- Training BHP profile ---
\begin{figure}[H] % Producer: Test case 128
	\begin{subfigure}  {0.50\linewidth} \centering
		\includegraphics[width= 0.75\textwidth]	{{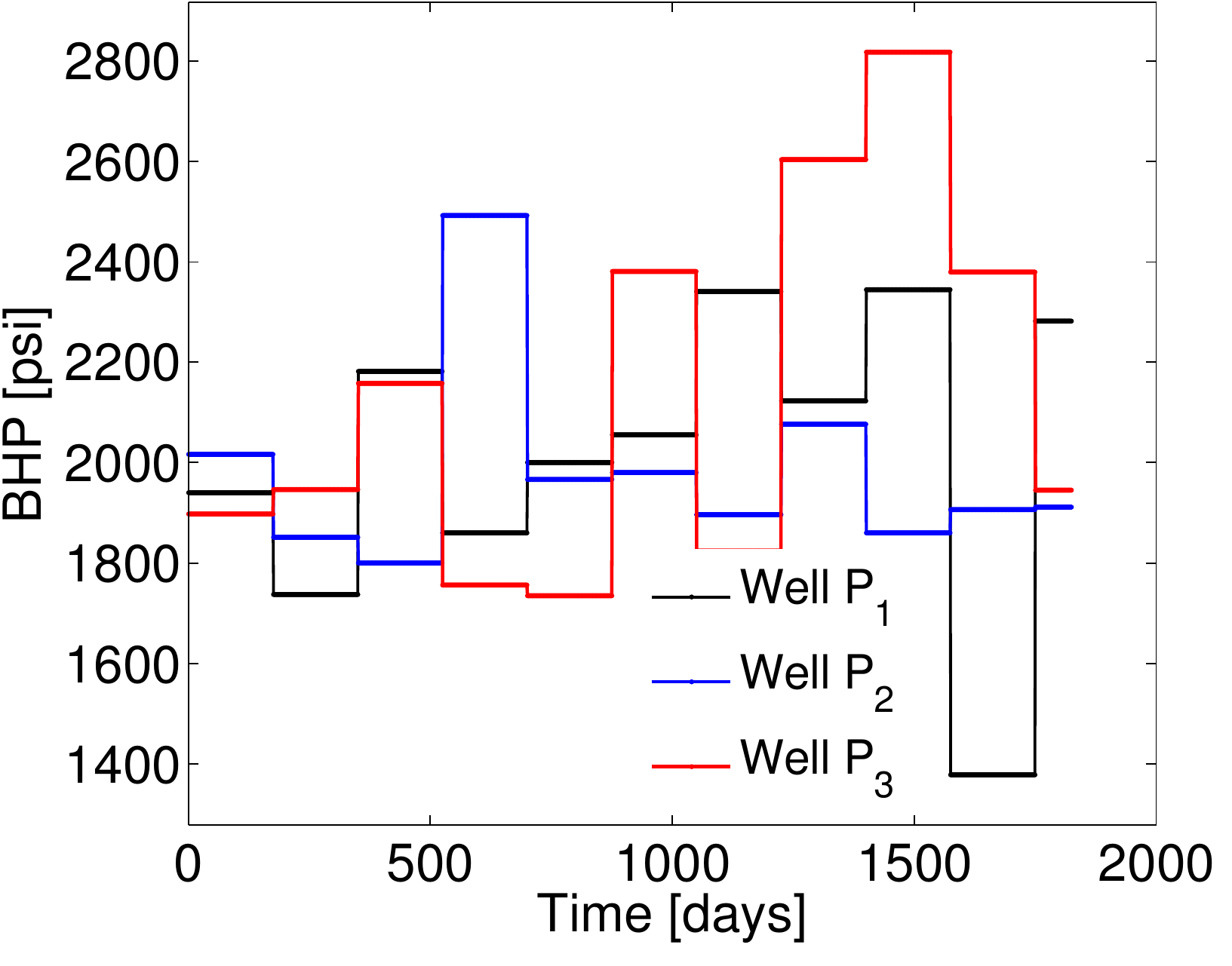}}	       
		\caption{Producer BHP profiles}
	\end{subfigure}
\quad  % Injector: Test case 128
	\begin{subfigure}  {0.50\linewidth} \centering
		\includegraphics[width= 0.75\textwidth]	{{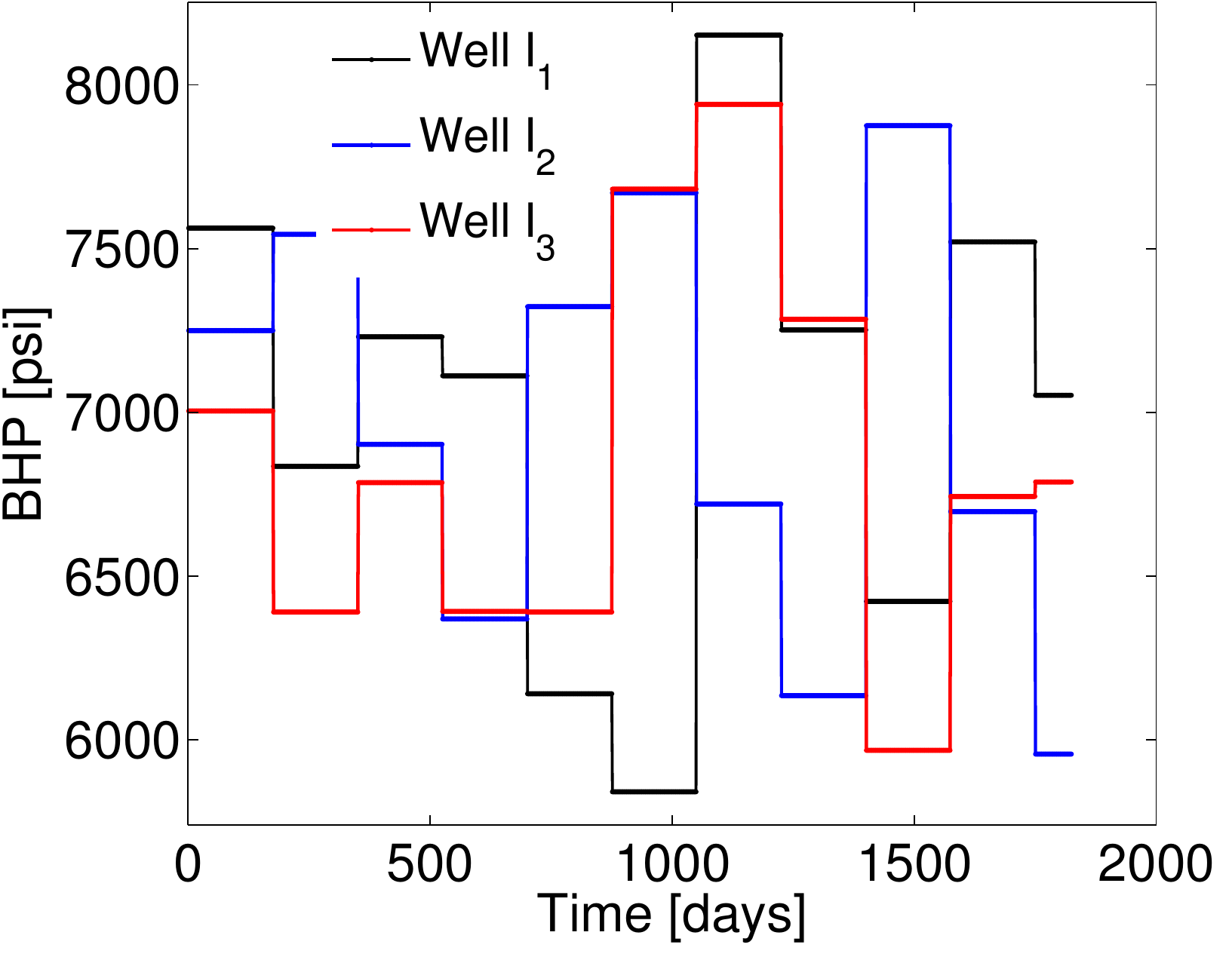}}	       
		\caption{Injector BHP profiles}
	\end{subfigure}
 \caption {Test Case 1 -- BHP profiles.}
 \label{fig:BHPprofileTestCases}
\end{figure}

We consider $\numTestSamples=200$ sets of BHP controls to construct the EMML training
$\paramTrain$ and test sets $\paramTest$. As described in Section~\ref{subsection:TrainingTestDataset:EMML}, each of these sets
is characterized by time-varying BHPs $\uNPlusOne$, $n=1,\ldots,\numTimesteps$,
obtained by adding a unique (time-varying) random perturbation to the primary training BHP
profiles. The time-varying BHPs for a particular case (Case~1) are
shown in \figRefOne{fig:BHPprofileTestCases}. The frequency of change in the
primary training BHP schedule is every 200~days (\figRefOne{fig:case1BHPTraining}),
while the frequency of change in the BHP schedule for Case~1 is every 175~days
(\figRefOne{fig:BHPprofileTestCases}).
%
% Test BHP profile pictures

We consider a total of $\numFeatures=84$ features, which include those listed in Table \ref{table:features}. After neglecting highly
correlated features, the total number of retained features is reduced to around
40 (note that each QoI may retain a different subset of features).
If we consider a memory of $\howMuchMemory{1}$, we obtain
$\numFeatures=84\times2=168$ features; this is reduced to about 64 after neglecting highly correlated features.

Based on extensive numerical experiments, we observed that the highest EMML accuracy was obtained on 
the test set $\paramTest$ using Approach~1 (i.e., 
error-modeling Method 4 in Section \ref{sec:method4} for the well-block pressure, and
Method 3 in Section \ref{sec:method3} for the well-block saturation),
a memory of $\howMuchMemory{1}$, 
$\numHFS{30}$ EMML training points,
classification for \reviewer{determining regression-model locality for
production-well QoIs},
and random forests (RF) for regression.
 While performing
classification, we set $\epsilon_{A}=0.05$ and $\epsilon_{C}=0.6$. \reviewer{These values are somewhat heuristic, but they appropriately identify the basic behaviors (solution stages) we wish to capture through classification.}
%(see Figure\ref{fig:labels}). 
\reviewer{We compute the hyper-parameters for random forests by minimizing the out-of-bag error.}
We first report the numerical results corresponding to these best-case
parameters. Then, in
Section~\ref{section:AdditionalAlgorithmictreatments}, we quantify EMML
performance for other choices of algorithmic parameters\reviewer{; for
example, we assess the effect of employing
clustering to determine regression-model locality, using only global
regression models, and applying LASSO for regression.} From Section \ref{subsection:application} we recall that there are two possible
applications of the EMML QoI-error prediction: (1) as a correction to the
surrogate-model QoI, or (2) as an error indicator to be used within the
ROMES framework. Here we consider the first
application, and in Section \ref{subsection:EMML_ROMES_Exepriments} the second.

%
%------------------------------------
\subsection{EMML for QoI correction: Test Case 1}\label{sec:EMMLtestCase1}
\begin{figure}[H] %-- Pressure --
	\begin{subfigure}  {1\linewidth} \centering 
	\label{fig:TestCase1_PrSw_Tr30_M1_Classification_RF_WellBlockPressure_P1}
	        \includegraphics[width=0.73\textwidth, angle= 0,origin=c]	{{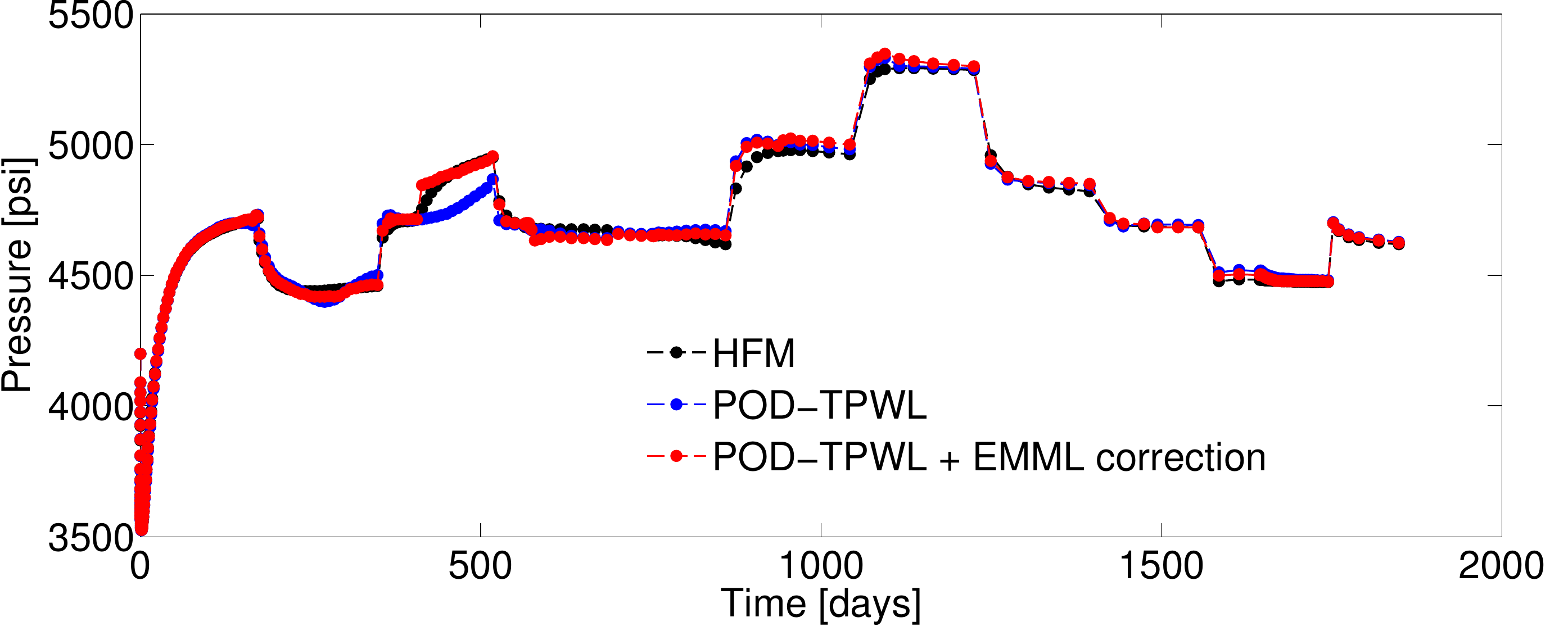}}   
  	        \caption {Well-block pressure at production well $\wellName{P}{1}$}	
	\end{subfigure}
\\
% Saturation
	\begin{subfigure}  {1\linewidth} \centering 
	\label{fig:TestCase1_PrSw_Tr30_M1_Classification_RF_WellBlockSaturation_P1}
	        \includegraphics[width=0.73\textwidth, angle= 0,origin=c]	{{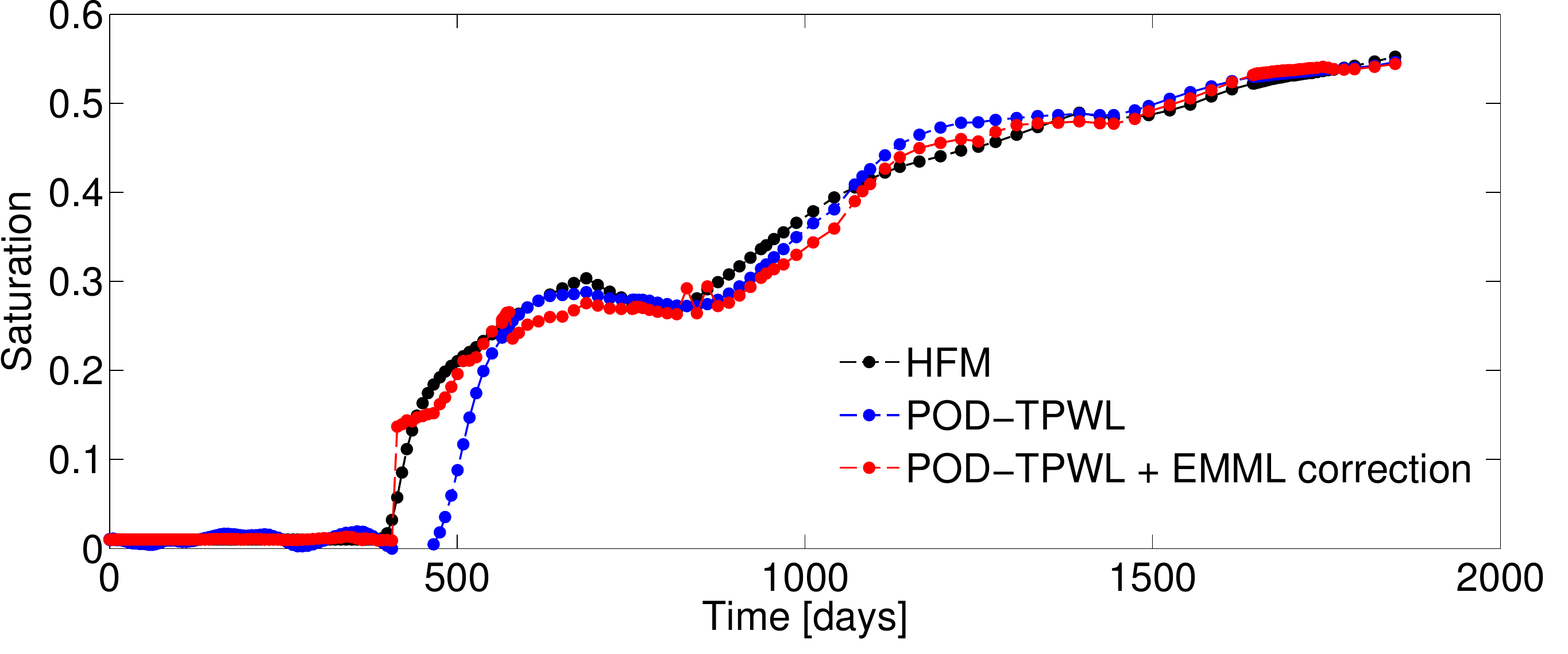}}   
  	        \caption {Well-block saturation at production well $\wellName{P}{1}$}
	\end{subfigure}
\\
	\begin{subfigure}  {1\linewidth} \centering 
	\label{fig:TestCase1_PrSw_Tr30_M1_Classification_RF_WellBlockPressure_I3}
	        \includegraphics[width=0.73\textwidth, angle= 0,origin=c]	{{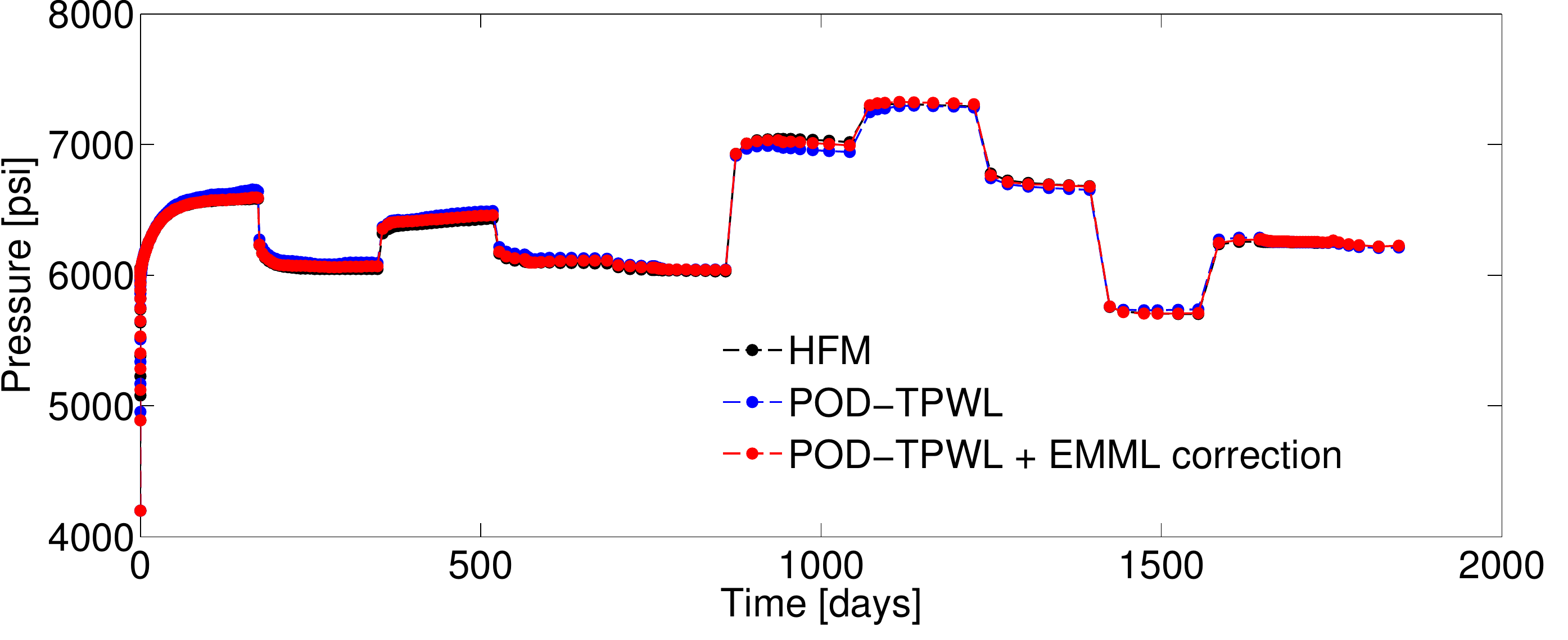}}   
  	        \caption {Well-block pressure at injection well $\wellName{I}{3}$}	
	\end{subfigure}
	\caption {EMML for sampled-state correction -- Test Case 1. Well-block pressure and saturation 
	predicted by various models. Best-performing EMML parameters:
	$\numHFS{30},~\howMuchMemory{1}$, classification \reviewer{for determining regression-model locality}, random-forest regression. }
\label{fig:TestCase1_PrSw_Tr30_M1_Classification_RF}
\end{figure}

\begin{figure}[H]
		 	%-- Oil --
	\begin{subfigure}  {1\linewidth} \centering 
	\label{fig:TestCase1Rates_Tr30_M1_Classification_RF_Oil}
	        \includegraphics[width=0.73\textwidth, angle= 0,origin=c]	{{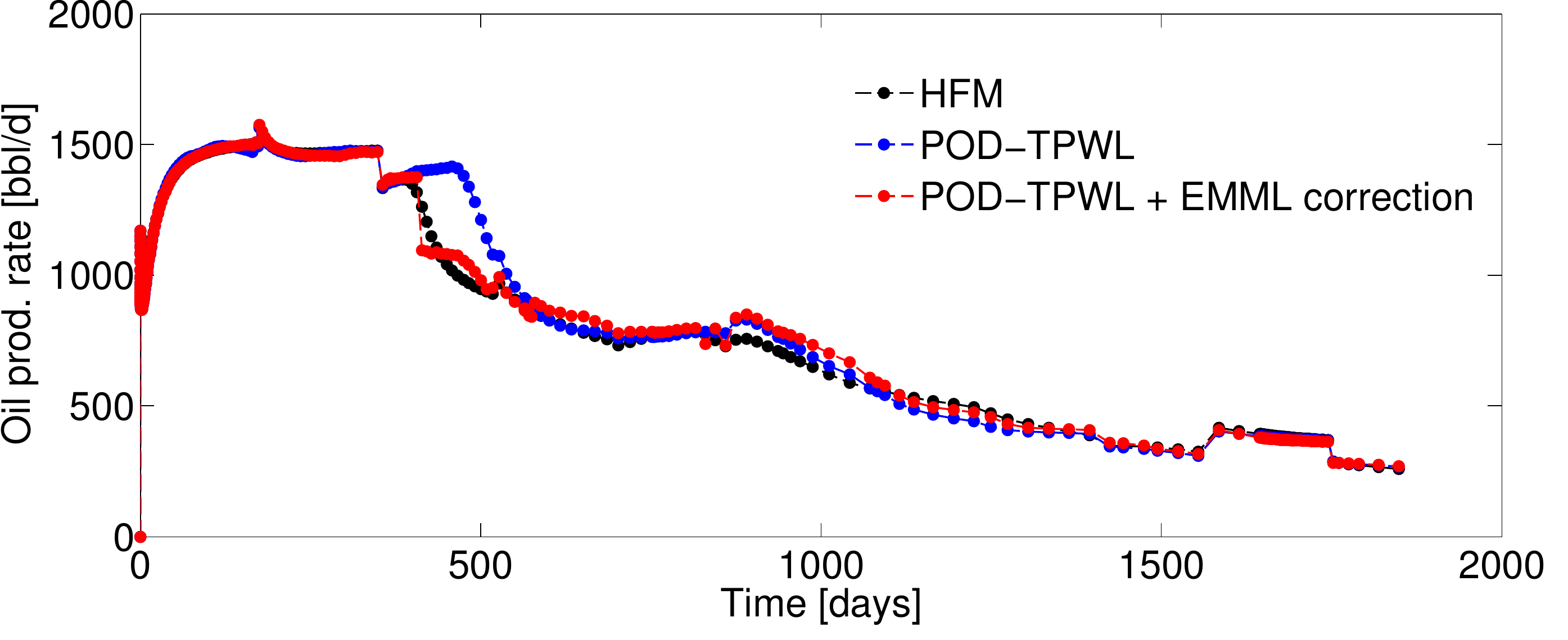}}   
  	        \caption {Oil production rate at well $\wellName{P}{1}$}
	\end{subfigure}
\\  %-- Water --
	\begin{subfigure}  {1\linewidth} \centering
	\label{fig:TestCase1Rates_Tr30_M1_Classification_RF_Water}
	        \includegraphics[width=0.73\textwidth, angle= 0,origin=c]	{{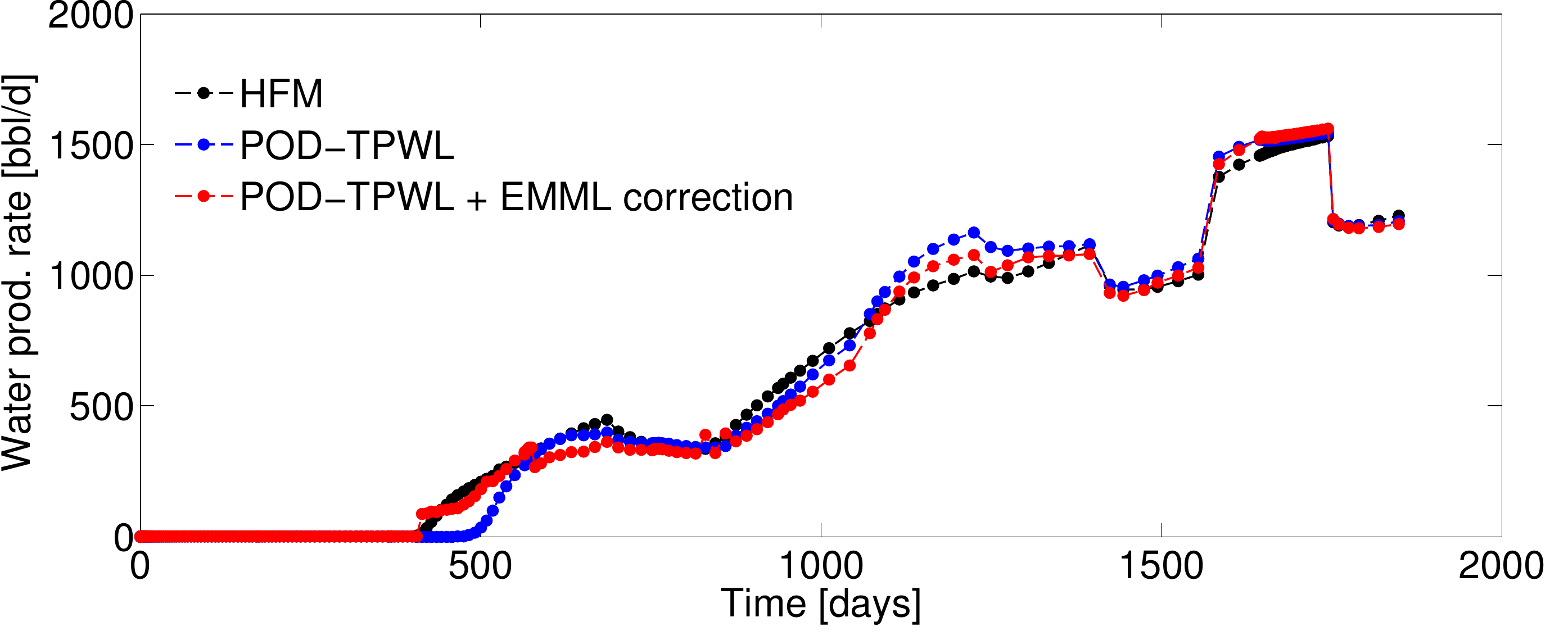}}
  	        \caption {Water production rate at well $\wellName{P}{1}$}
	\end{subfigure}
\\		 	%-- Inj --
	\begin{subfigure} {1\linewidth} \centering
	\label{fig:TestCase1Rates_Tr30_M1_Classification_RF_Inj}
	        \includegraphics[width=0.73\textwidth, angle=0,origin=c]	{{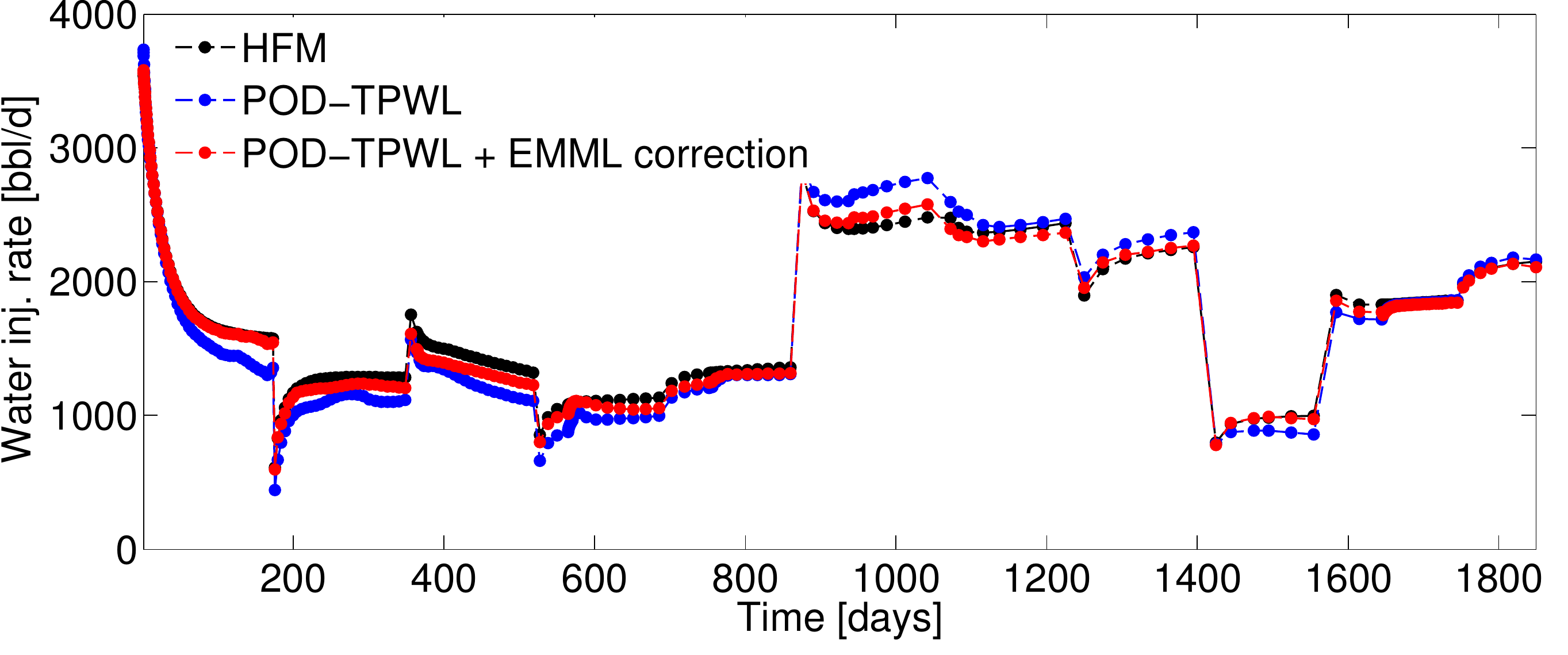}}
  	        \caption {Water injection rate at well $\wellName{I}{3}$}
	\end{subfigure}
\caption {EMML for QoI correction -- Test Case 1. Production and injection rates
predicted by various models. Best-performing EMML parameters:
$\numHFS{30},~\howMuchMemory{1}$, classification \reviewer{for determining regression-model locality},
random-forest regression.}
\label{fig:TestCase1Rates_Tr30_M1_Classification_RF}
\end{figure}

%------------------------------------
We first present results for Test Case $1$, represented by $\paramsTestArg{1}
\in \paramTest$. As will be described in 
Section~\ref{subsection:additionalTestCase},
this case corresponds to the median
time-integrated POD--TPWL error in the test set
$\paramTest$. 
\figRefOne{fig:TestCase1_PrSw_Tr30_M1_Classification_RF} reports
results for the 
pressure for wells $\wellName{P}{1}$ and $\wellName{I}{3}$ and the
saturation for well $\wellName{P}{1}$; note that these quantities associate
with the sampled state that is corrected as an intermediate step in Methods 3 and 4.
\figRefOne{fig:TestCase1Rates_Tr30_M1_Classification_RF} presents several QoI; these correspond to the oil and water
production rates in well $\wellName{P}{1}$, and the water injection rate for
well $\wellName{I}{3}$. We focus on wells
$\wellName{P}{1}$ and $\wellName{I}{3}$, as they are the wells with the highest
cumulative liquid production and injection. 
%-----------------------------------------------------
%------------ Test case 1  ------------
%-----------------------------------------------------

%
% -------Flow rate figures start--------------------- 
% --- Test case 423  ---

For Test Case $1$, \reviewer{the} POD--TPWL prediction (blue curve) for the
saturation for well $\wellName{P}{1}$ has an error that is most
noticeable at around $500$~days in
\figRefOne{fig:TestCase1_PrSw_Tr30_M1_Classification_RF}b. Similarly, the
POD--TPWL error in the production rates is evident at around the same time in
\figRefOne{fig:TestCase1Rates_Tr30_M1_Classification_RF}a,b. In
\figRefOne{fig:TestCase1_PrSw_Tr30_M1_Classification_RF}b, we observe that the
EMML-corrected well-block saturation (red curve) demonstrates improved accuracy around the time of water 
breakthrough ($\sim 410$~days). Similarly, in
\figRefOne{fig:TestCase1Rates_Tr30_M1_Classification_RF}, we see that the
EMML-corrected flow rates display better accuracy than the POD--TPWL results. The
improvement is most apparent in the breakthrough prediction in
\figRefOne{fig:TestCase1Rates_Tr30_M1_Classification_RF}b, and in oil
production rate (\figRefOne{fig:TestCase1Rates_Tr30_M1_Classification_RF}a) at
a time of about 500~days.  

%

%
%
%
%------------------------------------
\subsection{EMML for QoI correction: additional test cases}
\label{subsection:additionalTestCase}
%------------------------------------
We now present results for two additional test cases with control
vectors $\paramsGen \in \paramTest$, which correspond to 
different POD--TPWL prediction errors. We then assess EMML performance for an
ensemble containing the entire EMML test set $\paramTest$ ($|\paramTest|=170$ cases).

\begin{figure}[H] % Producer: Test case 497
	\begin{subfigure}  {0.50\linewidth} \centering
		\includegraphics[width= 0.75\textwidth]	{{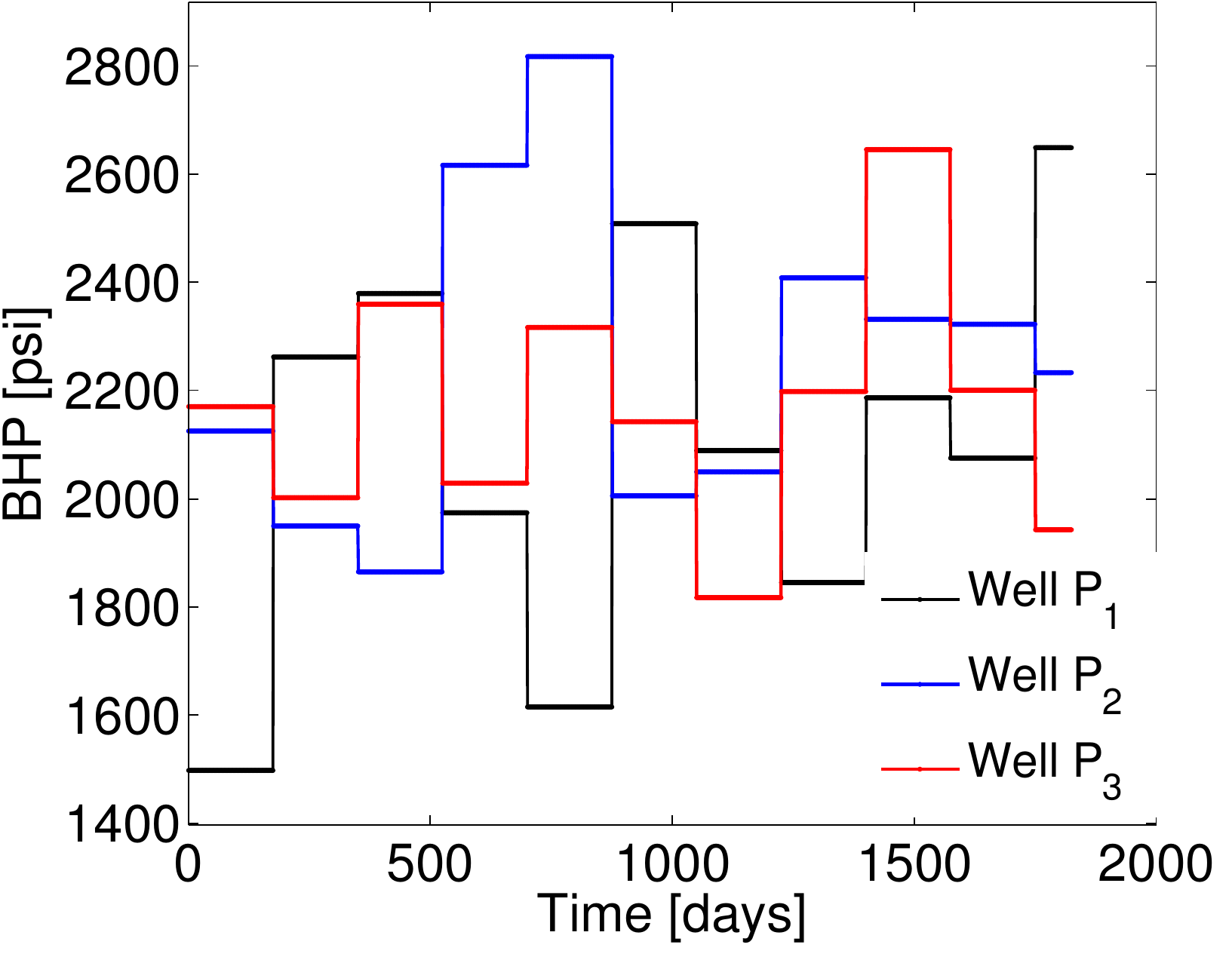}}	       
		\caption{Producer BHP profiles}
	\end{subfigure}
\quad  % Injector: Test case 497
	\begin{subfigure}  {0.50\linewidth} \centering
		\includegraphics[width= 0.75\textwidth]	{{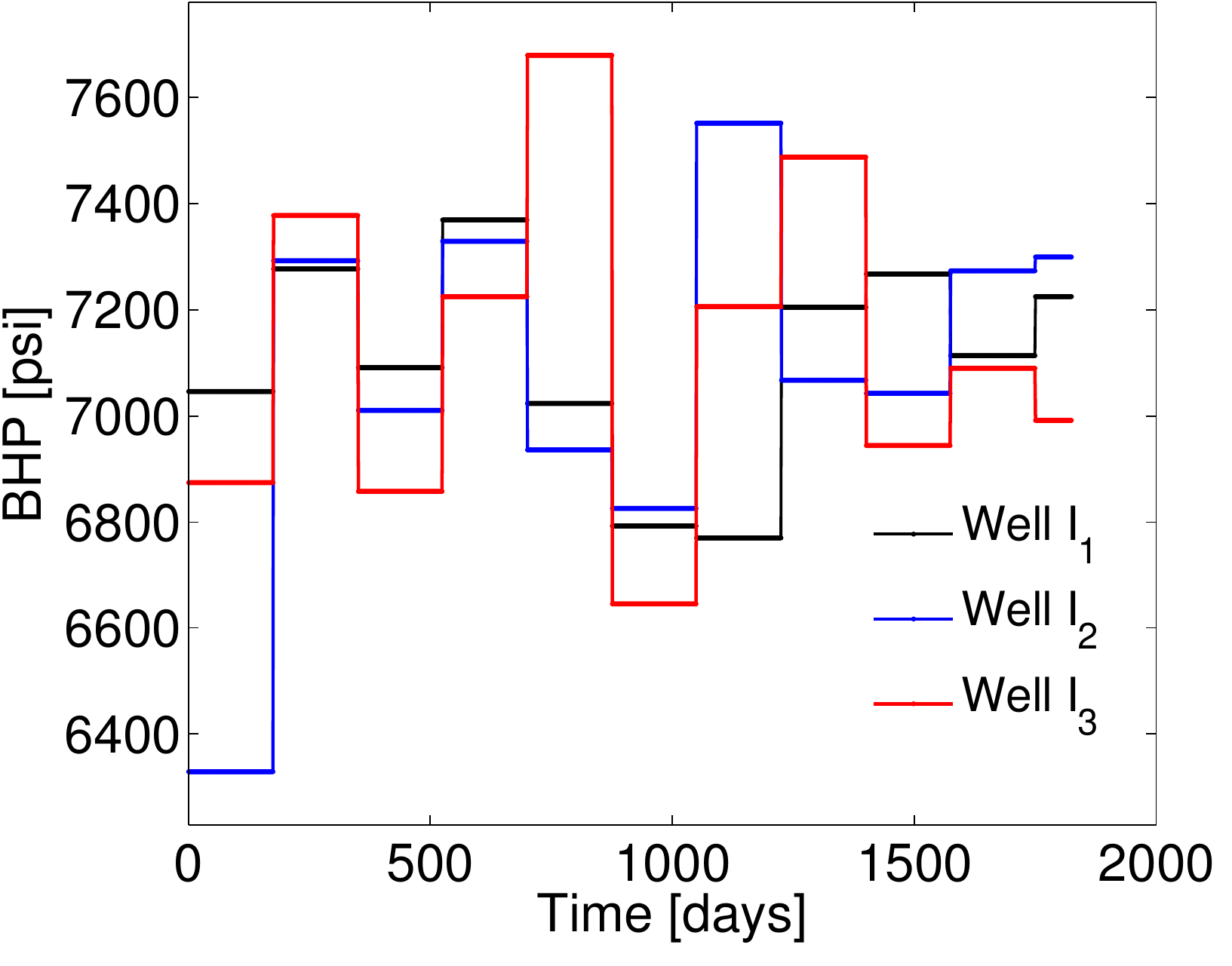}}	       
		\caption{Injector BHP profiles}
	\end{subfigure}
 \caption {Test Case 2 -- BHP profiles.}
 \label{fig:BHPprofileTestCase2}
\end{figure}

\begin{figure}[H]	%-- Water --
		 	%-- Oil --
	\begin{subfigure}  {1\linewidth} \centering 
	        \includegraphics[width=0.73\textwidth, angle= 0,origin=c]	{{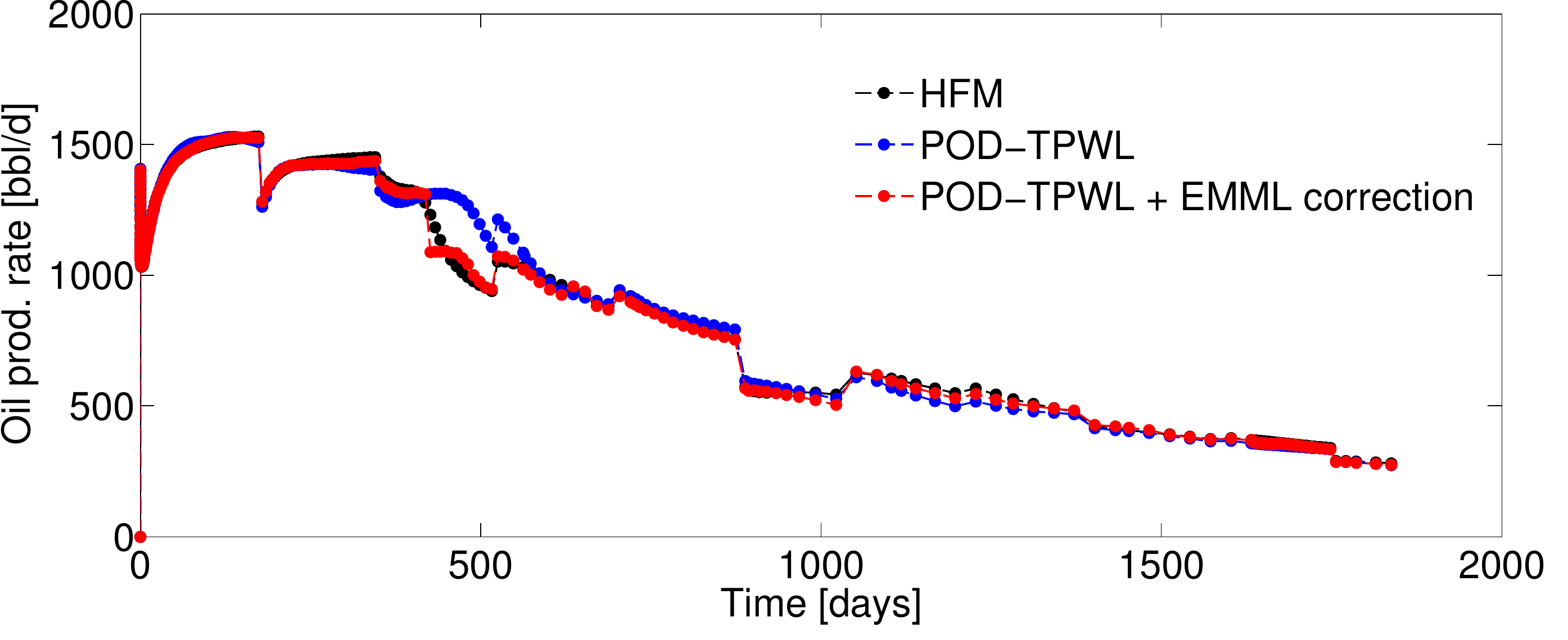}}   
  	        \caption {Oil production rate at well $\wellName{P}{1}$}
	\end{subfigure}
\\  %-- Water --
	\begin{subfigure}  {1\linewidth} \centering
	        \includegraphics[width=0.73\textwidth, angle= 0,origin=c]	{{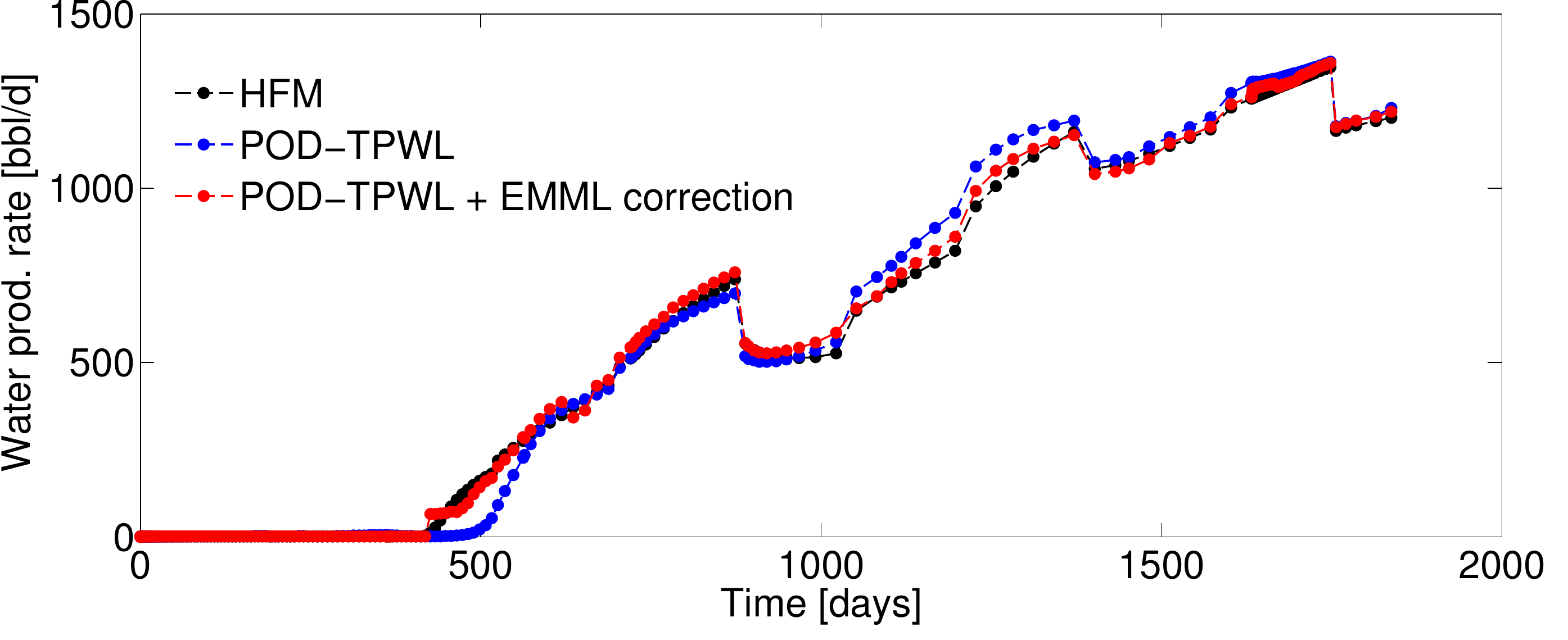}}
  	        \caption {Water production rate at well $\wellName{P}{1}$}
	\end{subfigure}
\\		 	%-- Inj --
	\begin{subfigure} {1\linewidth} \centering
	        \includegraphics[width=0.73\textwidth, angle=0,origin=c]	{{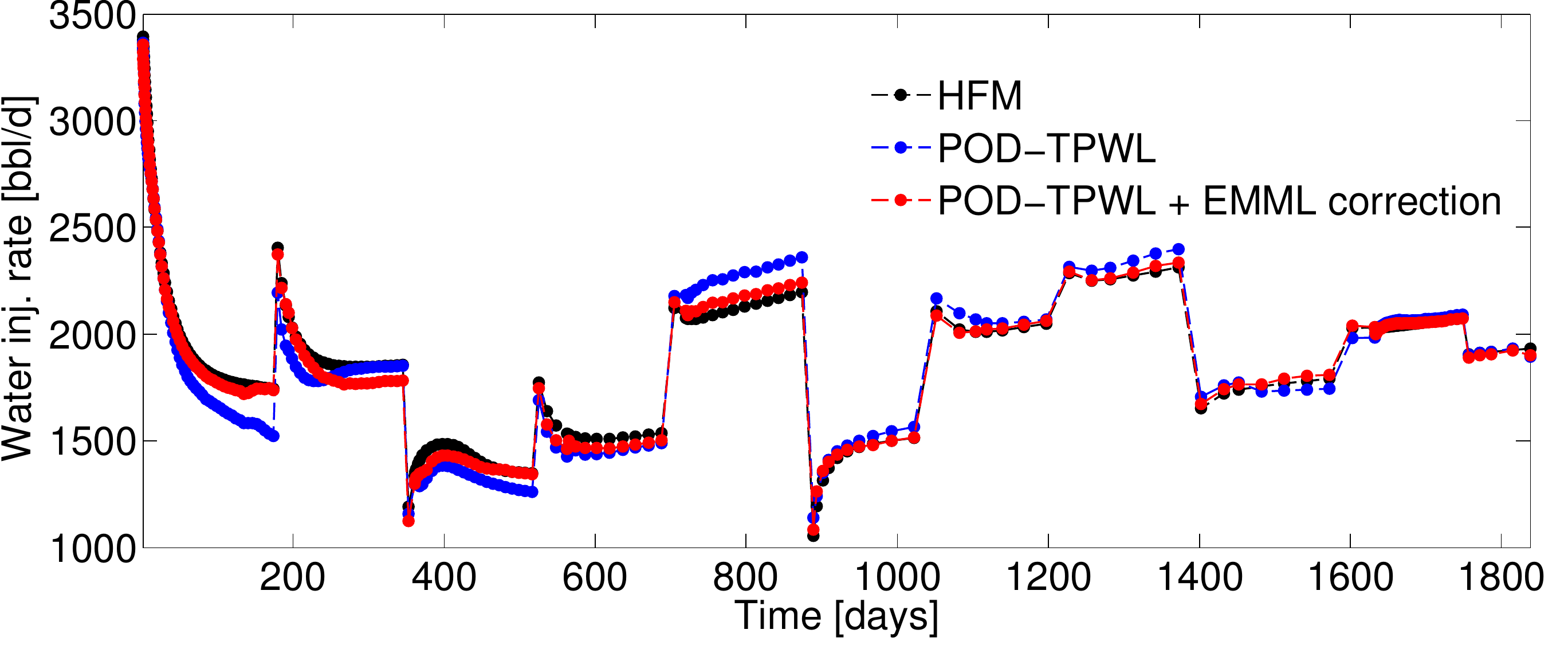}}
  	        \caption {Water injection rate at well $\wellName{I}{3}$}
	\end{subfigure}
\caption {EMML for QoI correction -- Test Case 2. Production and injection rates
predicted by various models. Best-performing EMML parameters:
$\numHFS{30},~\howMuchMemory{1}$, classification \reviewer{for determining regression-model locality},
random-forest regression.}
\label{fig:TestCase2Rates_Tr30_M1_Classification_RF}
\end{figure}

BHP schedules for Test Cases 2 and 3---represented by control vectors
$\paramsTestArg{2}$ and $\paramsTestArg{3}$, respectively---are shown in
\figRefTwo{fig:BHPprofileTestCase2}{fig:BHPprofileTestCase3}. Test Case 2, for
which $\Delta \BHPnoBf^{P}(\paramsTestArg{2}) = 0.54$ and $\Delta
\BHPnoBf^{I}(\paramsTestArg{2}) = 0.10$,
corresponds to a smaller perturbation in the BHPs relative to the primary training run BHPs compared to that in Test
Case 1 ($\Delta \BHPnoBf^{P}(\paramsTestArg{1}) = 0.68,~\Delta \BHPnoBf^{I}(\paramsTestArg{1}) = 0.23$). It also
corresponds to lower POD--TPWL error compared to Test Case 1. The results for
production and injection rates for Test Case 2 are shown in
\figRefOne{fig:TestCase2Rates_Tr30_M1_Classification_RF}. The POD--TPWL error
is again most noticeable at around 500~days for the oil and water production
rates at well $P_1$. The
correction is clearly evident in
\figRefOne{fig:TestCase2Rates_Tr30_M1_Classification_RF}a at around 500~days
and in \figRefOne{fig:TestCase2Rates_Tr30_M1_Classification_RF}b at around 500
and 1250~days. Slight improvement in water injection rate (\figRefOne{fig:TestCase2Rates_Tr30_M1_Classification_RF}c) is also apparent just before 200~days. 

\begin{figure}[H]	%-- Water --
		 	%-- Oil --
	\begin{subfigure}  {1\linewidth} \centering 
	        \includegraphics[width=0.73\textwidth, angle= 0,origin=c]	{{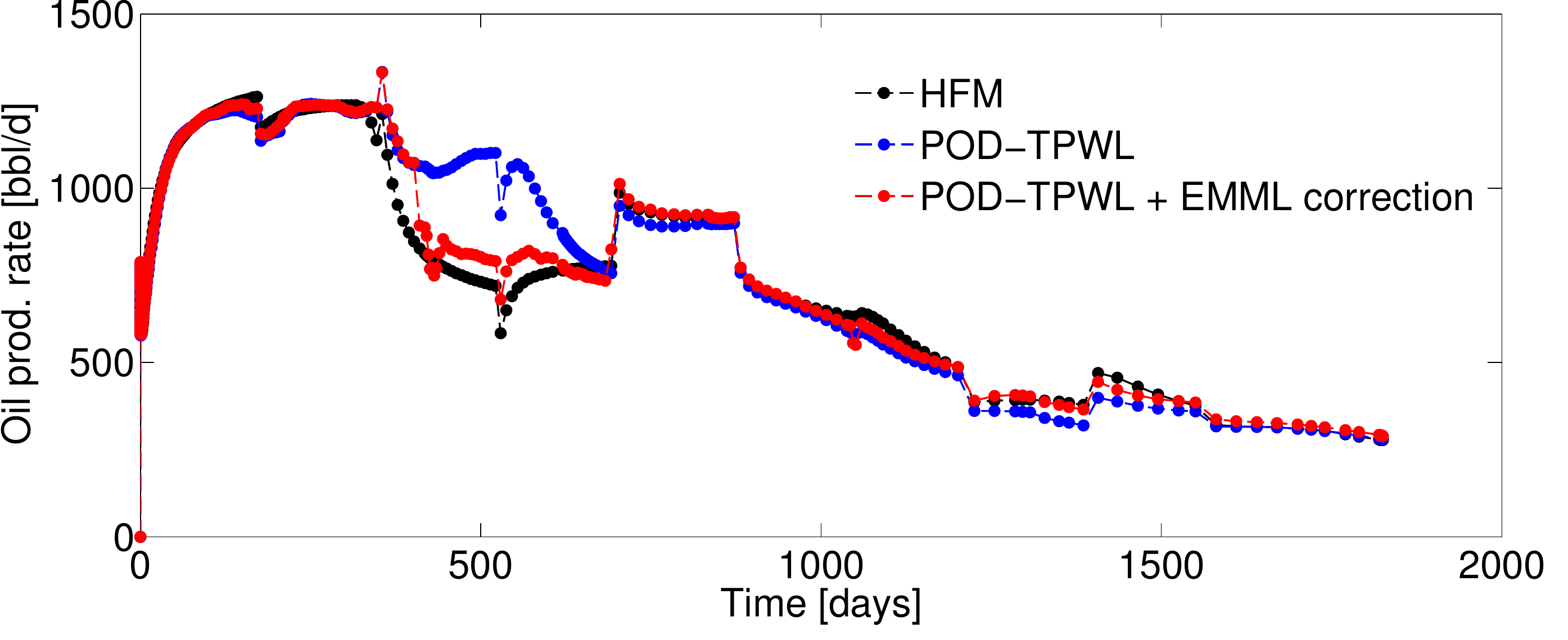}}   
  	        \caption {Oil production rate at well $\wellName{P}{1}$}
	\end{subfigure}
\\  %-- Water --
	\begin{subfigure}  {1\linewidth} \centering
	        \includegraphics[width=0.73\textwidth, angle= 0,origin=c]	{{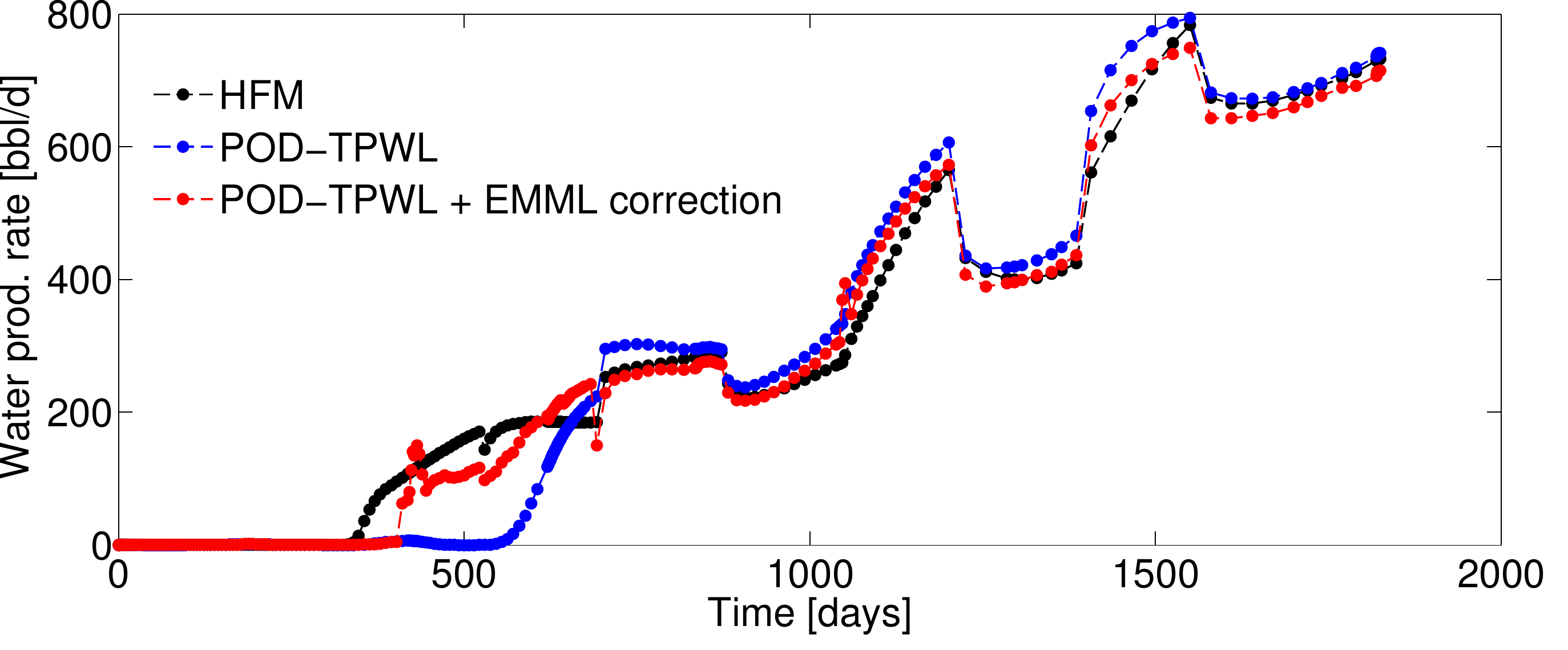}}
  	        \caption {Water production rate at well $\wellName{P}{1}$}
	\end{subfigure}
\\		 	%-- Inj --
	\begin{subfigure} {1\linewidth} \centering
	        \includegraphics[width=0.73\textwidth, angle=0,origin=c]	{{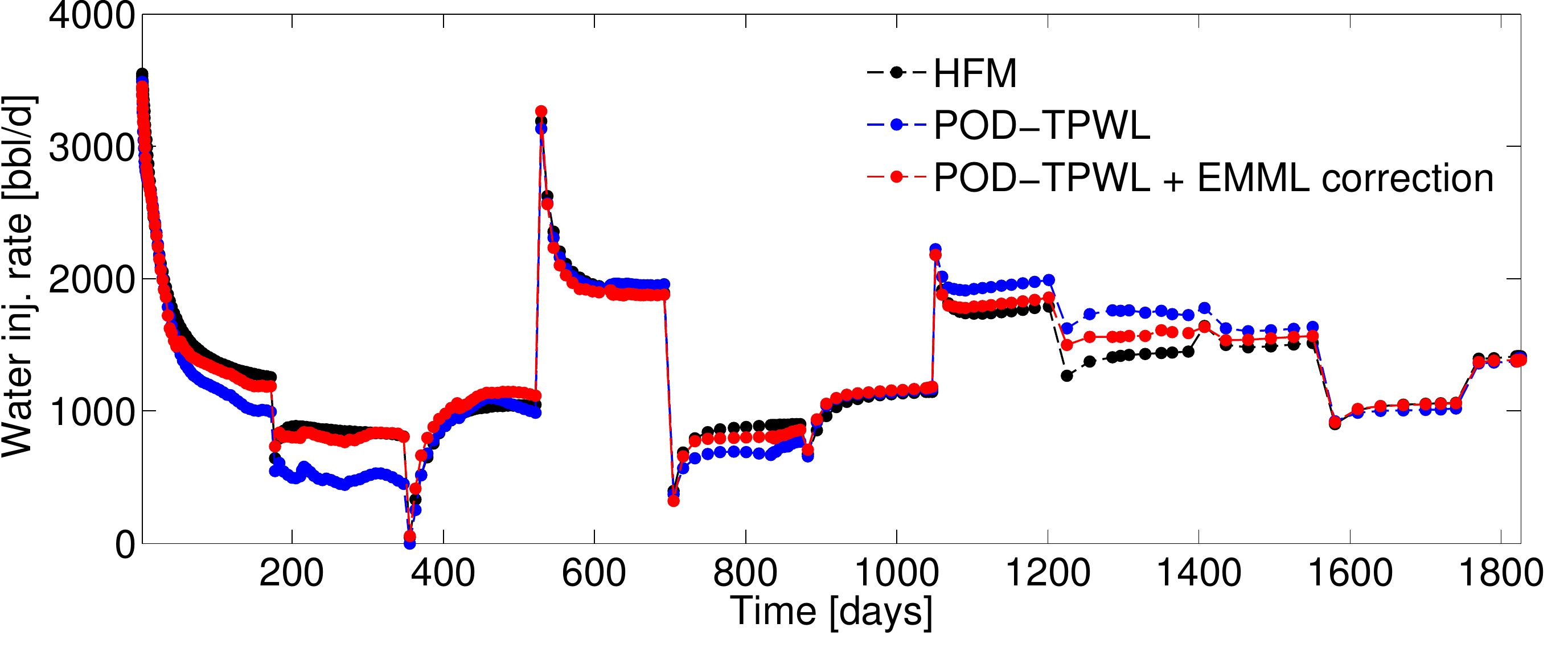}}
  	        \caption {Water injection rate at well $\wellName{I}{3}$}
	\end{subfigure}
\caption {EMML for QoI correction -- Test Case 3. Production and injection rates
predicted by various models. Best-performing EMML parameters:
$\numHFS{30},~\howMuchMemory{1}$, classification \reviewer{for determining regression-model locality},
random-forest  regression.}
\label{fig:TestCase3Rates_Tr30_M1_Classification_RF}
\end{figure}

Test Case 3, with $\Delta \BHPnoBf^{P}(\paramsTestArg{3}) = 0.46$ and $\Delta
\BHPnoBf^{I}(\paramsTestArg{3}) = 0.32$, corresponds to a higher perturbation
in the injector BHPs compared to that in Test Case 1, and it leads to a larger
POD--TPWL error. The POD--TPWL error is again most evident at around 500~days
for both the oil and water production rates as shown in
	\figRefOne{fig:TestCase3Rates_Tr30_M1_Classification_RF}a,b. These results
	are again significantly improved by the proposed EMML-based correction. We note finally that the corrected
solutions for the production and injection rates display fluctuations at some times. This is because, when
constructing the corrections, we treat each time step as independent, consistent with the i.i.d.
assumption. 

To quantify EMML performance over the entire test set $\paramTest$ of
$|\paramTest|=170$ cases, we define the following relative time-integrated error
measures for the POD--TPWL and corrected solutions:
\begin{subequations}
\label{eqn:errorDefinitionTPWL:EMML}
\begin{align}
	E_{\RL}(j,\wellBlockSetDummy) &=  \frac{1}{\lvert \wellBlockSetDummy \rvert} \mathlarger{\mathlarger{\sum}}\limits_{
	q\in\{(q_j)_d|d\in\wellBlockSetDummy\}}   \frac {\mathlarger{\sum}_{n=1}^{\numTimesteps}  \left\lvert   
	\errorOperatorQoITimeN
	\right\rvert \Delta t^n } { \mathlarger{\sum}_{n=1}^{\numTimesteps} q^n    \Delta t^n } \times 100\%,  \\
	E_{\rm{corr}}(j,\wellBlockSetDummy) &=  \frac{1}{\lvert \wellBlockSetDummy \rvert} \mathlarger{\mathlarger{\sum}}\limits_{
	q\in\{(q_j)_d|d\in\wellBlockSetDummy\}}
		   \frac
{\mathlarger{\sum}_{n=1}^{\numTimesteps}   \left\lvert
\errorOperatorQoIApproxTimeN
-  \errorOperatorQoITimeN  \right\rvert  \Delta t^n } {
\mathlarger{\sum}_{n=1}^{\numTimesteps} q^n   \Delta t^n }
\times 100\%.
\end{align}
\end{subequations}
Here, $E_{\RL}$ denotes the relative average time-integrated error in the POD--TPWL
solution and $E_{\rm{corr}}$ designates the average time-integrated error in the
EMML-corrected solution. 
Note that $\wellBlockSetDummy =
\wellBlockSetProducer$ (with $j=o$ or $j=w$) for production wells, while 
$\wellBlockSetDummy = \wellBlockSetInjector$ (with $j=w$) for injection wells.

\begin{figure}[]	 % -- Oil --
 	\begin{subfigure}  {0.5\linewidth} \centering %{0.5\linewidth} \centering 
		\includegraphics[width=1\textwidth]	{{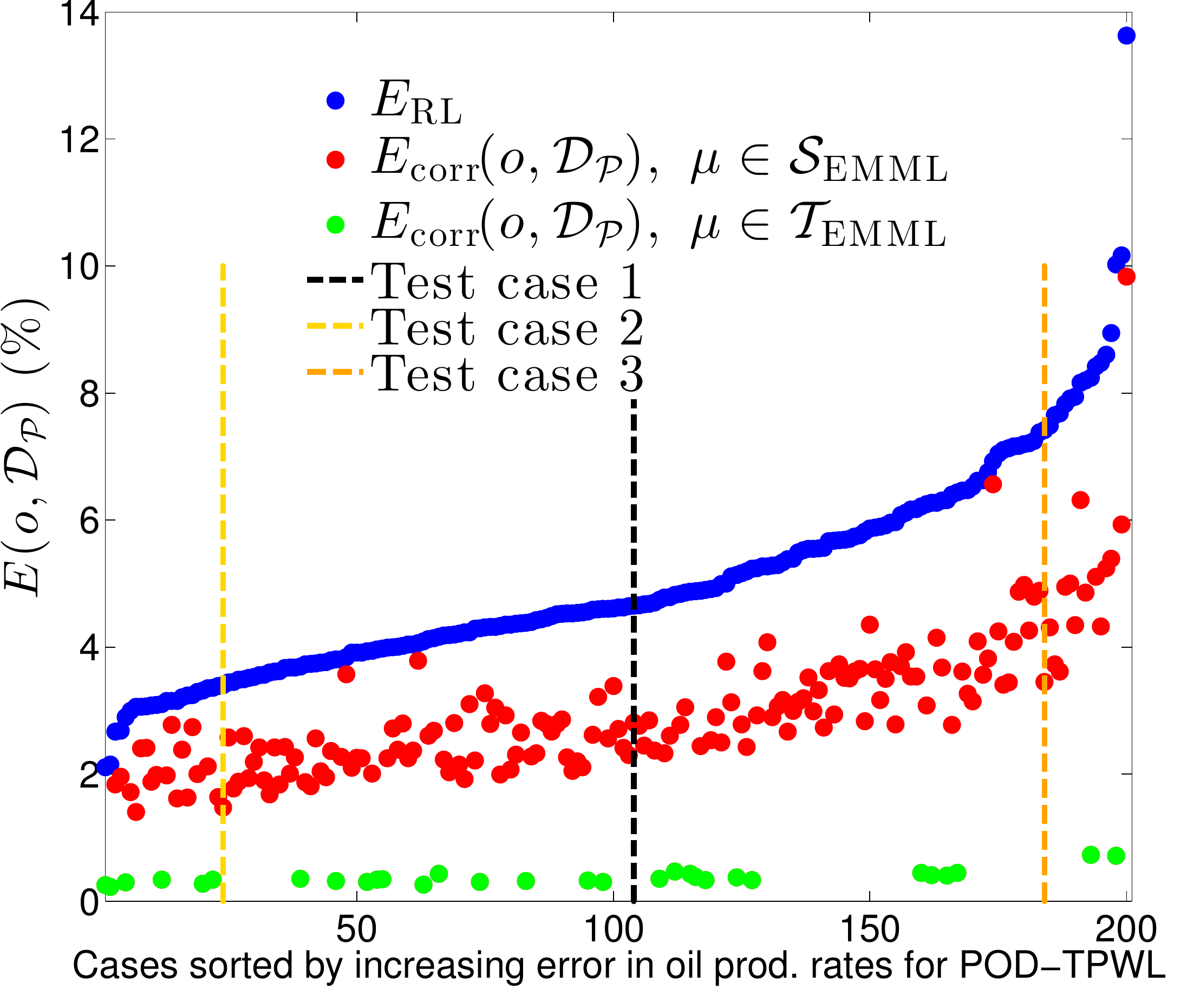}}	       
	        \caption{Oil production rate}
	\end{subfigure}
% \\
\quad
	\begin{subfigure}  {0.5\linewidth} \centering %{0.5\linewidth} \centering
		\includegraphics[width=1\textwidth]	{{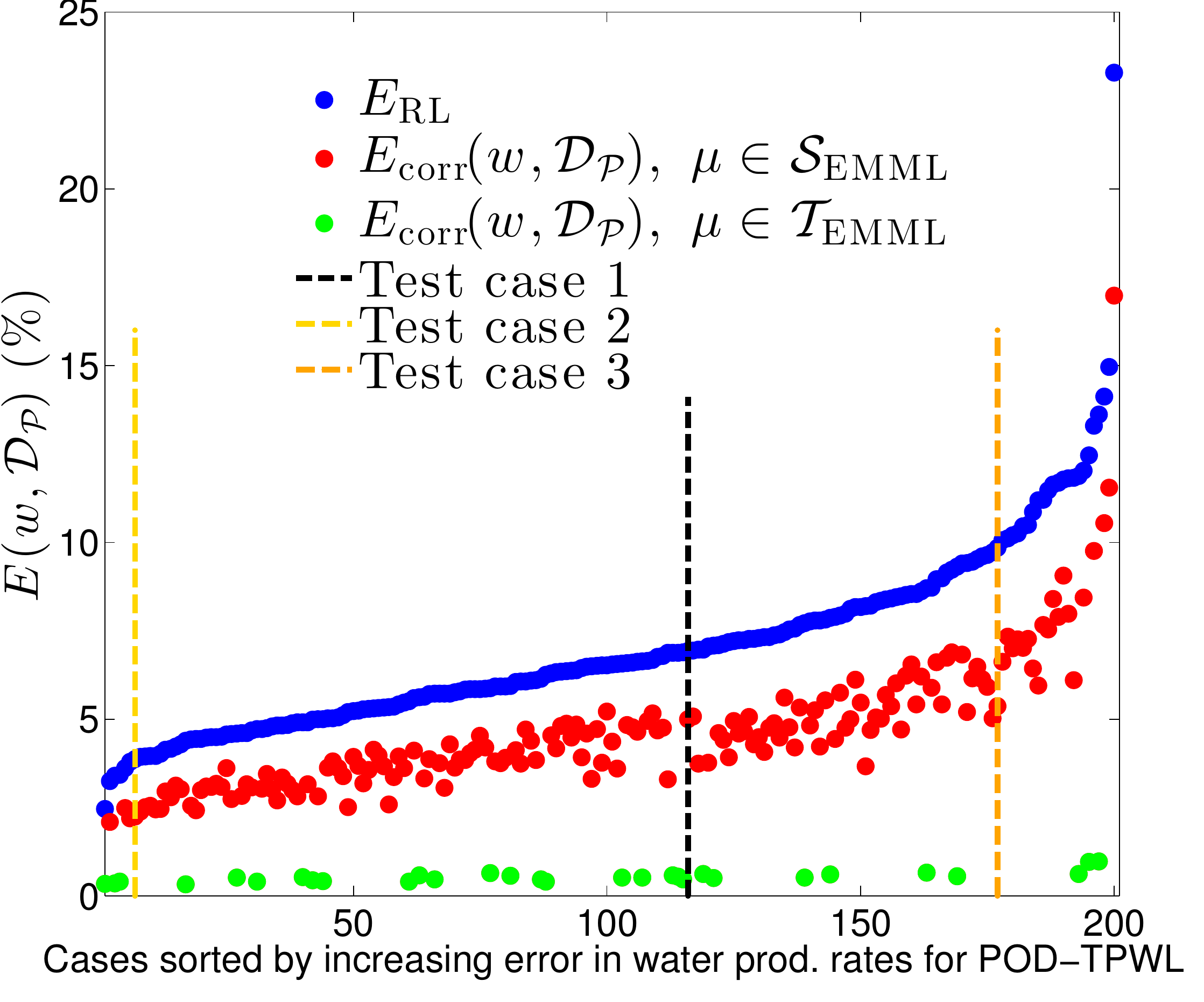}}	       
	        \caption{Water production rate}
	\end{subfigure}
\\
	\begin{subfigure}  {1\linewidth} \centering
		\includegraphics[width=0.5\textwidth]	{{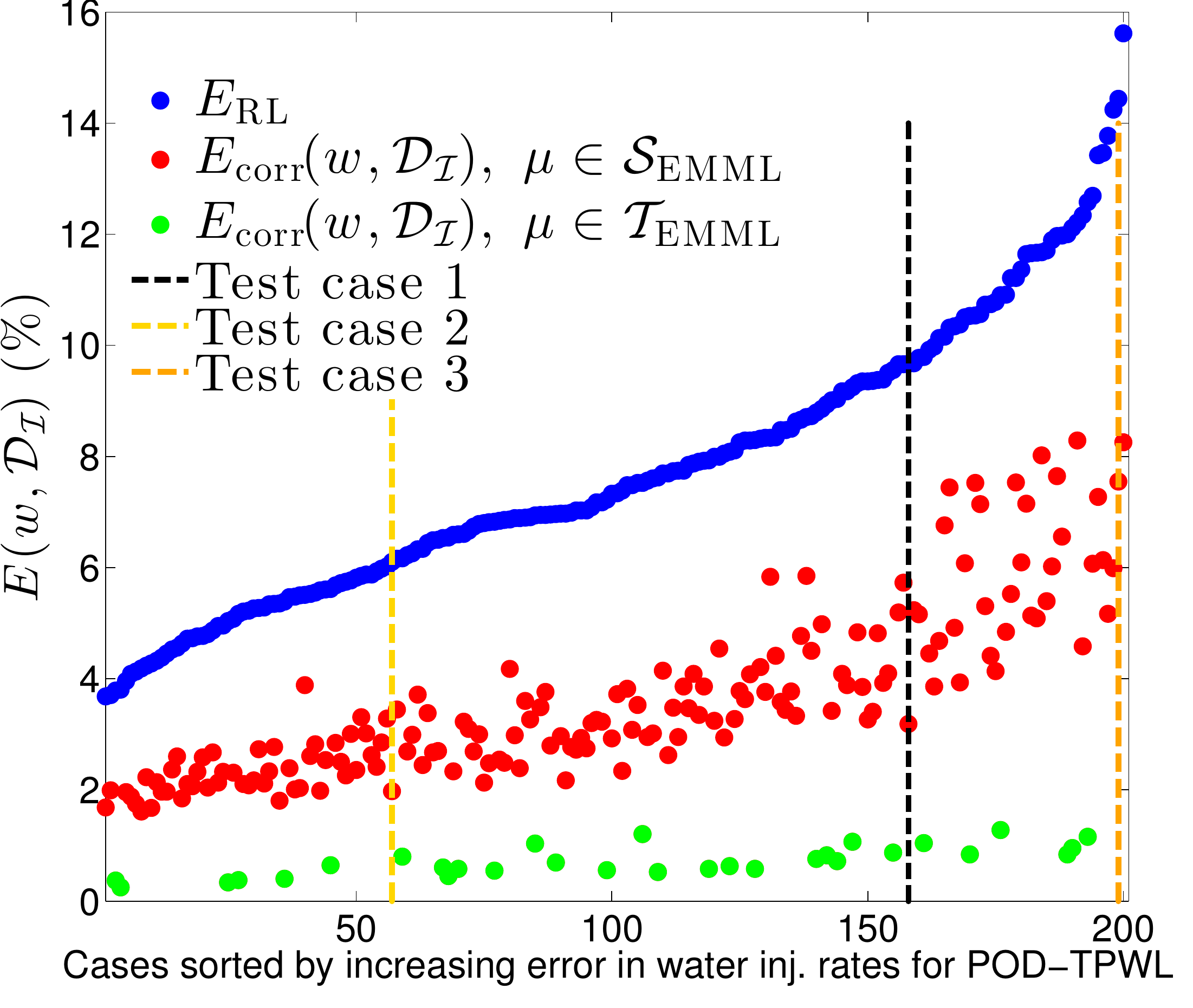}}	       
	        \caption{Water injection rate}
	\end{subfigure}
\caption {EMML for QoI correction: additional test cases.
	Relative time-integrated  error in production and injection rates as defined by
\eqnRefOne{eqn:errorDefinitionTPWL:EMML} for $\numHFS{30},~\howMuchMemory{1}$, classification
+ RF.
}
\label{fig:errorPlotMultipleCases_Tr30_M1_Classification_RF}
\end{figure}

\figRefOne{fig:errorPlotMultipleCases_Tr30_M1_Classification_RF} displays the
time-integrated errors for the entire set of 200 cases in the EMML training $\paramTrain$ and
test $\paramTest$ sets. 
The cases are sorted by increasing
POD--TPWL error. 
For each
case in the ensemble, the figure reports the time-integrated error in the POD--TPWL prediction
(blue) and the EMML-corrected POD--TWPL predictions, which are further
distinguished by whether they correspond to cases in the EMML training set
$\paramTrain$ (green) or test set $\paramTest$ (red). Consistent with the
results presented in Section \ref{sec:EMMLtestCase1}, the time-integrated
errors for \textit{all} test cases are reduced after
application of the EMML-driven correction. Test Cases 1, 2 and 3, discussed
above, are identified in
\figRefOne{fig:errorPlotMultipleCases_Tr30_M1_Classification_RF}. These three test
cases correspond to the $10$th, $50$th and $90$th
percentiles in $E_{\RL}(o,\wellBlockSetProducer)$. Note that the
time-integrated error in the EMML training data (green points) is small but
nonzero, which indicates that the random forest model
$\errorModelPrimaryVarApprox{i}$ does not perfectly fit the data. This is
intentional, as it prevents overfitting, which can potentially lead to large
errors in EMML test-case predictions.

%
%
% --- Flow rate for test case 2, 3 start ---
%-----------------------------------------------------
%------------ Test case 2,3  ------------
%-----------------------------------------------------
% Test case 497
% Test BHP profile pictures

% Test case 616
% Test BHP profile pictures
\begin{figure}[H] % Producer: Test case 616
	\begin{subfigure}  {0.50\linewidth} \centering
		\includegraphics[width= 0.75\textwidth]	{{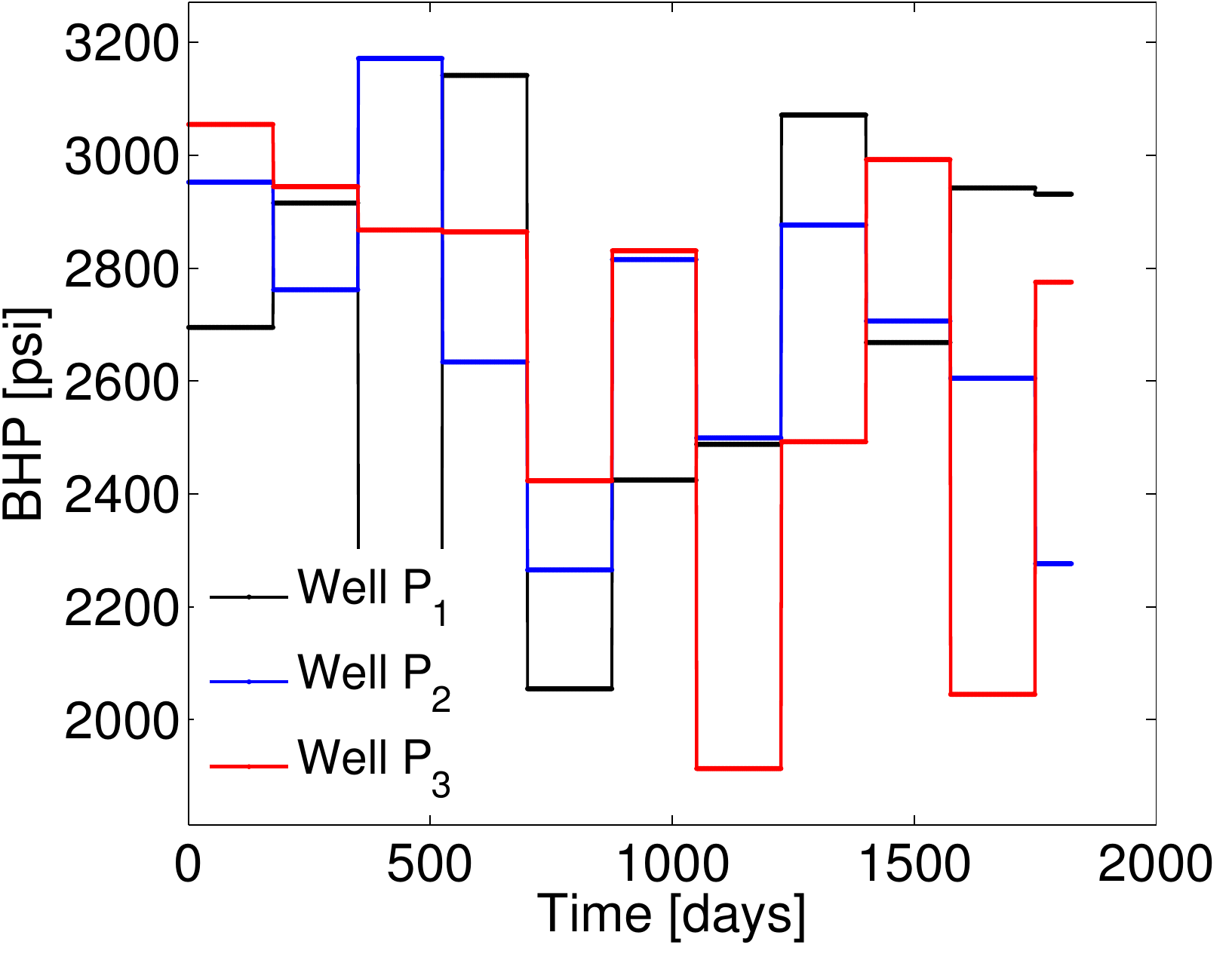}}	       
		\caption{Producer BHP profiles}
	\end{subfigure}
\quad  % Injector: Test case 616
	\begin{subfigure}  {0.50\linewidth} \centering
		\includegraphics[width= 0.75\textwidth]	{{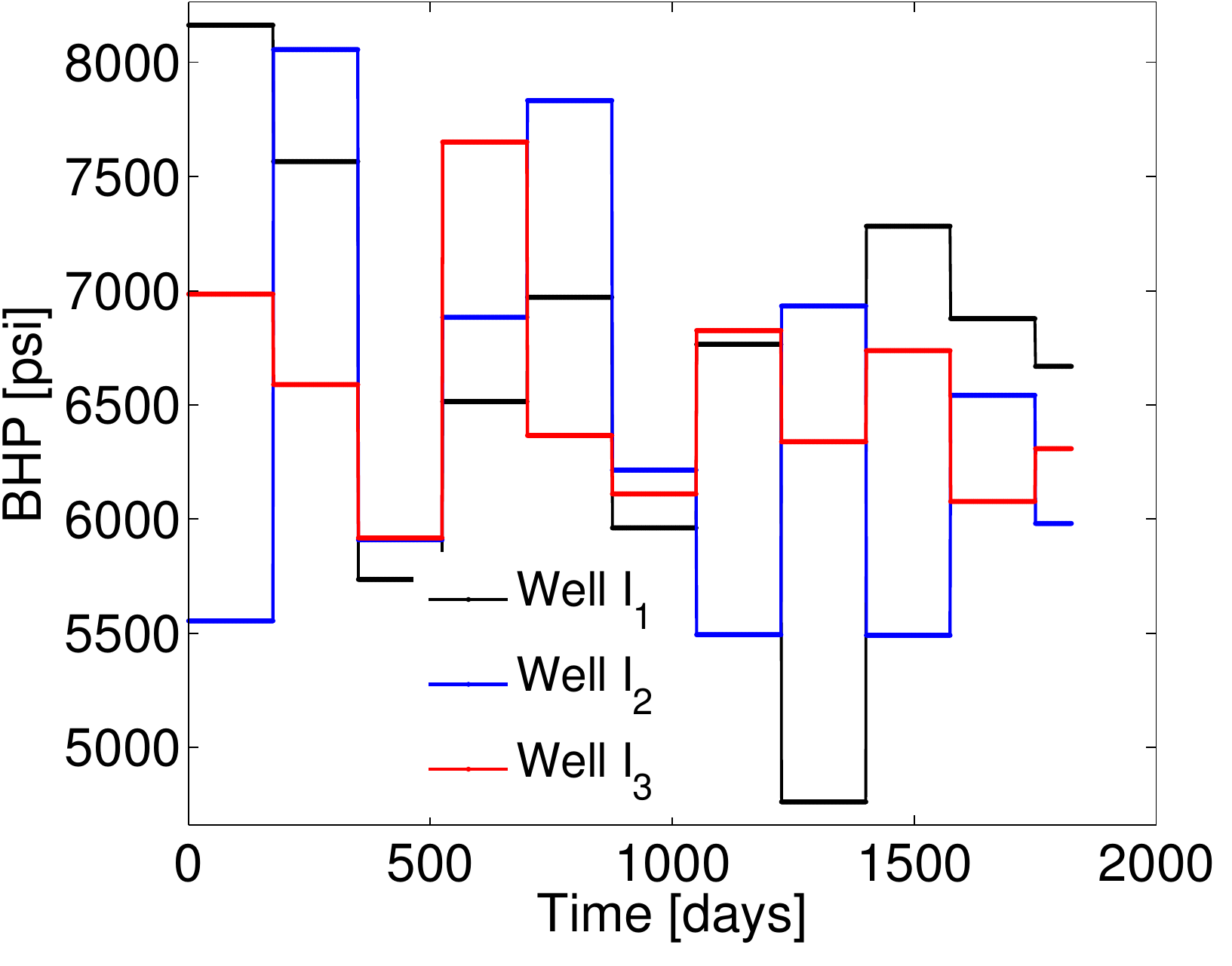}}	       
		\caption{Injector BHP profiles}
	\end{subfigure}
 \caption {Test Case 3 -- BHP profiles.}
 \label{fig:BHPprofileTestCase3}
\end{figure}  

%-----------------------------------------------------
%------------ Time integrated errors - all 200 cases ----- 
%-----------------------------------------------------
%
%
% --- Cumulative error in oil, water, and injection rates--- 
Table~\ref{table:errorComparison:EMML:TPWL} presents the median errors for the test set $\paramTest$ results displayed in
\figRefOne{fig:errorPlotMultipleCases_Tr30_M1_Classification_RF}. We observe
that by applying the EMML correction, we reduce the three errors,
$E_{\RL}(o,\wellBlockSetProducer)$, $E_{\RL}(w,\wellBlockSetProducer)$ and
$E_{\RL}(w,\wellBlockSetInjector)$, by about $38\%$ on average.
Although we achieve substantial improvements at certain time instances using
EMML (as is evident in
\figRefOne{fig:TestCase1Rates_Tr30_M1_Classification_RF}a,b),
small errors in time persist. The EMML procedure reduces these errors but it does not completely
eliminate them.
\begin{table}[]
\caption{EMML for QoI correction: additional test cases. Median value of the time-integrated errors in 
POD--TPWL and corrected
solutions over cases in EMML test set $\paramsGen\in\paramTest$}
\centering
\label{table:errorComparison:EMML:TPWL}
\begin{tabular}{c| c c c }
\toprule
\hline
 Method & {$E_{i}(o,\wellBlockSetProducer)$} & {$E_{i}(w,\wellBlockSetProducer)$} &
{$E_{i}(w,\wellBlockSetInjector)$} \\  \hline 
POD--TPWL ($i=\RL$)  & 4.5\% & 6.3\% & 6.6\% \\
EMML correction ($i=$corr) &2.8\%    &4.4\%    &3.7\%  \\
\bottomrule
\end{tabular}
\end{table} 
 We note that one source of error in the EMML predictions is
 misclassification, which in turn leads to using the local regression
 model from the incorrect category. The
average misclassification error---defined as the ratio of the number of EMML
test samples misclassified to the total number of EMML test samples---over cases in the EMML test set
$\paramTest$ is 3\% in this set of experiments.
The misclassification error is primarily from misclassifying samples
whose actual category is $B^+$ as $B^-$, and vice-versa.

\reviewer{Finally, from the variable importance plots generated by the random forest regression model (see~~\cite{breiman2001random} for more details on variable importance plots), we observe some of the key features in this application to be $\sampleRomOutputState_{2d}\zN - \sampleRomOutputState_{2d}\zI$, $\sampleRomOutputState_{2d-1}\zN - \sampleRomOutputState_{2d-1}\zI$, $\sampleRomOutputState_{2d}\zNMinusOne - \sampleRomOutputState_{2d}\zIMinusOne$, $\sampleRomOutputState_{2d}\zN$, $\sampleRomOutputState_{2d}\zI$, $\dSdTn$, $\BHPnoBf^{\timestep}_{d}$, ${\rm{PVI}}^{\timestep}$, and $ \dfrac{\innerProduct{\zN}{\zI }} 	{\norm{\zN}{2}~\norm{\zI}{2}}$.}

%------------------------------------
\subsection{EMML for QoI correction: alternative EMML parameters}
\label{section:AdditionalAlgorithmictreatments}
%------------------------------------

For completeness, we now analyze EMML performance using different algorithmic
parameters. Table~\ref{table:errorComparison_allMethods} reports the median time-integrated errors
$E_{\rm{corr}}(j,\wellBlockSetDummy)$ over the 170 test cases corresponding to 
$\paramsGen\in\paramTest$. The EMML
parameters used here differ from the best-case parameters employed in Sections
\ref{sec:EMMLtestCase1}--\ref{subsection:additionalTestCase}. We vary the 
memory $\memory$, the number of
high-fidelity simulations $\numHFSsymbol$ used to construct the EMML training
data, and the \reviewer{method for determining regression-model locality} clustering or
classification). For a detailed
discussion of these results, we refer the reader to~\cite{trehanThesis}. 

We observe from Table~\ref{table:errorComparison_allMethods} that the EMML-based
corrections obtained using Approach~1 lead to more accurate results than those
obtained by Approach~2 (these two approaches are defined in Section~\ref{subsection:Features}). Additionally, a decrease in the number of high-fidelity simulations ($\numHFSsymbol$) used to build the EMML training
dataset leads to a noticeable decrease in accuracy. Reduced accuracy is also
observed when LASSO regression is used instead of random-forest regression.
In addition, using classification (a
supervised machine learning technique) \reviewer{to determine regression-model locality} prior to
\reviewer{constructing local RF models}  performs better than the use of
clustering (an automated unsupervised learning approach).
Even in the absence of \reviewer{employing local regression models} (which adds complexity to the EMML
method), EMML \reviewer{with a global error model} still provides improved accuracy relative to POD--TPWL, though
more accurate results are achieved when classification is used \reviewer{to
determine locality}. We note
finally that the impact of memory is very small for this test set.
\begin{table}[H]
\caption{
EMML for QoI correction: alternative EMML parameters. Median value of the time-integrated errors in 
POD--TPWL and corrected
solutions over cases in EMML test set $\paramsGen\in\paramTest$ with different
EMML parameters
}
\centering
\label{table:errorComparison_allMethods}
%\begin{tabular}{l c c c c c c c c c c c c}
\begin{tabular}{ c c c c c |  c c c }
%\toprule  \hline 
\toprule
\hline
\multicolumn{5}{c|}{\multirow{2}{*}{Method}}   & $E_{i}(o,\wellBlockSetProducer)$ &
$E_{i}(w,\wellBlockSetProducer)$  &
$E_{i}(w,\wellBlockSetInjector)$ \\
 &&&&& (\%)&
(\%) &
(\%)\\[1mm]  \hline 
\multicolumn{5}{c|}{POD--TPWL ($i=\RL$)}  & 4.49 & 6.34 & 6.62 \\
\hline
\multicolumn{5}{c|}{EMML correction ($i=$corr)}  & &  &  \\
Approach & $\memory$ & $\nTrain$ & \reviewer{Locality} & Regression & & & \\
\cline{1-5}
1 & 1 & 30 & classification & RF & 2.78 & 4.39 & 3.67\\
1 & 0 & 30 & classification & RF & 2.79 & 4.38 & 3.68\\
1 & 1 & 15 & classification & RF & 3.58 & 5.63 & 4.12\\
1 & 1 & 30 & classification & LS & 3.33 & 5.67 & 11.32\\
1 & 1 & 30 & clustering & RF & 3.23 & 5.22 & 3.68\\
1 & 1 & 30 & none & RF & 3.26 & 5.18 & 3.67\\
2 & 1 & 30 & classification & RF & 2.78 & 5.62 & 11.06\\
\bottomrule
\end{tabular}
\end{table} 

%------------------------------------
\subsection{EMML as an error indicator for ROMES}
\label{subsection:EMML_ROMES_Exepriments}
%------------------------------------
We now consider the second application of the EMML error model: as an error
indicator for the ROMES method \cite{drohmann2015romes}. 
In particular, we apply the framework proposed in  Section
\ref{sec:EMML_ROMES} with the following
scalar-valued function of the surrogate QoI errors:
\begin{equation}
\errorFun(\errorOperatorQoIArg{1},\ldots,\errorOperatorQoIArg{\numTimesteps})
=  \frac {\mathlarger{\sum}_{n=1}^{\numTimesteps}  \lvert
{\delta}_{d}^{\timestep}  \rvert \Delta t^{n} } {
	\mathlarger{\sum}_{n=1}^{\numTimesteps}  \left (
	\left(q^{\timestep}\right)_{\RL} +    {\delta}_{d}^{\timestep} \right)
	\Delta t^{n} }.
\end{equation}
We also denote the average value of $\errorFun$ associated with QoI errors
$\delta_q$
for a
given phase $j$ 
and a
specified set of wells $\bar{\mathcal D}$ as
\begin{equation} \label{eq:averageValueQoIError}
\bar \errorFun(\delta_q;j,\bar{\mathcal D}) = \frac{1}{\bar{\mathcal D}}
\sum_{q\in\{(q_j)_d|d\in\wellBlockSetDummy\}}\errorFun(\errorOperatorQoIArg{1},\ldots,\errorOperatorQoIArg{\numTimesteps})
\times 100\%.
\end{equation} 

%
% --- Estimated time-integerated error ---

%-----------------------------------------------------
%------------ Estimated Time integrated errors - all 200 cases 
%-----------------------------------------------------
%
%
\begin{figure}[H]	
 % -- Oil --
	\begin{subfigure}  {0.5\linewidth} \centering %{0.5\linewidth} \centering
		\includegraphics[width= 0.85\textwidth]	{{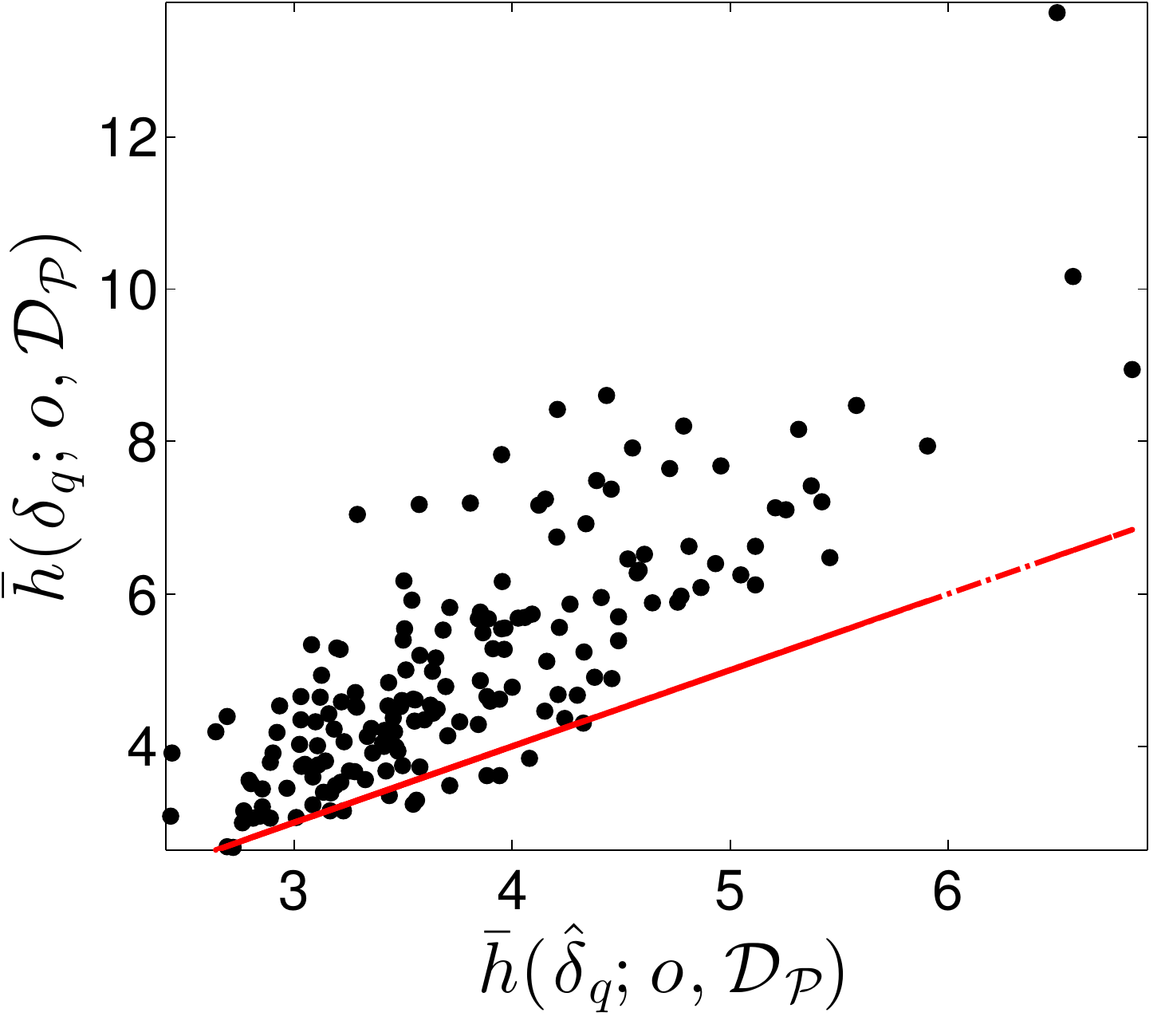}}	       
	        \caption{Oil production rate}
	\end{subfigure}
%\quad
\quad
	\begin{subfigure}  {0.5\linewidth} \centering %{0.5\linewidth} \centering
		\includegraphics[width= 0.85\textwidth]	{{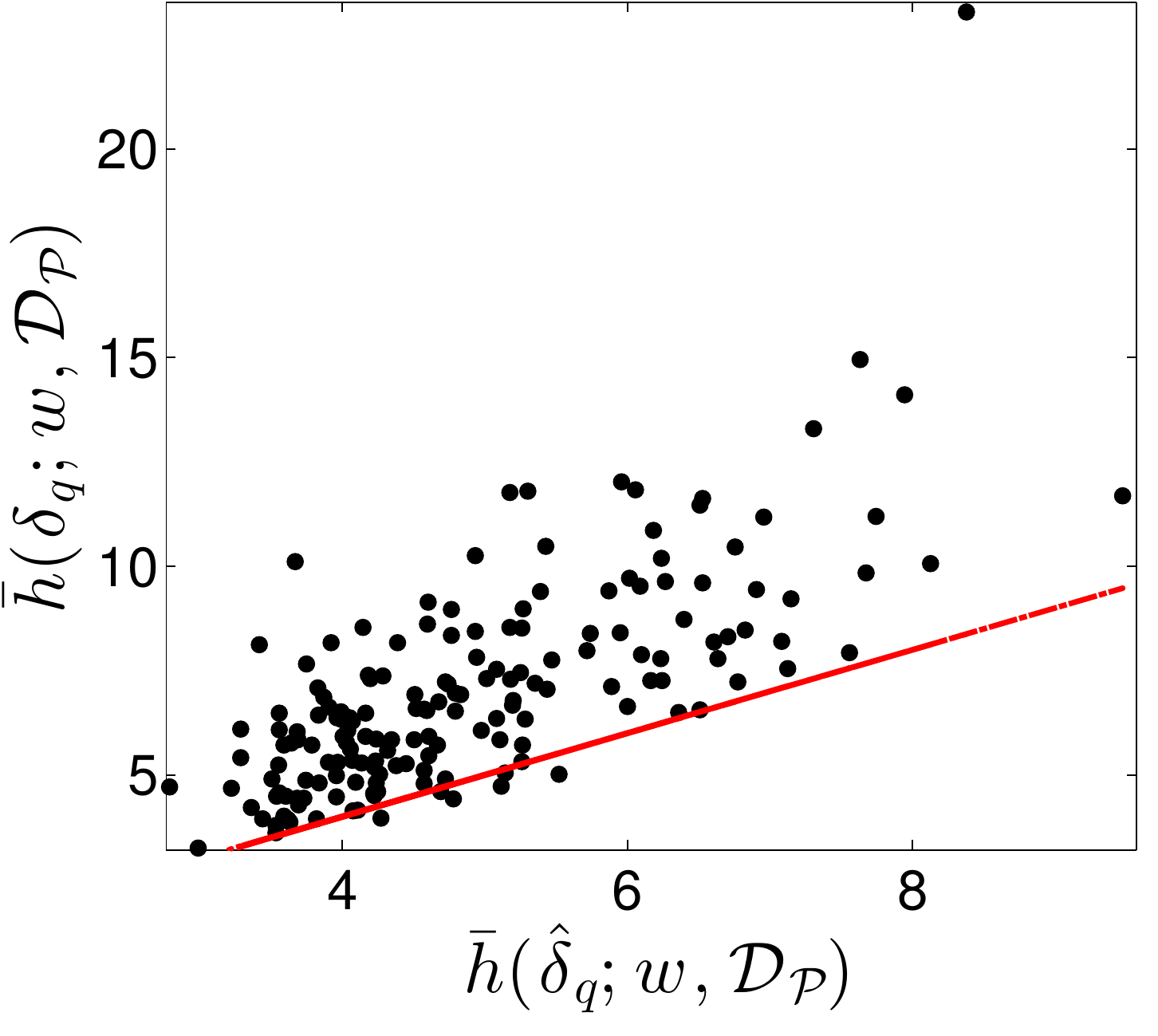}}	       
	        \caption{Water production rate}
	\end{subfigure}
\\
	\begin{subfigure}  {1\linewidth} \centering
		\includegraphics[width= 0.43\textwidth]	{{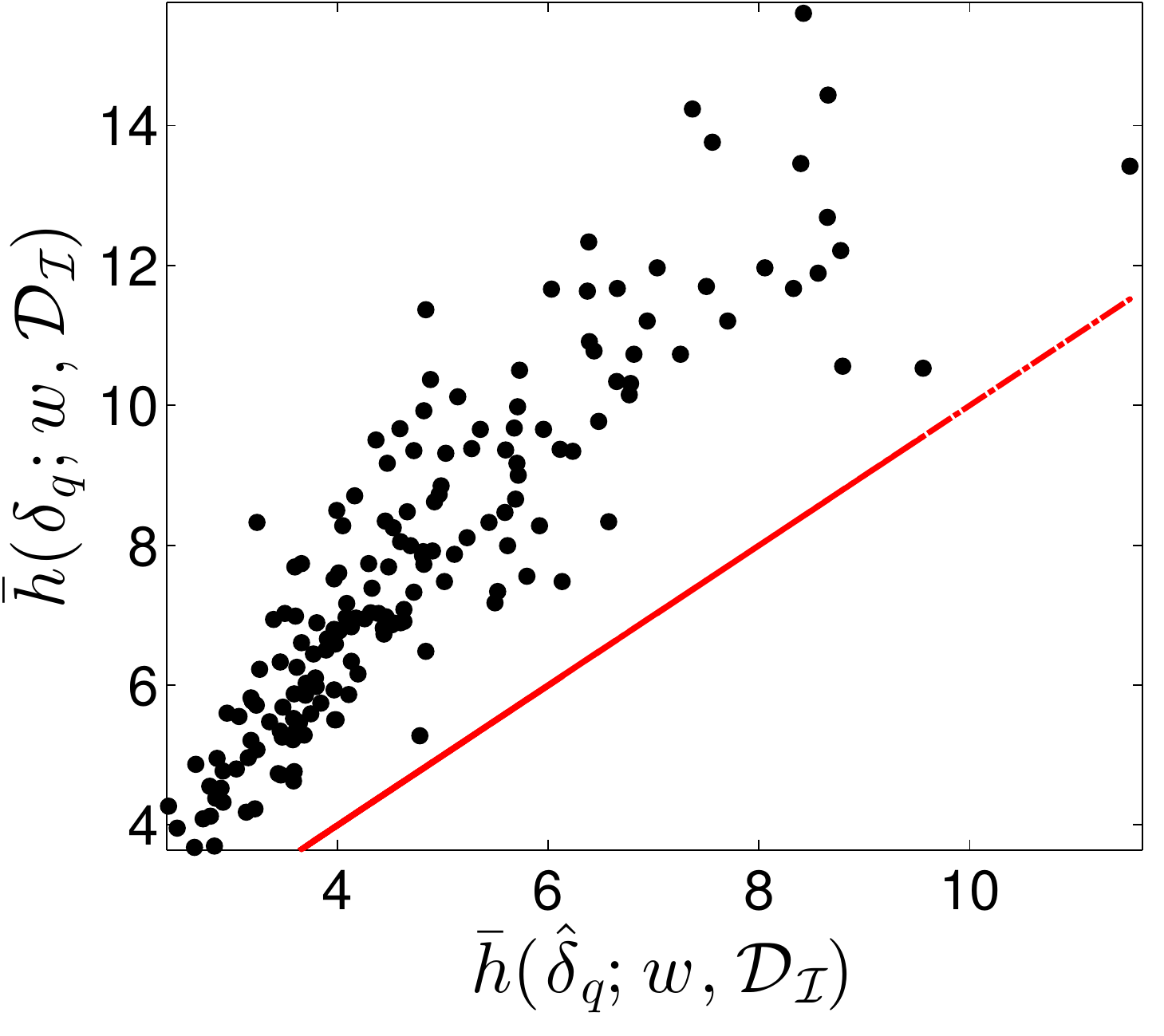}}	       
	        \caption{Water injection rate}
	\end{subfigure}
\caption {Relationship between
the true time-integrated error ($
 \bar \errorFun(\delta_q;o,\wellBlockSetProducer)
 $, $
 \bar \errorFun(\delta_q;w,\wellBlockSetProducer)
 $, and $
 \bar \errorFun(\delta_q;w,\wellBlockSetInjector)
 $)
and the EMML-\reviewer{approximated} time-integrated error 
($
 \bar \errorFun(\hat\delta_q;o,\wellBlockSetProducer)
 $,
$
 \bar \errorFun(\hat\delta_q;w,\wellBlockSetProducer)
 $, and
$
 \bar \errorFun(\hat\delta_q;w,\wellBlockSetInjector)
 $).
EMML parameters:	
$\nTrain=30$, $\howMuchMemory{1}$, classification \reviewer{for determining regression-model locality},
random-forest regression.
The red line corresponds to the prediction associated with the
EMML-\reviewer{approximated} error alone; this illustrates the bias in the
EMML-\reviewer{approximated}
error.}
\label{fig:predictedErrorPlotMultipleCases_Tr30_M1_Classification_RF}
\end{figure}  
\figRefOne{fig:predictedErrorPlotMultipleCases_Tr30_M1_Classification_RF}
displays cross-plots of the true time-integrated error 
($\bar \errorFun(\delta_q;o,\wellBlockSetProducer)$, $\bar \errorFun(\delta_q;w,\wellBlockSetProducer)$, and $\bar \errorFun(\delta_q;w,\wellBlockSetInjector)$)
versus the EMML-\reviewer{approximated} time-integrated error ($\bar \errorFun(\hat\delta_q;o,\wellBlockSetProducer)$, $\bar \errorFun(\hat\delta_q;w,\wellBlockSetProducer)$, and $\bar \errorFun(\hat\delta_q;w,\wellBlockSetInjector)$).
Although scatter is apparent, the relationship is essentially linear. However,
note that $\bar \errorFun(\delta_q;j,\bar{\mathcal D})$
is generally greater than 
$\bar \errorFun(\hat\delta_q;j,\bar{\mathcal D})$, which leads to a systematic bias. Thus, applying the EMML-computed quantities $\bar \errorFun(\hat\delta_q;o,\wellBlockSetProducer)$,
$\bar \errorFun(\hat\delta_q;w,\wellBlockSetProducer)$, and
$\bar \errorFun(\hat\delta_q;w,\wellBlockSetInjector)$
 to predict their `true' counterparts will be biased; this is reflected by the
 red line in
 \figRefOne{fig:predictedErrorPlotMultipleCases_Tr30_M1_Classification_RF},
 which corresponds to the prediction if the EMML-computed quantity alone is
 applied for prediction. 
Note that this bias is not trivial to fix within the regression method itself, as the regression was 
constructed for predicting \textit{time-instantaneous errors}, while the
observed bias
is present for \textit{time-
and well-averaged errors}.

As described in Section \ref{sec:EMML_ROMES}, we can address this issue using the 
ROMES method \cite{drohmann2015romes}. This technique applies
Gaussian-process
(GP)
regression to model the (generally unknown) average true error
$
 \bar \errorFun(\delta_q;j,\bar{\mathcal D})
 $
using an error indicator, which (in this case) corresponds to the (computable) average error predicted by EMML
 $
 \bar \errorFun(\hat\delta_q;j,\bar{\mathcal D})
 $.
To construct the GP, we employ 15 additional high-fidelity and POD--TPWL
simulations for parameter instances corresponding to
$\paramTrainROMES\subset\paramTest$, i.e.,
$\lvert \paramTrainROMES \rvert = 15$. These simulations 
provide the ROMES training data
$
	\{ ( 
 \bar \errorFun(\delta_q(\paramsGen);j,\bar{\mathcal D})
,
 \bar \errorFun(\hat\delta_q(\paramsGen);j,\bar{\mathcal D})
) \}_{\paramsGen\in\paramTrainROMES}.
$
We then construct a GP using the DACE package 
\cite{lophaven2002dace}
with 
%$\theta_L=\theta_U=1\times
%10^{10}$, 
a first-order-polynomial mean function, and a Gaussian covariance
function.
\begin{figure}[H]	
 % -- Oil --
	\begin{subfigure}  {0.5\linewidth} \centering
		\includegraphics[width= 0.85\textwidth]	{{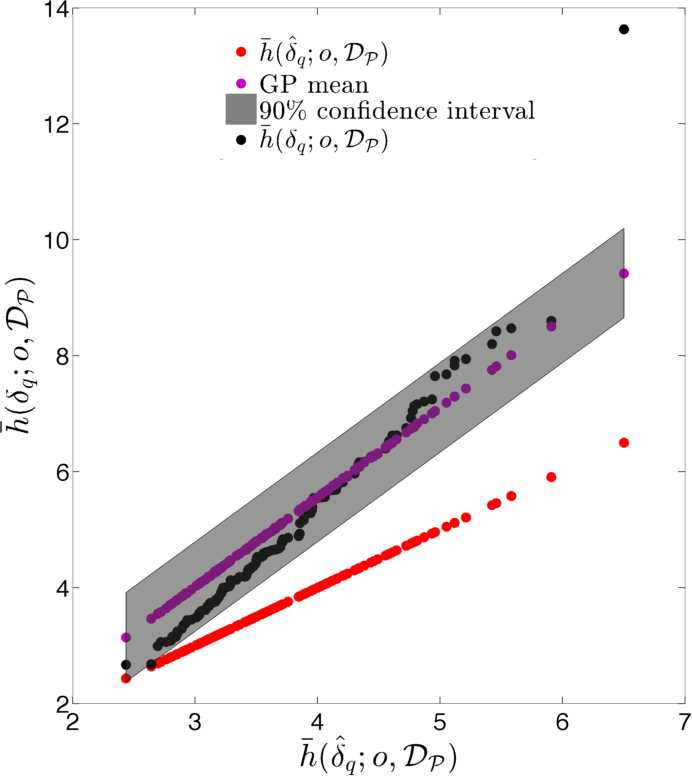}}	       
	        \caption{Oil}
	\end{subfigure}
\quad
	\begin{subfigure}  {0.5\linewidth} \centering
		\includegraphics[width= 0.85\textwidth]	{{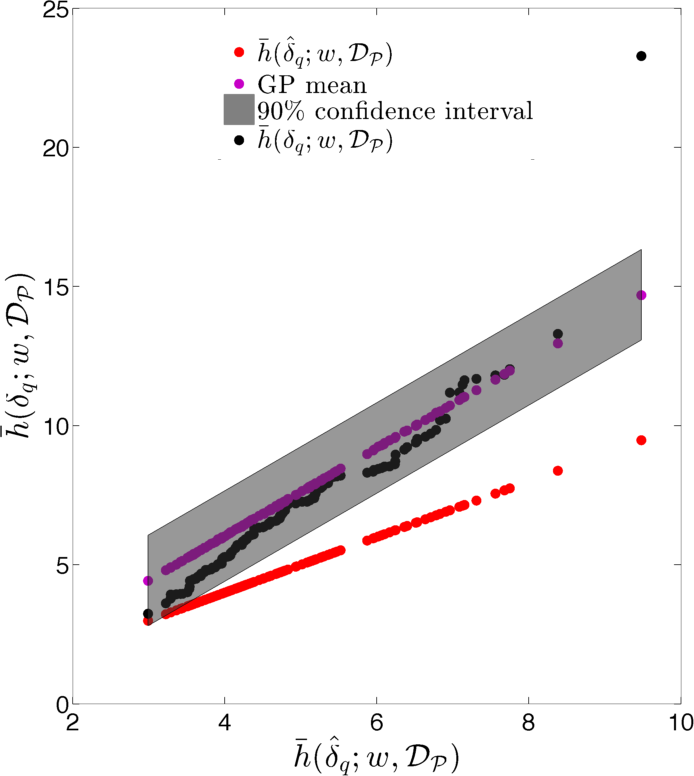}}	       
	        \caption{Water}
	\end{subfigure}
\\
	\begin{subfigure}  {1\linewidth} \centering
		\includegraphics[width= 0.425\textwidth]	{{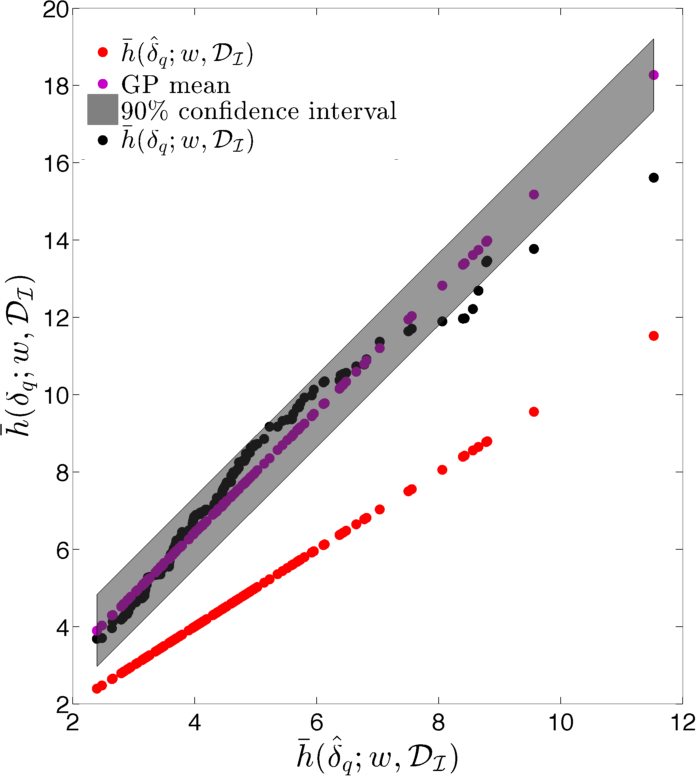}}	       
	        \caption{Water injection rate}
	\end{subfigure}
\caption {Post-corrected time-integrated EMML results using GP.}
\label{fig:ROMESpostProcessing}
\end{figure}

\figRefOne{fig:ROMESpostProcessing} shows the resulting GP. This figure displays the Gaussian
process with a 90\% confidence interval (shaded region), the true
average error (black), the average error predicted by EMML
alone
(red points), and the average error
predicted by EMML after post-processing with ROMES (purple points). Most importantly, these results show that the (mean) EMML prediction after ROMES
postprocessing is significantly more accurate (i.e., closer to
the true errors) than the EMML prediction alone.
Further, the true time-integrated error for the majority of the
test cases lies within the 90\% confidence interval predicted by the GP; this
demonstrates the importance of a \textit{statistical} prediction rather
than a deterministic prediction, as the confidence interval quantifies the prediction
uncertainty. We thus conclude that the method proposed in Section
\ref{sec:EMML_ROMES} is effective at modeling time-integrated errors in this
application.
%
% --- Post corrected time integrated error ---
%

%------------------------------------
\subsection{Computational costs}
\label{sec:costs}
%------------------------------------

We first discuss the computational costs incurred by POD--TPWL and EMML. 
All reported timings were obtained on a machine with dual E5520 Intel CPUs (4
cores, 2.26~Ghz) and 24~GB memory using a Matlab implementation of the high-fidelity and
surrogate models and an $R$~\cite{team2013r} implementation of
EMML.
The offline computational costs for POD--TPWL entail (1) executing 
$|\paramTrainTPWL|=3$ high-fidelity simulations---which can be done in parallel---for
$\paramsGen\in\paramTrainTPWL$;
the cost of a single high-fidelity simulation is 370~seconds, and (2) assembling POD--TPWL operators
via \eqnRefOne{eqn:PODTPWLcoeff}, which consumes 23.5~seconds.

The offline computational costs for EMML training entail (1) executing
$|\paramTrain|=\nTrain=30$ high-fidelity and POD--TPWL simulations (in parallel) for
$\paramsGen\in\paramTrain$;
the POD--TPWL solution takes only 0.5 seconds, which constitutes a speedup of
approximately 700 relative to the HFM, (2) constructing a
classification model for each of the $|\wellBlockSetProducer|=3$ producer
wells $d\in\wellBlockSetProducer$; this consumes 1010~seconds
per producer well, and (3) constructing a random-forest regression model for
each of the 9
QoI defined in \eqnRefOne{eqn:PeacemanModel}; this consumes 412~seconds per
regression model.
The total offline cost of steps 2 and 3---assuming serial
computation---is
6735~seconds. We note that this is approximately 18 times costlier than a
single high-fidelity simulation. Elapsed timings can be readily reduced through use of parallel processing (each QoI can be treated by a different processor). The ratio of offline costs relative to HFM simulation cost will be smaller for larger-dimensional and more complex HFMs, as the EMML training cost in steps 2
and 3 is independent of the complexity of the HFM. However, the EMML training cost does
scale linearly with the number of QoI, assuming serial processing.

In terms of online costs, the POD--TPWL simulation consumes only $0.5$~seconds
(as mentioned above), and the online EMML error prediction
takes about $3\times 10^{-4}$~seconds for
querying the regression model. Thus, online costs are very small relative to
offline costs. We note finally that production optimization computations in
this setting may require $O(100-1000)$ flow simulations, so an offline cost of
$O(10)$ HFM simulations, as required by the EMML-based framework, represents
an acceptable overhead.

%----------------------------------------------
\section{Concluding remarks}
\label{section:conclusion}
%----------------------------------------------

In this work we introduced a general method for error \reviewer{modeling} using machine learning
(EMML). We applied the EMML framework for \reviewer{modeling} the error introduced by 
surrogate models of dynamical systems.
The framework employs high-dimensional regression methods from machine
learning to map a set of inexpensively computed error indicators (features)
to a prediction of the
(time-dependent) surrogate-model quantity-of-interest error (response).
The method
requires first constructing an EMML training dataset by simulating both the
surrogate model and the high-fidelity model for some instances of the input
parameters.
In particular, we proposed:
 \begin{itemize} 
  \item four different methods for modeling the error (Section \ref{subsection:ErrorModeling}), 
\item two methods (classification and clustering) \reviewer{for
determining the notion of locality that is employed to construct} local regression models
(Section
\ref{subsection:spacePartitioning}),
\item two techniques (random forests and LASSO) for performing regression
(Section \ref{sec:regressionMethods}), and
\item two applications of the resulting error models: as a \textit{correction}
	to the surrogate-model QoI prediction, and as a way to model functions of
	the QoI error using the ROMES method (Section \ref{subsection:application}).
	 \end{itemize}

We specialized the method to one particular application: subsurface-flow
modeling with a POD--TPWL reduced-order model as a surrogate. For this application we proposed
specific EMML method ingredients, such as particular choices for error modeling (Section \ref{subsection:Features}), feature design (Section
\ref{sec:feature_design}), training and test data (Section
\ref{subsection:TrainingTestDataset:EMML}), and classification features to use
\reviewer{for determining regression-model locality} (Section
\ref{subsection:spacePartitioning:TPWL}).

In the numerical experiments, we observed that the EMML method performed the best using
the following algorithmic parameters: error-modeling using Approach~1,
memory $\howMuchMemory{1}$, $\numHFS{30}$ simulations to train the EMML
model, classification to \reviewer{determine regression-model locality}, and
random-forest regression. When the EMML error models are used as a correction
to the surrogate-model prediction, we demonstrated improved
accuracy in the output quantities of interest relative to the original
POD--TPWL surrogate model for a large number of test cases (Sections
\ref{sec:EMMLtestCase1}--\ref{subsection:additionalTestCase}). 
When the EMML error models are used to model functions of the QoI error via
ROMES, we observed that the EMML prediction---when combined with a ROMES-based
Gaussian-process model---produced an \reviewer{accurate} prediction with statistical confidence intervals.
It is important to note, however, that the EMML
offline cost is not negligible (Section
\ref{sec:costs}), as it entails (1) constructing the EMML
training dataset, which requires executing $\numHFSsymbol$ high-fidelity and
surrogate-model simulations, (2) \reviewer{determining regression-model locality} for
every QoI, and (3) constructing a regression model for every QoI. Thus, this framework is only appropriate for use in many-query problems such as optimization and uncertainty quantification. 

The EMML framework provides a general error \reviewer{modeling} methodology for
surrogate models of dynamical systems; it assumes only that the surrogate
produces a large set of features that can be mined for potential error
indicators. Thus, it is applicable to a wide variety of surrogates, such as reduced-order models (considered in this work) and those presented in~\cite{carlberg2013gnat, chaturantabut2011application,ghommem2015complexity, peherstorfer2014localized, alotaibi2015global}, and coarsened
models, for example. The application of EMML with upscaled (effectivized) subsurface flow models, within the context of uncertainty quantification, was considered by~\cite{trehanThesis}. Future work should be directed
toward \reviewer{modeling} error with other types of surrogate models (for a range of applications), and also for modeling errors introduced while performing geological parameterization~\cite{vo2014new, vo2016regularized}. Algorithmic improvements within the EMML framework could also be considered. For example, our approach entails univariate regression
for each QoI considered; future work could investigate the use of multivariate 
regression that accounts for interactions
between the QoI.
It may also be worthwhile to explore other regression methods, such as artificial neural networks
with long short-term memory~\cite{hochreiter1997long}. However, our
initial (and ongoing) investigations into the use of deep artificial neural
networks for constructing error models demonstrate that the required training time is very
large relative to that incurred by high-fidelity simulations, so other approaches should also be considered.

\section{Acknowledgments}
We thank the industrial affiliates of the Stanford University Smart Fields and Reservoir Simulation Research (SUPRI-B) Consortia for financial support. 
Sandia National Laboratories is a multi-program laboratory managed and
operated by Sandia Corporation, a wholly owned subsidiary of Lockheed Martin
Corporation, for the U.S. Department of Energy's National Nuclear Security
Administration under contract DE-AC04-94AL85000.

\bibliography{References/reference.bib}        
\appendix
\section{Random-forest regression}
\label{subsection:randomForest}
%-------------------------------------------

\begin{figure}[H]
	\begin{subfigure}  {0.5\linewidth} \centering
		\includegraphics[page= 1, width=0.85 \textwidth]
		{{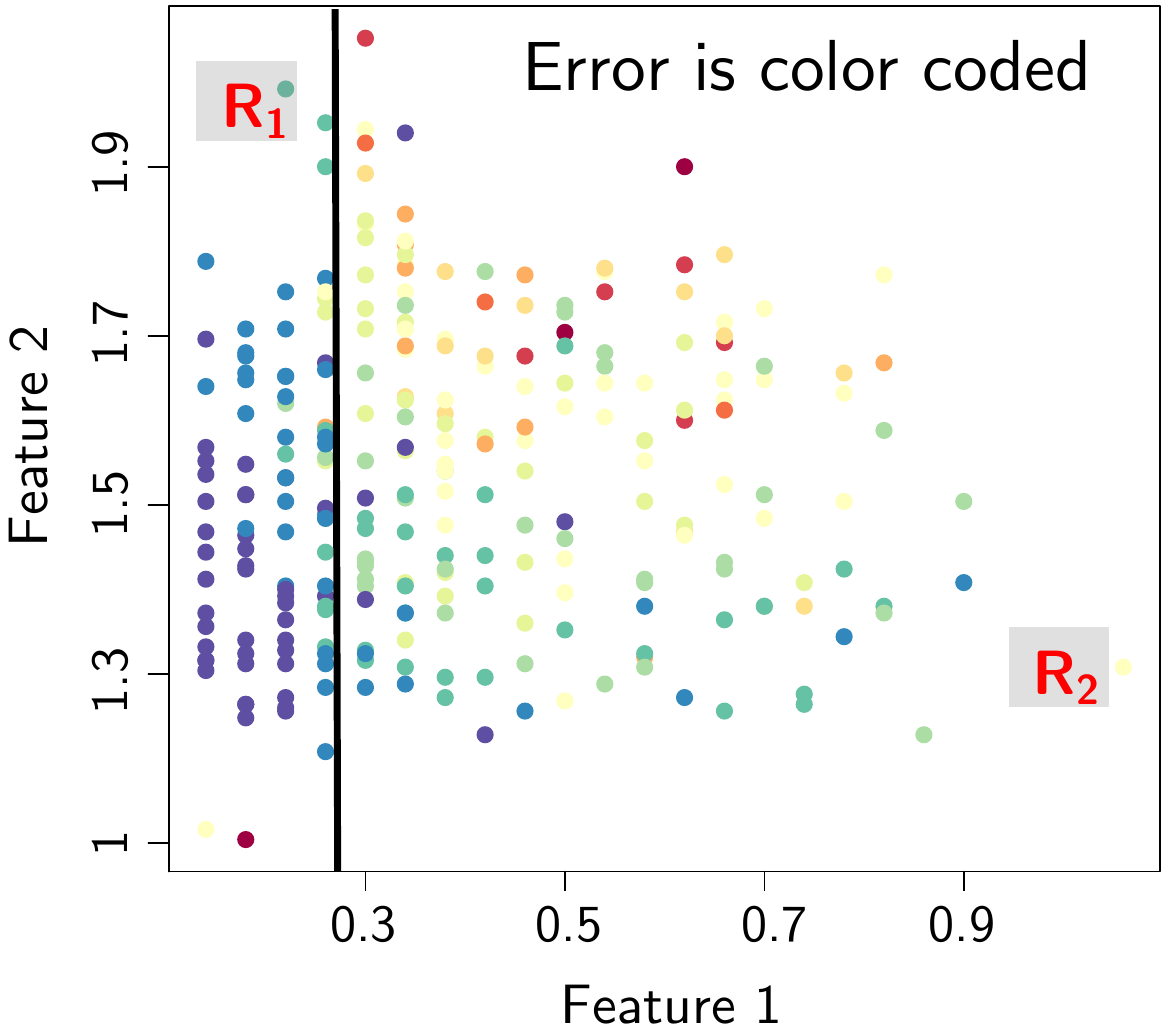}}			\caption{Segmenting domain into
		regions $R_1$ and $R_2$ with $j=1$, $s=0.27$}
	\end{subfigure}
\hspace{1mm}
	\begin{subfigure}  {0.5\linewidth} \centering
		\includegraphics[page= 2, width=0.85 \textwidth]
		{{figures/segmentingDomain.pdf}}			\caption{Segmenting region $R_2$
		into $R_2$ and $R_3$ with $j=2$, $s=1.6$}
	\end{subfigure}
\caption{Schematic of segmentation of the feature space using a decision tree.
Figure adapted and modified from~\cite{james2013introduction}.}
\label{fig:segmentingFeatureSpace}
\end{figure}

\begin{figure}[H]
\centering
	\includegraphics[page= 3, width=0.4 \textwidth]	{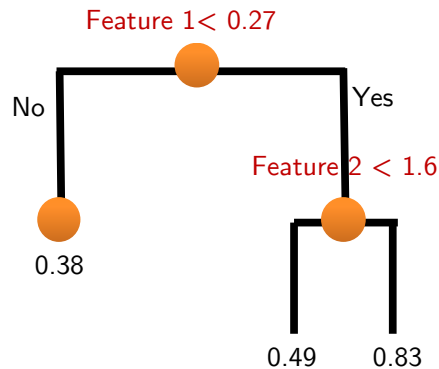}
	\caption{Regression tree.}
\label{fig:segmentingFeatureSpaceDecisionTree}
\end{figure}

Random forest is a supervised machine-learning technique based on constructing
an ensemble of decision trees; it can be used for
both classification and regression. Here, decision trees are constructed by
segmenting the feature space along (canonical) directions corresponding to
individual features. In \figRefOne{fig:segmentingFeatureSpace}, which is adapted
and modified from \cite{james2013introduction}, we plot the (color-coded) true error---which is the response we aim to predict---as a function of two features corresponding to synthetic training
data. 
As shown in \figRefOne{fig:segmentingFeatureSpace}a, segmentation of the
feature space involves splitting the domain into regions $R_1$ and $R_2$. This
is achieved by computing the feature index $j\in\{1,\ldots,\numFeatures\}$
and cutpoint $s\in\RR{}$ corresponding to the segmentation that 
minimizes the residual sum of squares 
\begin{equation}
\label{eqn:randomforest}
	\sum_{n,k: \features^n(\paramsGenArg{k}) \in R_1(j,s)}
	(\errorGen^n(\paramsGenArg{k}) - \hat{\errorGen}_{R_1})^2 + \sum_{n,k:
	\features^n(\paramsGenArg{k}) \in R_2(j,s)} (\errorGen^n(\paramsGenArg{k}) -
	\hat{\errorGen}_{R_2})^2.   
\end{equation}
Here, we define feature-space regions 
$R_1(j,s) = \{\features \ |\ f_j<s\}$ and $R_2(j,s) =
\{\features\ |\ f_j\ge s\}$  and we denote the mean value of the response across the training samples in $R_k$ as
$\hat{\errorGen}_{R_k} \in \RR{},\ k = 1,2$. This
segmentation is carried out recursively, as shown in
\figRefOne{fig:segmentingFeatureSpace}b.
Because recursive segmentation can be summarized in a tree structure as depicted in
\figRefOne{fig:segmentingFeatureSpaceDecisionTree}, this type of approach is
known as a \textit{decision-tree method}. The values assigned to leaves of the
tree correspond to the response predicted in the feature-space region
associated with those leaves, which is the mean value of the response across training data in that region. Note that segmenting the feature space in this
way enables nonlinear interactions among the features to be captured.
Unfortunately, decision-tree models often exhibit high variance and low bias,
which can result in overfitting the training
data~\cite{james2013introduction}. To resolve this problem, random-forest
regression constructs an \textit{ensemble} of decision trees, and the
prediction corresponds to the average prediction across all trees in the
ensemble. While constructing a given decision tree, the
random-forest technique considers only a randomly selected subset of the
feature indices as candidates for performing a split. This randomization serves to
decorrelate the trees, thereby reducing the variance incurred by averaging the
predictions from different trees, which acts to improve prediction quality.

\begin{figure}[H]
	\begin{subfigure}  {0.5\linewidth} \centering
		\includegraphics[page= 1, width=0.65 \textwidth]	{{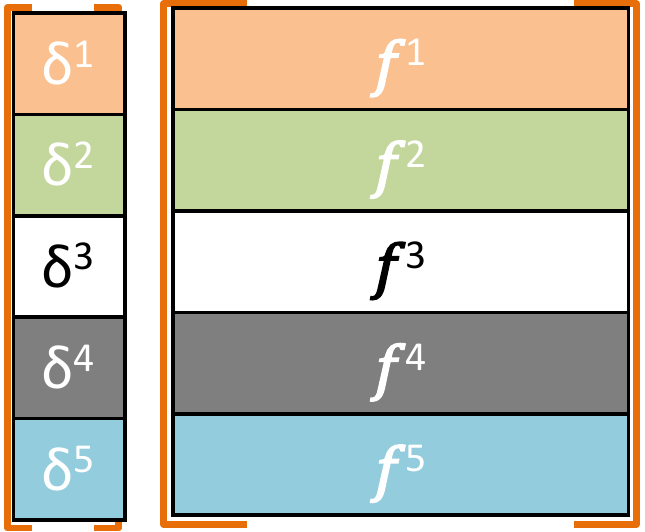}}	
		\caption{Data}
	\end{subfigure}
\quad
	\begin{subfigure}  {0.5\linewidth} \centering
		\includegraphics[page= 2, width=0.65 \textwidth]	{{figures/randomForest.pdf}}	
		\caption{Bootstrapped data used to construct a tree}
	\end{subfigure}
\caption{Bootstrapping of the data.}
\label{fig:bootstrapping}
\end{figure}

Additionally, random-forest regression grows trees on a bootstrap-sampled
version of the training data. Bootstrapping (a resampling technique involving
sampling with replacement) is illustrated in \figRefOne{fig:bootstrapping}a
for a data set containing five samples. To construct the first decision tree,
the random-forest technique samples, with replacement, five training samples
from the data shown in \figRefOne{fig:bootstrapping}a. However, due to
sampling with replacement, it is probable that some samples in the resulting
data set are repeated, as shown in \figRefOne{fig:bootstrapping}b. Bootstrapping further
assists in variance reduction without increasing the bias.

While random forests demonstrate improved accuracy relative to a single decision
tree, they lead to a loss of interpretation. 
\reviewer{We note that the performance of random forests can be improved by various
mechanisms, such as pruning or recursively dropping the least-important features.}
For additional details on random
forests, see Ref.~\cite{breiman2001random}.

%
%-------------------------------------------
\section{LASSO regression}
\label{subsection:lasso}
%-------------------------------------------
Least absolute shrinkage and selection operator (LASSO) is a regression
method that fits a linear model to the feature-to-error mapping, i.e.,
$\errorGen^n(\paramsGen) =
\beta_0 +
\sum\limits_{j=1}^{\numFeatures}f^n_j(\paramsGen)\beta_{j}+\errorMapNoise$,
$n=1,\ldots,\numTimesteps$  and $\beta_j\in\RR{}$,
$j=0,\ldots,\numFeatures$, while reducing
the number of nonzero coefficients. It does so by computing coefficients
$\beta_j$, $j=0,\ldots,\numFeatures$, that minimize an objective
function composed of the sum of
squared errors and an $L_1$-penalty on the coefficients
\begin{align}
	\sum_{n=1}^{\numTimesteps} \sum_{k=1}^{\nTrain}\left(
	\errorGen^{n}(\paramsGenArg{k}) -\beta_0 -
\sum\limits_{j=1}^{\numFeatures}f^n_j(\paramsGenArg{k})\beta_{j} \right)^{2} + \lambda
	\sum_{j=1}^{\numFeatures} \left | \beta_{j} \right |,
\end{align}
where $\lambda\in\RR{}$ is a penalty parameter, and larger values of 
$\lambda$ promote sparsity in the computed coefficients.

The value of $\lambda$ is typically chosen by $k$-fold cross validation, which is a
process that randomly partitions the training data into $k$ subsamples
of equal size. One subsample is withheld while a
linear model is constructed using the remaining $k-1$ subsamples. The model is
then tested on the withheld sample, and the prediction errors are recorded. To
reduce variability, the process is repeated $k$ times, with a different
subsample withheld each time. Finally, the prediction errors are averaged, and
this average is referred to as the cross-validation error. To determine the optimal
value of $\lambda$ in practice, we compute the cross-validation error on only
a subset of candidate values for $\lambda$, and select the value yielding the
smallest cross-validation error.
%
%-------------------------------------------

\end{document}